\documentclass[12pt,twoside, leqno]{article}
\usepackage{amssymb}
\usepackage{amsmath}
\usepackage{amsthm}
\usepackage{amssymb,mathrsfs,amsmath}
\usepackage[perpage,symbol*]{footmisc}

\allowdisplaybreaks

\def\ls{\lesssim}
\def\gs{\gtrsim}
\def\fz{\infty}

\def\r{\right}
\def\lfz{\left}

\def\ls{\lesssim}
\def\gs{\gtrsim}

\def\paz{{\partial}}

\def\supp{{\mathop\mathrm{\,supp\,}}}

\def\rr{{\mathbb R}}
\def\rn{{{\rr}^n}}
\def\zz{{\mathbb Z}}
\def\nn{{\mathbb N}}
\def\cc{{\mathbb C}}

\newcommand{\wz}{\widetilde}

\newcommand{\ca}{{\mathcal A}}

\newcommand{\cd}{{\mathcal D}}

\newcommand{\cg}{{\mathcal G}}

\newcommand{\cn}{{\mathcal N}}

\newcommand{\cs}{{\mathcal S}}

\def\az{\alpha}
\def\lz{\lambda}
\def\blz{\Lambda}

\def\bfai{\Phi}
\def\dz{\delta}

\def\epz{\epsilon}

\def\bz{\beta}
\def\rz{\rho}
\def\pz{\psi}

\def\gz{{\gamma}}
\def\bgz{{\Gamma}}

\def\tz{\theta}
\def\sz{\sigma}
\def\zez{\zeta}
\def\wz{\widetilde}
\def\wt{\widehat}
\def\ls{\lesssim}
\def\gs{\gtrsim}

\def\ol{\overline}

\def\boz{\Omega}

\def\uc{{\varepsilon}}

\def\esup{\mathop\mathrm{\,esssup\,}}

\def\divz{{{\mathop\mathrm {div}}}}

\def\hs{\hspace{0.3cm}}

\newtheorem{theorem}{Theorem}[section]
\newtheorem{lemma}{Lemma}[section]

\newtheorem{proposition}{Proposition}[section]

\newtheorem{remark}{Remark}[section]
\newtheorem{definition}{Definition}[section]

\def\supp{{\mathop\mathrm{\,supp\,}}}
\def\diam{{\mathop\mathrm{diam}}}
\def\dist{{\mathop\mathrm{\,dist\,}}}

\def\lf{\lfloor}
\def\rf{\rfloor}
\numberwithin{equation}{section}

\begin{document}

\arraycolsep=1pt

\title{{\vspace{-5cm}\small\hfill\bf Rev. Mat. Iberoam. (to appear)}\\
\vspace{4.5cm}\Large Real-variable characterizations
of Orlicz-Hardy spaces on strongly Lipschitz domains of
$\mathbb{R}^n$ \footnotetext{\hspace{-0.35cm} 2010 {\it Mathematics
Subject Classification}: {Primary: 42B30; Secondary: 42B35, 42B20,
42B25, 35J25, 42B37, 47B38.}
\endgraf{\it Keywords}: Orlicz-Hardy space,  divergence
form elliptic operator, strongly Lipschitz domain, Dirichlet
boundary condition, Gaussian property, nontangential maximal
function, Lusin area function, atom.}}
\author{\bf Dachun Yang and Sibei Yang}
\date{ }
\maketitle

\begin{center}
\begin{minipage}{4.35in}\small
\begin{center}{\bf Abstract}\end{center}
\quad Let $\Omega$ be a strongly Lipschitz domain of $\mathbb{R}^n$,
whose complement in $\mathbb{R}^n$ is unbounded. Let $L$ be a second
order divergence form elliptic operator on $L^2 (\Omega)$ with the
Dirichlet boundary condition, and the heat semigroup generated by
$L$ have the Gaussian property $(G_{\mathrm{diam}(\Omega)})$ with
the regularity of their kernels measured by $\mu\in(0,1]$, where
$\mathrm{diam}(\Omega)$ denotes the diameter of $\Omega$. Let $\Phi$
be a continuous, strictly increasing, subadditive and positive
function on $(0,\infty)$ of upper type 1 and of strictly critical
lower type $p_{\Phi}\in(n/(n+\mu),1]$. In this paper, the authors
introduce the Orlicz-Hardy space $H_{\Phi,\,r}(\Omega)$ by
restricting arbitrary elements of the Orlicz-Hardy space
$H_{\Phi}(\mathbb{R}^n)$ to $\boz$ and establish its atomic
decomposition by means of the Lusin area function associated with
$\{e^{-tL}\}_{t\ge0}$. Applying this, the authors obtain two
equivalent characterizations of $H_{\Phi,\,r}(\boz)$ in terms of the
nontangential maximal function and the Lusin area function
associated with the heat semigroup generated by $L$.
\end{minipage}
\end{center}

\section{\hspace{-0.7cm}{\bf.}  Introduction\label{s1}}

The theory of Hardy spaces on the $n$-dimensional Euclidean space
$\rn$, was originally initiated by Stein and Weiss in \cite{sw60}.
Later, Fefferman and Stein \cite{fs72} systematically developed a
real-variable theory for the Hardy spaces $H^p (\rn)$ with
$p\in(0,1]$, which plays an important role in various fields of
analysis; see, for example, \cite{s93,clms,m94,s94}. It is well
known that the Hardy space $H^p (\rn)$ with $p\in(0,1]$ is a good
substitute of $L^p (\rn)$ in the study of the boundedness of
operators; for example, the classical Riesz transform is bounded on $H^p
(\rn)$, but not on $L^p (\rn)$ with $p\in(0,1]$. An important
feature of $H^p (\rn)$ is their atomic decomposition
characterizations, which were established by Coifman \cite{co74}
when $n=1$ and Latter \cite{l78} when $n>1;$ see also \cite{w85}.

On the other hand,  as a generalization of $L^p (\rn)$, 
the Orlicz space was introduced by
Birnbaum-Orlicz in \cite{bo31} and Orlicz in \cite{o32}; since then,
the theory of the Orlicz spaces themselves has been well developed
and these spaces have been widely used in probability, statistics,
potential theory, partial differential equations, as well as
harmonic analysis and some other fields of analysis; see, for
example, \cite {rr91,rr00,byz08,mw08,io02}. Moreover, Orlicz-Hardy
spaces are also suitable substitutes of the Orlicz spaces in the
study of boundedness of operators; see, for example,
\cite{ja80,vi87,jy10,jyz09,jy}. Recall that Orlicz-Hardy spaces and
their dual spaces were studied by Janson \cite{ja80} on $\rn$ and
Viviani \cite{vi87} on spaces of homogeneous type in the sense of
Coifman and Weiss \cite{cw71}.

It is known that Hardy spaces $H^p (\rn)$ are essentially related to
the Laplacian 
$$\Delta:=\sum^n_{i=1}\frac {\partial^2} {\partial
x_i^2}.$$ 
In recent years, the study of the real-variable theory of
various function spaces associated with different differential
operators has inspired great interests; see, for example,
\cite{adm,amr08,dy05,y08,dxy,hm09,hlmmy,jy,jy10,jyz09,yyj}.
In particular, Orlicz-Hardy spaces associated with some differential
operators and their dual spaces were introduced and studied in
\cite{jyz09,jy10,jy}.

One important aspect of the development in the theory of Hardy
spaces is the study of Hardy spaces on domains of $\rn$; see, for
example, \cite{m90,cks93,cds99,tw96,ar03,dy04,h,h09}. Especially,
Chang, Krantz and Stein \cite{cks93} introduced the Hardy spaces
$H^p_r (\boz)$ and $H^p_z (\boz)$ on the domain $\boz$ for
$p\in(0,1]$, respectively, by restricting arbitrary elements of $H^p
(\rn)$ to $\boz$, and restricting elements of $H^p (\rn)$ which are
zero outside $\ol{\boz}$ to $\boz$, where and in what follows,
$\ol{\boz}$ denotes the {\it closure} of $\boz$ in $\rn$.
We point out that the Hardy spaces $H^p_r (\boz)$ and $H^p_z (\boz)$,
when $\boz$ is a bounded smooth domain of $\rn$ and
$p\in(0,1]$, naturally appeared in the study of
the regularity of the Green operators, respectively,
for the Dirichlet boundary problem and the Neumann
boundary problem in \cite{cks93,cds99}.
For these Hardy spaces, atomic decompositions have been
obtained in \cite{cks93} when $\boz$ is a special Lipschitz domain
or a bounded Lipschitz domain of $\rn$. Let $\boz$ be a strongly
Lipschitz domain, $H^1_r (\boz)$ and $H^1_z (\boz)$ be defined as in
\cite{cks93}. Auscher and Russ \cite{ar03} proved that $H^1_r(\boz)$
and $H^1_z (\boz)$ can be characterized by the non-tangential
maximal function and the Lusin area function associated with
$\{e^{-t\sqrt L}\}_{t\ge0}$, respectively, under the so-called
Dirichlet and the Neumann boundary conditions, where $L$ is an
elliptic second-order divergence operator such that for all $t\in(0,
\fz)$, the kernel of $e^{-tL}$ has the Gaussian property
$(G_{\infty})$ in the sense of Auscher and Russ \cite[Definition
3]{ar03} (see also Definition \ref{d2.1} below). Moreover, for these
Hardy spaces, Huang \cite{h} established a characterization in terms
of the Littlewood-Paley-Stein function associated with $L$. Assume
that the regularity of the kernel of the heat semigroup generated by
$L$ is measured by $\mu\in(0,1]$. When $\boz$ is a special Lipschitz
domain of $\rn$, $p\in(n/(n+\mu),1]$ and $L$ satisfies the Neumann
boundary condition, Duong and Yan \cite{dy04} gave a simple proof of
the atomic decomposition for elements in $H^p_z (\boz)$ via the
nontangential maximal function associated with the Poisson semigroup
generated by $L$.

Let $\boz$ be a strongly Lipschitz domain of $\rn$, whose complement
in $\rn$ is unbounded. Let $L$ be a second order divergence form
elliptic operator on $L^2 (\boz)$ with the Dirichlet boundary
condition, and the heat semigroup generated by $L$ have the Gaussian
property $(G_{\mathrm{diam}(\boz)})$ with the regularity of their
kernels measured by $\mu\in(0,1]$ (see Definition \ref{d2.1} below
for the definition), where $\mathrm{diam}(\Omega)$ denotes the {\it
diameter} of $\boz$. Let $\Phi$ be a continuous, strictly
increasing, subadditive and positive function on $(0,\infty)$ of
upper type 1 and of strictly critical lower type
$p_{\Phi}\in(n/(n+\mu),1]$ (see \eqref{2.4} below for the definition
of $p_{\Phi}$).  A typical example of such functions is 
$$\Phi(t):=t^p$$ 
for all $t\in(0,\fz)$ and $p\in (n/(n+\mu),1]$. Motivated by
\cite{ar03,cks93,jyz09,jy10,vi87}, in this paper, we introduce the
Orlicz-Hardy space $H_{\Phi,\,r}(\Omega)$ by restricting elements of
the classical Orlicz-Hardy space $H_{\Phi}(\mathbb{R}^n)$ to $\boz$,
and give its atomic decomposition by means of the Lusin area
function associated with the heat semigroup generated by $L$.
Applying this, we obtain two equivalent characterizations of
$H_{\Phi,\,r}(\boz)$ in terms of the nontangential maximal function
and the Lusin area function associated with the heat semigroup
generated by $L$. Let $H^1_{S_P} (\boz)$ be the Hardy space defined
by the Lusin area function associated with the Poisson semigroup
generated by $L$. As a byproduct, by applying the method used in
this paper for the atomic decomposition of elements in
$H_{\Phi,\,r}(\Omega)$ via the Lusin area function associated with
the heat semigroup generated by $L$ (see Proposition \ref{p3.4}
below), we also give a direct proof of the atomic decomposition for
all $f\in H^1_{S_P}(\boz)$ in Proposition \ref{p3.5} below, which
answers the question asked by Duong and Yan \cite[p.\,485, Remarks
(iii)]{dy04} in the case that $p=1$.

To state the main result of this paper, we first recall some necessary notions.
Throughout the whole paper, we always assume that $\boz$ is a {\it
strongly Lipschitz domain} of $\rn$; namely, $\boz$ is a proper open
connected set in $\rn$ whose boundary is a finite union of parts of
rotated graphs of Lipschitz maps, at most one of these parts
possibly unbounded. It is well known that strongly Lipschitz domains
include special Lipschitz domains, bounded Lipschitz domains and
exterior domains; see, for example, \cite{ar03,at01a} for their
definitions and properties.

Throughout the whole paper, for the sake of convenience, we choose
the norm on $\rn$ to be the {\it supremum norm}; namely, for any
$$x=(x_1,\,x_2,\,\cdots,\,x_n)\in\rn,\ \
|x|:=\max\{|x_1|,\,\cdots,\,|x_n|\},$$ for which balls determined by
this norm are cubes associated with the usual Euclidean norm with
sides parallel to the axes.

\begin{remark}\rm\label{r1.1}
Let $\boz$ be a strongly Lipschitz domain of $\rn$. Then $\boz$ is a
space of homogeneous type in the sense of Coifman and Weiss
\cite{cw71}. Furthermore, as a space of homogeneous type, the
collection of all balls of $\boz$ is given by the set
$$\lfz\{Q\cap\boz:\ \text{cube}\
Q\subset\rn\ \text{satisfying}\  x_{Q}\in\boz\  \text{and}\
l(Q)\le2\diam(\boz)\r\},$$ where $x_Q$ denotes the {\it center} of
$Q$, $l(Q)$ the {\it sidelength} of $Q$ and $\diam(\boz)$ the {\it
diameter} of $\Omega$, namely,
$$\diam(\boz):=\sup\{|x-y|:\
x,\,y\in\boz\};$$
see, for example, \cite{ar03}.
\end{remark}

Motivated by \cite{cks93}, we introduce the Orlicz-Hardy space
$H_{\Phi,\,r}(\Omega)$ as follows. We first recall the definition of the
Orlicz-Hardy space $H_{\Phi}(\rn)$ introduced by Viviani
\cite{vi87}. Let $\cs(\rn)$ denote the {\it space of all Schwartz
functions} with the classical topology and $\cs' (\rn)$ its {\it
topological dual} with the weak $\ast$-topology. For all
$f\in\cs' (\rn)$, let $\cg (f)$ denote its {\it grand maximal function};
see \cite[p.\,90]{s93}.

\begin{definition}\rm\label{d1.1}
Let $\Phi$ be a function of type $(p_0,p_1)$, where $0<p_0\le
p_1\le1$ (see Section \ref{s2.2} below for the definition of type
$(p_0,p_1)$). Define
$$H_{\Phi}(\rn):=\lfz\{f\in\cs' (\rn):\ \int_{\rn}
\Phi(\cg (f)(x))\,dx<\fz\r\}
$$
and
$$\|f\|_{H_{\Phi}(\rn)}:=\inf\lfz\{\lz\in (0,\fz):\ \int_{\rn}
\Phi\lfz(\frac{\cg (f)(x)} {\lz}\r)\,dx\le1\r\}.$$
\end{definition}

In what follows, let $\cd(\boz)$ denote the {\it space of all
infinitely differentiable functions with compact support in $\boz$}
endowed with the inductive topology, and $\cd'(\boz)$ its {\it
topological dual} with the weak $\ast$-topology which is called the
{\it space of distributions on $\boz$}.

\begin{definition}\rm\label{d1.2}
Let $\Phi$ be as in Definition \ref{d1.1} and $\boz$ a subdomain in
$\rn$. A distribution $f$ on $\boz$ is said to be in the {\it
Orlicz-Hardy space} $H_{\Phi,\,r}(\boz)$ if $f$ is the restriction
to $\boz$ of a distribution $F$ in $H_{\Phi}(\rn)$; namely,
\begin{eqnarray*}
H_{\Phi,\,r}(\boz):&=& \{f\in \cd'(\boz):\ \text{there exists an }\
F\in H_{\Phi}(\rn)\ \text{such that}\ F|_{\boz}=f\}\\
&=&H_{\Phi}(\rn)/\{F\in H_{\Phi}(\rn):\ F=0\ \text{on}\ \boz\}.
\end{eqnarray*}
Moreover, for all $f\in H_{\Phi,\,r}(\boz)$, the {\it quasi-norm} of
$f$ in $H_{\Phi,\,r}(\boz)$ is defined by
$$\|f\|_{H_{\Phi,\,r}(\boz)}:=\inf\lfz\{\|F\|_{H_{\Phi}(\rn)}:\
F\in H_{\bfai}(\rn)\ \text{and}\ F|_{\boz}=f\r\},$$
where the infimum is taken over all $F\in H_{\Phi}(\rn)$ satisfying
$F=f$ on $\boz$.
\end{definition}

\begin{remark}\rm\label{r1.2}
 Let $p\in(0,1]$. When $\Phi(t):= t^p$ for all $t\in(0,\fz)$,
the space $H_{\Phi,\,r}(\boz)$ was introduced by Chang, Krantz and
Stein \cite{cks93}. In this case, we denote the Orlicz-Hardy
spaces $H_{\Phi}(\rn)$ and $H_{\Phi,\,r}(\boz)$, respectively,
by $H^p (\rn)$ and $H^p_r(\boz)$.
\end{remark}

We now describe the divergence form elliptic operators considered in
this paper and the most typical example is the Laplace operator on
the Lipschitz domain of $\rn$ with the Dirichlet boundary condition. If
$\boz$ is a strongly Lipschitz domain of $\rn$, we denote by
$W^{1,\,2}(\boz)$ the usual {\it Sobolev space on $\boz$} equipped
with the norm
$$\lfz(\|f\|^2_{L^2 (\boz)}+\|\nabla f\|^2_{L^2
(\boz)}\r)^{1/2},$$
where $\nabla f$ denotes the {\it distributional
gradient} of $f$. In what follows, $W^{1,\,2}_0 (\boz)$ stands for
the {\it closure of $C^{\fz}_c (\boz)$ in $W^{1,\,2}(\boz)$}, where
$C^{\fz}_c(\boz)$ denotes the {\it set of all $C^\fz(\rn)$
functions on $\boz$ with compact support}.

If $A:\ \rn\to M_n (\cc)$ is a measurable function, define
$$\|A\|_\fz:=\esup_{x\in\rn,\,|\xi|=|\eta|=1}|A(x)\xi
\cdot\ol{\eta}|,$$
where $M_n (\cc)$ denotes the
{\it set of all $n\times
n$ complex-valued matrixes}, $\xi,\,\eta\in\cc^n$ and $\ol{\eta}$
denotes the {\it conjugate vector} of $\eta$. For all $\dz\in(0,1]$,
denote by $\ca (\dz)$ the {\it class of all measurable
functions $A:\
\rn\to M_n (\cc)$ satisfying the ellipticity condition}; namely, for
all $x\in\rn$ and $\xi\in\cc^n$,
\begin{equation}\label{1.1}
\|A\|_\fz\le \dz^{-1} \ \text{and}\ \Re
(A(x)\xi\cdot\xi)\ge\dz|\xi|^2,
\end{equation}
where and in what follows, $\Re (A(x)\xi\cdot\xi)$ denotes the {\it real
part} of $ A(x)\xi\cdot\xi$. Denote by $\ca$
the {\it union of all $\ca (\dz)$ for $\dz\in(0,1]$}.

When $A\in\ca$ and $V$ is a {\it closed subspace of
$W^{1,\,2}(\boz)$ containing $W^{1,\,2}_0 (\boz)$}, denote by $L$ the
{\it maximal-accretive operator} (see \cite[p.\,23, Definition
1.46]{o04} for the definition) on $L^2 (\boz)$ with largest domain
$D(L)\subset V$ such that for all $f\in D(L)$ and $g\in V$,
\begin{equation}\label{1.2}
\langle Lf,g\rangle=\int_{\boz}A(x)\nabla f(x)\cdot\ol{\nabla
g(x)}\,dx,
\end{equation}
where $\langle \cdot,\cdot\rangle$ denotes the {\it interior product} in
$L^2 (\boz)$. In this sense, for all $f\in D(L)$, we write
\begin{equation}\label{1.3}
Lf:=-\divz(A\nabla f).
\end{equation}

We recall the following Dirichlet and Neumann boundary conditions
of $L$ from \cite[p.\,152]{ar03}.

\begin{definition}\rm\label{d1.3}
Let $\boz$ be a strongly Lipschitz domain and $L$ as in \eqref{1.3}.
The operator $L$ is called to satisfy the {\it Dirichlet boundary
condition} (for simplicity, DBC) if $V:= W^{1,\,2}_0 (\boz)$ and the
{\it Neumann boundary condition} (for simplicity, NBC) if $V:=
W^{1,\,2}(\boz)$.
\end{definition}

Let $\boz$ be a strongly Lipschitz domain
of $\rn$. Recall that for an Orlicz function $\Phi$ on $(0,\fz)$, a
measurable function $f$ on $\boz$ is called to be in the {\it space
$L^{\Phi}(\boz)$} if $\int_{\boz}\Phi(|f(x)|)\,dx<\fz$. Moreover,
for any $f\in L^{\Phi}(\boz)$, define
$$\|f\|_{L^{\Phi}(\boz)}:=\inf\lfz\{\lz\in(0,\fz):\
\int_{\boz}\Phi\lfz(\frac{|f(x)|}{\lz}\r)\,dx\le1\r\}.$$ If
$p\in(0,1]$ and $\bfai(t):= t^p$ for all $t\in(0,\fz)$,  we then
denote $L^{\bfai}(\boz)$ simply by $L^p (\boz)$.

\begin{definition}\rm\label{d1.4}
Let $\Phi$ satisfy Assumption (A) (see Section \ref{s2.2} for the
definition of Assumption (A)), $\boz$ be a strongly Lipschitz domain
of $\rn$ and $L$ as in \eqref{1.3}.
For all $f\in L^2 (\boz)$ and $x\in\boz$, let
$$\cn_h (f)(x):=\sup_{y\in\boz,\,t\in(0,2\diam(\boz)),\,|y-x|<t}
\lfz|e^{-t^2 L}(f)(y)\r|.$$ A function $f\in L^2 (\boz)$ is said to be
in $\wz{H}_{\Phi,\,\cn_h}(\boz)$ if $\cn_h (f)\in L^{\Phi}(\boz)$;
moreover, define
\begin{eqnarray*}
\|f\|_{H_{\Phi,\,\cn_h}(\boz)}&&:=\lfz\|\cn_h
(f)\r\|_{L^{\Phi}(\boz)}\\
&&:=\inf\lfz\{\lz\in(0,\fz):\ \int_{\boz}\Phi\lfz(\frac{\cn_h
(f)(x)}{\lz}\r)\,dx\le1\r\}.
\end{eqnarray*}
The {\it Orlicz-Hardy space} $H_{\Phi,\,\cn_h}(\boz)$ is defined to
be the completion of the space $\wz{H}_{\Phi,\,\cn_h}(\boz)$ in the
quasi-norm $\|\cdot\|_{H_{\Phi,\,\cn_h}(\boz)}$.
\end{definition}

\begin{remark}\rm\label{r1.3}
(i) Since $\Phi$ is of strictly lower type $p_{\Phi}$ (see
\eqref{2.4} for its definition), we have that for all $f_1,\,f_2\in
H_{\Phi,\,\cn_h}(\boz)$,
$$\|f_1 +f_2\|^{p_{\Phi}}_{H_{\Phi,\,\cn_h}(\boz)}\le\|f_1 \|^{p_{\Phi}}
_{H_{\Phi,\,\cn_h}(\boz)}+\|f_2\|^{p_{\Phi}}_{H_{\Phi,\,\cn_h}(\boz)}.$$

(ii) From the theorem of completion of Yosida \cite[p.\,56]{yo95},
it follows that $\wz{H}_{\Phi,\,\cn_h}(\boz)$ is dense in
$H_{\Phi,\,\cn_h}(\boz)$; namely, for any $f\in
H_{\Phi,\,\cn_h}(\boz)$, there exists a Cauchy sequence
$\{f_k\}_{k=1}^{\fz} \subset \wz{H}_{\Phi,\,\cn_h}(\boz)$ such that
$$\lim_{k\to\fz}\|f_k -f\| _{H_{\Phi,\,\cn_h}(\boz)}=0.$$
Moreover, if $\{f_k\}_{k=1}^{\fz}$ is a Cauchy sequence in
$\wz{H}_{\Phi,\,\cn_h}(\boz)$, then there  exists a uniquely $f\in
H_{\Phi,\,\cn_h}(\boz)$ such that
$$\lim_{k\to\fz}\|f_k -f\|_{H_{\Phi,\,\cn_h}(\boz)}=0.$$
\end{remark}

In what follows, $Q(x,t)$ denotes the {\it closed cube of $\rn$
centered at $x$ and of the sidelength $t$ with sides parallel to the
axes}. Similarly, given $Q:= Q(x,t)$ and $\lz\in(0,\fz)$, we write
$\lz Q$ for the {\it$\lz$-dilated cube}, which is the cube with the
same center $x$ and with sidelength $\lz t$. For any $f\in L^2
(\boz)$ and $x\in\boz$, the {\it Lusin area functions} $S_h$ and
$\wz{S}_h$ associated with $\{e^{-t^2 L}\}_{t\ge0}$ are respectively
defined by
\begin{equation*}
S_h (f)(x):=\lfz\{\int_{\bgz(x)}\lfz|t^2 Le^{-t^2 L}
(f)(y)\r|^2\,\frac{dy\,dt}{t|Q(x,t)\cap\boz|}\r\}^{1/2}
\end{equation*}
and
\begin{equation*}
\wz{S}_h (f)(x):=\lfz\{\int_{\bgz(x)}\lfz|t \nabla e^{-t^2 L}
(f)(y)\r|^2\, \frac{dy\,dt}{t|Q(x,t)\cap\boz|}\r\}^{1/2},
\end{equation*}
 where $\bgz(x)$ is the {\it cone} defined by
$$\bgz(x):=\{(y,t)\in\boz\times(0,2\diam(\boz)):\ |y-x|<t\}.$$

\begin{definition}\rm\label{d1.5}
Let $\Phi$ satisfy Assumption (A), $\boz$ be a strongly Lipschitz
domain of $\rn$ and $L$ as in \eqref{1.3}.  Assume that $L$
satisfies DBC and the semigroup generated by $L$ has the Gaussian
property $(G_{\diam(\boz)})$. A function $f\in L^2 (\boz)$ is said
to be in $\wz{H}_{\Phi,\,S_h}(\boz)$ if $S_h (f)\in L^{\Phi}(\boz)$.
Recall that
\begin{equation}\label{1.4}
\|S_h (f)\|_{L^{\Phi}(\boz)}:=\inf\lfz\{\lz\in(0,\fz):\
\int_{\boz}\Phi\lfz(\frac{S_h (f)(x)}{\lz}\r)\,dx\le1\r\}.
\end{equation}
Furthermore, define
$$\|f\|_{H_{\Phi,\,S_h}(\boz)}:=\|S_h
(f)\|_{L^{\Phi}(\boz)}.$$
The {\it Orlicz-Hardy space}
$H_{\Phi,\,S_h}(\boz)$ is defined to be the completion of
$\wz{H}_{\Phi,\,S_h}(\boz)$ in the quasi-norm
$\|\cdot\|_{H_{\Phi,\,S_h}(\boz)}$.

If $\boz$ is bounded, a function $f\in L^2 (\boz)$ is said to be in
$\wz{H}_{\Phi,\,S_h,\,d_{\boz}}(\boz)$ if $S_h (f)\in
L^{\Phi}(\boz)$; moreover, define
\begin{eqnarray}\label{1.5}
\quad&&\|f\|_{H_{\Phi,\,S_h,\,d_{\boz}}(\boz)}\\
&&\hs:=\!\|S_h(f)\|_{L^{\Phi}(\boz)}+\inf\lfz\{\lz\in(0,\fz):\
\bfai\lfz(\frac{\|e^{-d_{\boz}^2 L}(f)\|_{L^1
(\boz)}}{\lz}\r)\le1\r\},\nonumber
\end{eqnarray}
where and in what follows, $d_{\boz}:=2\diam(\boz)$ and $\|S_h
(f)\|_{L^{\Phi}(\boz)}$ is as in \eqref{1.4}. The {\it Orlicz-Hardy
space} $H_{\Phi,\,S_h,\,d_{\boz}}(\boz)$ is defined to be the
completion of the space $\wz{H}_{\Phi,\,S_h,\,d_{\boz}}(\boz)$ in
the quasi-norm $\|\cdot\|_{H_{\Phi,\,S_h,\,d_{\boz}}(\boz)}$.

The {\it Orlicz-Hardy
spaces} $H_{\Phi,\,\wz{S}_h}(\boz)$ and
$H_{\Phi,\,\wz{S}_h,\,d_{\boz}}(\boz)$ when $\boz$ is bounded are
defined via replacing $S_h$, respectively, in the definitions of
$H_{\Phi,\,S_h}(\boz)$ and $H_{\Phi,\,S_h,\,d_{\boz}}(\boz)$ by
$\wz{S}_h$.
\end{definition}

If $\boz$ is bounded, by $|\boz|<\fz$, we know that $L^2
(\boz)\subset L^1 (\boz)$, which, together with the Gaussian
property $(G_{\diam(\boz)})$ and Fubini's theorem, implies that for
all $f\in L^2 (\boz)$, $e^{-d_{\boz}^2 L}(f)\in L^1 (\boz)$. Thus,
if $f\in L^2 (\boz)$ and $S_h (f)\in L^{\Phi}(\boz)$,  then
$\|f\|_{H_{\Phi,\,S_h,\,d_{\boz}}(\boz)}$ and
$\|f\|_{H_{\Phi,\,\wz{S}_h,\,d_{\boz}}(\boz)}$ make sense.

In what follows, we denote by $\boz^\complement$ the {\it complement
of $\boz$ in $\rn$}. The main result of this paper is as follows.

\begin{theorem}\label{t1.1}
Let $\Phi$ satisfy Assumption (A) and $L$ be as in \eqref{1.3}. Let
$\boz$ be a strongly Lipschitz domain of $\rn$ such that
$\boz^\complement$ is unbounded. Assume that $L$ satisfies DBC and
the semigroup generated by $L$ has the Gaussian property
$(G_{\diam(\boz)})$.

\begin{enumerate}
  \item[$\mathrm{(i)}$] If $\boz$ is unbounded, then the spaces
$H_{\Phi,\,r}(\boz)$, $H_{\Phi,\,\cn_h}(\boz)$,
$H_{\Phi,\,\wz{S}_h}(\boz)$ and $H_{\Phi,\,S_h}(\boz)$
coincide with equivalent norms.

  \item[$\mathrm{(ii)}$] If $\boz$ is bounded, then the spaces
$H_{\Phi,\,r}(\boz)$, $H_{\Phi,\,\cn_h}(\boz)$,
$H_{\Phi,\,\wz{S}_h,\,d_{\boz}}(\boz)$ and
$H_{\Phi,\,S_h,\,d_{\boz}}(\boz)$
coincide with equivalent norms. Moreover, if, in addition,  $n\ge3$ and
$(G_{\fz})$ holds, then the spaces
$H_{\Phi,\,\wz{S}_h,\,d_{\boz}}(\boz)$,
$H_{\Phi,\,S_h,\,d_{\boz}}(\boz)$, $H_{\Phi,\,\wz{S}_h}(\boz)$ and
$H_{\Phi,\,S_h}(\boz)$ coincide with equivalent norms.
\end{enumerate}
\end{theorem}

We first point out that the coincidence between $H_{\Phi,\,r}(\boz)$
and $H_{\Phi,\,\cn_h}(\boz)$ of Theorem \ref{t1.1} when $\Phi(t):=
t$ for all $t\in (0,\fz)$ was already obtained by Auscher and Russ
in \cite[Proposition 19, Theorems 1 and 20]{ar03}.

We also remark that although a strongly Lipschitz domain can be
regarded as a space of homogeneous type, Theorem \ref{t1.1} can not
be deduce from a general theory of Hardy spaces on spaces of
homogeneous type, since its proof strongly depends on the
geometrical property of strongly Lipschitz domains and
the divergence structure of the considered operator $L$.

The following chains of inequalities give the strategy of the proof
of Theorem \ref{t1.1}(i). For all $f\in H_{\Phi,\,r}(\boz)\cap L^2
(\boz)$, we have
\begin{equation}\label{1.6}
\|f\|_{H_{\Phi,\,r}(\boz)}\gs\|f\|_{H_{\Phi,\,\cn_h}(\boz)}\gs
\|f\|_{H_{\Phi,\,\wz{S}_h}(\boz)}\gs
\|f\|_{H_{\Phi,\,S_h}(\boz)}\gs\|f\|_{H_{\Phi,\,r}(\boz)},
\end{equation}
where the implicit constants are independent of $f$.  The proof of
the first inequality in \eqref{1.6} is standard by applying the
atomic decomposition of $H_{\Phi}(\rn)$ established by Viviani
\cite{vi87} and the relation between $H_{\Phi,\,r}(\boz)$ and
$H_{\Phi}(\rn)$; see Proposition \ref{p3.1} below. We prove the
second and the third inequalities, respectively, in Propositions
\ref{p3.2} and \ref{p3.3} below. We point out that Proposition
\ref{p3.2} plays an important role in the proof of Theorem
\ref{t1.1} and the key step in the proof of Proposition \ref{p3.2}
is to establish a ``good-$\lz$ inequality" concerning $\cn_h (f)$
and $\wz{S}_h (f);$ see Lemma \ref{l3.5} below. To show the last
inequality of \eqref{1.6} in Proposition \ref{p3.4}(i) below, for
all $f\in H_{\Phi,\,S_h}(\boz)\cap L^2 (\boz)$, we establish its
atomic decomposition by using a Calder\'on reproducing formula  
on $L^2 (\boz)$ associated with $L$ (see
\eqref{3.42} below), the atomic
decomposition of functions in the tent space on $\boz$, and the
reflection technology related to Lipschitz domains on $\rn$ which
was proved by Auscher and Russ in \cite[p.\,183]{ar03} and plays a
key role in the proof of Theorem \ref{t1.1} (see also Lemma
\ref{l3.9} below). But, this reflection technology was not necessary
in the study of the Orlicz-Hardy space $H_{\Phi,\,z}(\boz)$ in
\cite{yys} (see also \cite{ar03}).

Similarly to the proof of Theorem \ref{t1.1}(i), the following
chains of inequalities give the strategy of the proof of Theorem
\ref{t1.1}(ii), namely, we shall show that 
for all $f\in H_{\Phi,\,r}(\boz)\cap L^2 (\boz),$
\begin{eqnarray*}
\|f\|_{H_{\Phi,\,r}(\boz)}&&\gs\|f\|_{H_{\Phi,\,\cn_h}(\boz)}\gs
\|f\|_{H_{\Phi,\,\wz{S}_h,\,d_{\boz}}(\boz)}\\
&&\gs
\|f\|_{H_{\Phi,\,S_h,\,d_{\boz}}(\boz)}\gs\|f\|_{H_{\Phi,\,r}(\boz)},
\end{eqnarray*} where the implicit constants are independent of $f$.
In this case that $\boz$ is bounded, the Calder\'on reproducing
formula \eqref{3.42} on $L^2 (\boz)$ associated with $L$ used in the
proof of Theorem \ref{t1.1}(i) is never valid. Thus, instead of
\eqref{3.42}, we use a local Calder\'on reproducing formula 
on $L^2 (\boz)$ associated with $L$ (see
\eqref{3.72} below). Moreover,
if $\boz$ is bounded, $n\ge 3$ and $(G_{\fz})$ holds, using the fact
that the operator $L^{-1}$ is bounded from $L^p (\boz)$ into $L^q
(\boz)$ for some $p,\,q\in(1,\fz)$ satisfying $1<p<q<\fz$ and
$\frac1p-\frac1q=\frac2n,$
which can be proved by a way similar to the proof of
\cite[p.\,42, Proposition 5.3]{a07}, we further show that the second
term in \eqref{1.5} can be controlled by the Orlicz norm of the
Lusin area function $S_h (f)$, which implies the second part of
Theorem \ref{t1.1}(ii).

Let $\Phi$ satisfy Assumption (A), $\boz$ be a unbounded strongly
Lipschitz domain of $\rn$, and $L$ an
elliptic second-order divergence operator on $L^2(\boz)$ satisfying
the Neumann boundary condition and the Gaussian property
$(G_{\infty})$. As mentioned above, the Orlicz-Hardy space
$H_{\Phi,\,z}(\Omega)$ was
introduced in \cite{yys} and its several equivalent characterizations,
including the nontangential maximal function characterization
and the Lusin area function characterization associated with
$\{e^{-t\sqrt{L}}\}_{t\ge0}$, the vertical and the nontangential
maximal function characterizations associated with
$\{e^{-tL}\}_{t\ge0}$, and the Lusin area function characterization
associated with $\{e^{-tL}\}_{t\ge0}$, were also obtained therein.

For all $f\in L^2 (\boz)$ and $x\in\boz$, let
$$S_P (f)(x):=\lfz\{\int_{\wz{\bgz}(x)}\lfz|t\paz_t e^{-t\sqrt L}
(f)(y)\r|^2\frac{dy\,dt}{t|Q(x,t)\cap\boz|}\r\}^{1/2},$$
where
$$\wz{\bgz}(x):=\{(y,t)\in\boz\times(0,\fz):\,|x-y|<t\}.$$
Let
$$\wz{H}^1_{S_P}(\boz):=\lfz\{f\in L^2 (\boz):\
\|f\|_{H^1_{S_P}(\boz)}:=\|S_P (f)\|_{L^1 (\boz)}<\fz\r\}.$$
The {\it Hardy space} $H^1_{S_P}(\boz)$ is defined to be the completion
of $\wz{H}^1_{S_P}(\boz)$ in the norm $\|\cdot\|_{H^1_{S_P}(\boz)}$.
By applying the method used in the proof of Proposition
\ref{p3.4}(i) below, we also give a direct proof for the atomic
decomposition of elements in $H^1_{S_P}(\boz)$ in Proposition
\ref{p3.5} below, which gives an answer to the question asked by
Duong and Yan \cite[p.\,485, Remarks (iii)]{dy04} in the case that
$p=1$. (We point out that the Lusin area function $S_P$ was also
given in \cite[p.\,154]{ar03} via replaced $|Q(x,t)\cap\boz|$ by
$t^n$. This may be problematic in obtaining some estimates, like the
estimate in line 1 from the bottom of \cite[p.\,164]{ar03}, by
regarding $\boz$ as a space of homogeneous type when $\boz$ is {\it
bounded}.)

The layout of this paper is as follows. In Section \ref{s2}, we
first recall some properties of the divergence form elliptic
operator $L$ on $\rn$ or a strongly Lipschitz domain $\boz$, and
then describe some basic assumptions on $L$; then we describe some
basic assumptions on Orlicz functions and present some properties of
these functions. In Section \ref{s3}, we give the proof of Theorem
\ref{t1.1}.

Finally we make some conventions on notation. Throughout the whole
paper, $L$ always denotes the {\it second order divergence form
elliptic operator} as in \eqref{1.3}. We denote by $C$ a {\it
positive constant} which is independent of the main parameters, but
it may vary from line to line. We also use $C(\gz,\bz,\cdots)$ to
denote a {\it positive constant} depending on the indicated
parameters $\gz,$ $\bz$, $\cdots$. The {\it symbol} $A\ls B$ means
that $A\le CB$. If $A\ls B$ and $B\ls A$, then we write $A\sim B$.
The {\it symbol} $\lf s\rf$ for $s\in\rr$ denotes the maximal
integer not more than $s$; $Q(x, t)$ denotes a {\it closed cube} in
$\rn$ with {\it center} $x\in\rn$ and
{\it sidelength} $l(Q):= t$ and
$$CQ(x,t):= Q(x,Ct).$$
For any given normed spaces $\mathcal A$
and $\mathcal B$ with the corresponding norms $\|\cdot\|_{\mathcal
A}$ and $\|\cdot\|_{\mathcal B}$, ${\mathcal A}\subset{\mathcal B}$
means that for all $f\in \mathcal A$, then $f\in\mathcal B$ and
$\|f\|_{\mathcal B}\ls \|f\|_{\mathcal A}$. For any subset $G$ of
$\rn$, we denote by $G^\complement$ the {\it set $\rn\setminus G$};
for a measurable set $E$, denote by $\chi_{E}$ the {\it
characteristic function} of $E$. We also set $\nn:=\{1,\, 2,\,
\cdots\}$ and $\zz_+:=\nn\cup\{0\}$. For any
$\tz:=(\tz_{1},\ldots,\tz_{n})\in\zz_{+}^{n}$, let
$$|\tz|:=\tz_{1}+\cdots+\tz_{n}$$
and
$$\partial^{\tz}_x:=\frac{\partial^{|\tz|}}{\partial
{x_{1}^{\tz_{1}}}\cdots\partial {x_{n}^{\tz_{n}}}}.$$
For any sets $E$, $F\subset\rn$ and $z\in\rn$, let
$$\dist(E, F):=\inf_{x\in
E,\,y\in F}|x-y|$$
and
$$\dist(z, E):=\inf_{x\in E}|x-z|.$$

\section{\hspace{-0.7cm}{\bf.} Preliminaries\label{s2}}

In Subsection \ref{s2.1}, we first recall some properties of the
divergence form elliptic operator $L$ on $\rn$ or a strongly
Lipschitz domain $\boz$, and then describe some basic assumptions on
$L$; in Subsection \ref{2.2}, we describe some basic assumptions of
Orlicz functions and then present some properties of these
functions.

\subsection{\hspace{-0.6cm}{\bf.}
The divergence form elliptic operator $L$\label{s2.1}}

Let $L$ be as in \eqref{1.3}. Then $L$ generates a semigroup
$\{e^{-tL}\}_{t\ge0}$ of operators that is {\it analytic} (namely,
it has an extension to a complex half cone $|\mathrm{arg} z|<\mu$
for some $\mu\in(0,\pi/2)$) and {\it contracting} on $L^2 (\boz)$
(namely, for all $f\in L^2 (\boz)$ and $t\in(0,\fz)$,
$\|e^{-tL}f\|_{L^2 (\boz)}\le\|f\|_{L^2 (\boz)}$); see, for example,
\cite{o04} for the details. Also, $L$ has a unique {\it maximal
accretive square root} $\sqrt{L}$ such that $-\sqrt{L}$ generates an
analytic and $L^2 (\boz)$-contracting semigroup $\{P_t\}_{t\ge0}$
with $P_t := e^{-t\sqrt{L}}$, the {\it Poisson semigroup} for $L$;
see, for example, \cite{k66} for the details.

Now we recall the Gaussian property of $\{e^{-tL}\}_{t\ge0}$
introduced by Auscher and Russ \cite[Definition 3]{ar03} on a
strongly Lipschitz domain; see also \cite{at98,at01a}.

\begin{definition}\rm\label{d2.1}
Let $\boz$ be a strongly Lipschitz domain of $\rn$ and $L$ as in
\eqref{1.3}. Let $\bz\in(0,\fz]$. The semigroup generated by $L$ is
called to have the \emph{Gaussian property} $(G_{\bz})$, if the following
(i) and (ii) hold:

\begin{enumerate}

\item[\rm{(i)}] The \emph{kernel} of $e^{-tL}$, denoted
by $K_t$, is a measurable function on $\boz\times\boz$ and there
exist positive constants $C$ and $\az$ such that for all
$t\in(0,\bz)$ and all $x,\,y\in\boz$,
\begin{equation}\label{2.1}
|K_t (x,y)|\le\frac{C}{t^{n/2}}e^{-\az\frac{|x-y|^2}{t}};
\end{equation}

\item[\rm{(ii)}] For all $x\in\boz$ and $t\in(0,\bz)$, the functions
$y\mapsto K_t (x,y)$ and
$$y\mapsto K_t (y,x)$$
are H\"older
continuous in $\boz$ and there exist positive constants $C$ and
$\mu\in(0,1]$ such that for all $t\in(0,\bz)$ and
$x,\,y_1,\,y_2\in\boz$,
\begin{equation}\label{2.2}
|K_t (x,y_1)-K_t (x,y_2)|+|K_t (y_1,x)-K_t
(y_2,x)|\le\frac{C}{t^{n/2}} \frac{|y_1 -y_2 |^{\mu}}{t^{\mu/2}}.
\end{equation}
\end{enumerate}
\end{definition}

\begin{remark}\rm\label{r2.1}
(i) The assumption $(G_{\fz})$ is always satisfied if $L$ is the
Laplacian or real symmetric operators (under DBC or NBC) on $\rn$ or
on Lipschitz domains except under NBC with $\boz$ bounded; see, for
example, \cite{at01a}.

(ii) The assumption $(G_{\fz})$ implies that for all
$\bz\in(0,\fz)$, $(G_{\bz})$ holds. If $\bz$ is finite, by
\cite[p.\,178, Lemma A.1]{ar03} and the property of semigroups, we
know that $(G_{\bz})$ and $(G_1)$ are equivalent.
\end{remark}

The following well-known fact is a simple corollary of the
analyticity of the semigroup $\{e^{-tL}\}_{t\ge0}$.
We omit the details.

\begin{lemma}\label{l2.1}
Let $\bz\in(0,\fz]$. Assume that $L$ has the Gaussian property
$(G_{\bz})$. Then the estimate \eqref{2.1} also holds for $t\paz_t
K_t$.
\end{lemma}

\subsection{\hspace{-0.6cm}{\bf.} Orlicz functions\label{s2.2}}

Let $\Phi$ be a positive function on $\rr_{+}:=(0,\fz)$. The
function $\Phi$ is said to be of {\it upper type $p$} (resp. {\it
lower type $p$}) for some $p\in[0,\fz)$, if there exists a positive
constant $C$ such that for all $t\in[1,\fz)$ (resp. $t\in(0,1]$) and
$s\in(0,\fz)$,
\begin{equation}\label{2.3}
\Phi(st)\le Ct^p \Phi(s).
\end{equation}
Obviously, if $\Phi$ is of lower type $p$ for some $p\in(0,\fz)$,
then
$$\lim_{t\to0^+}\Phi(t)=0.$$
Thus, for the sake of convenience, if
it is necessary, we may assume that $\Phi(0)=0$. If $\Phi$ is of
both upper type $p_1$ and lower type $p_0$, then $\Phi$ is said to
be of {\it type $(p_0,\,p_1)$}.
The function $\Phi$ is said to be of {\it strictly lower type $p$}
if for all $t\in(0,1)$ and $s\in(0,\fz)$,
$$\Phi(st)\le t^p\Phi(s),$$
and we define
\begin{eqnarray}\label{2.4}
p_{\Phi}:=\sup&&\{p\in(0,\fz):\ \Phi(st)\le t^p\Phi(s) \\
\nonumber &&\hs\text{holds for all}\ t\in(0,1)\ \text{and}\
s\in(0,\fz)\}.
\end{eqnarray}
In what follows, $p_{\Phi}$ is called the {\it strictly critical
lower type index of $\Phi$}. We point out that if $p_{\Phi}$ is
defined as in \eqref{2.4}, then $\Phi$ is also of strictly lower
type $p_{\Phi}$; see \cite{jy10} for the proof.

Throughout the whole paper, we always assume that $\Phi$ satisfies
the following assumptions.

\medskip

\noindent {\bf Assumption (A).} \emph{Let $\mu$ be as in
\eqref{2.2}, and $\Phi$ a positive function defined on $\rr_{+}$
which is of upper type 1 and  strictly critical lower type $p_{\Phi}
\in(\frac{n}{n+\mu},1]$. Also assume that $\Phi$ is continuous,
strictly increasing and subadditive.}

\medskip

Let $p\in (\frac{n}{n+\mu},1]$ and $\Phi(t):= t^p$ or $\Phi(t):=
t^{p}\ln(e^4+t)$ for all $t\in(0,\fz)$. Then $\bfai$ satisfies
Assumption (A) with $p_{\bfai}=p$; see \cite{jy10,lyy} for some
other examples.

Notice that if $\Phi$ satisfies Assumption (A), then $\Phi(0)=0$.
For any positive function $\wz\Phi$ of upper type 1 and
$p_{\wz\Phi}\in (\frac{n}{n+1},1]$, if we set
$$\Phi(t):=\int_0^t\frac{\wz\Phi(s)}{s}\,ds$$
for all $t\in[0,\fz)$, then by
\cite[Proposition 3.1]{vi87}, $\Phi$ is equivalent to $\wz\Phi$;
namely, there exists a positive constant $C$ such that
$$C^{-1}\wz\Phi(t)\le\Phi(t)\le C\wz\Phi(t)$$
for all $t\in[0,\fz)$;
moreover, $\Phi$ is a strictly increasing, subadditive and
continuous function of  upper type 1 and strictly
critical lower type
$$p_{\Phi}\equiv p_{\wz\Phi}\in\lfz(\frac{n}{n+\mu},1\r].$$
Notice that all our results are invariant
on equivalent functions satisfying
Assumption (A). From this, we deduce that all results with $\bfai$
as in Assumption (A) also hold for all positive functions $\wz\bfai$
of type 1 and strictly critical lower type
$p_{\wz\bfai}\in(\frac{n}{n+\mu},1]$.

Since $\Phi$ is strictly increasing, we define the function $\rz(t)$
on $\rr_{+}$
by setting, for all $t\in(0,\fz)$,
\begin{equation}\label{2.5}
\rz(t):=\frac{t^{-1}}{\Phi^{-1}(t^{-1})},
\end{equation}
 where $\Phi^{-1}$ is the inverse function of
$\Phi$. Then the types of $\Phi$  and $\rz$ have the following
relation: \emph{If $0<p_0\le p_1\le1$ and $\Phi$ is an increasing
function, then $\Phi$ is of type $(p_0,p_1)$ if and only if $\rz$ is
of type $(p_1^{-1}-1,p_0^{-1}-1);$} see \cite{vi87} for its proof.

\section{\hspace{-0.7cm}{\bf.} Proof
of Theorem \ref{t1.1}\label{s3}}

In this section, we present the proof of Theorem \ref{t1.1}. To this
end, we need some auxiliary area functions as follows. Recall that
$d_{\boz}:= 2\diam(\boz)$. Let
$\az\in(0,\fz),\,\epz,\,R\in(0,d_{\boz})$ and $\epz<R$. For all
given $f\in L^2 (\boz)$ and $x\in\boz$, let
\begin{equation*}
\wz{S}_h^{\az}(f)(x):=\lfz\{\int_{\bgz_{\az}(x)}\lfz|t \nabla e^{-t^2 L}
(f)(y)\r|^2\, \frac{dy\,dt}{t|Q(x,t)\cap\boz|}\r\}^{1/2}
\end{equation*}
and
\begin{equation*}
\wz{S}^{\epz,\,R,\,\az}_h
(f)(x):=\lfz\{\int_{\bgz^{\epz,\,R}_{\az}(x)}\lfz|t\nabla e^{-t^2 L}
(f)(y)\r|^2\, \frac{dy\,dt}{t|Q(x,t)\cap\boz|}\r\}^{1/2},
\end{equation*}
where and in what follows, for all $x\in\boz$, $\bgz_{\az}(x)$ and
$\bgz^{\epz,\,R}_{\az}(x)$ are the {\it cone} and the {\it truncated
cone}, respectively, defined by
\begin{equation*}
\bgz_{\az}(x):=\{(y,t)\in\boz\times(0,d_{\boz}):\ |y-x|<\az t\}
\end{equation*}
and
\begin{equation*}
\bgz^{\epz,\,R}_{\az}(x):=\{(y,t)\in\boz\times(\epz, R):\ |y-x|<\az
t\}
\end{equation*}
for $\az\in(0,\fz)$ and $0<\epz<R<d_{\boz}$. \emph{When $\az=1$, denote
$\wz{S}_h^{\az}(f),\,\wz{S}_h^{\epz,\,R,\,\az}(f)$ and
$\bgz_{\az}(x)$ simply, respectively, by
$\wz{S}_h(f),\,\wz{S}_h^{\epz,\,R}(f)$ and $\bgz(x)$.} 

To show Theorem \ref{t1.1}, we first establish 
the following Proposition \ref{p3.1}.

\begin{proposition}\label{p3.1}
Let $\Phi$ satisfy Assumption (A), $\boz$ be a strongly Lipschitz
domain of $\rn$ and $L$ as in \eqref{1.3}. Assume that the semigroup
generated by $L$ has the Gaussian property $(G_{\diam(\boz)})$. Then
under DBC,
$$(H_{\Phi,\,r}(\boz)\cap L^2 (\boz))\subset
(H_{\Phi,\,\cn_h}(\boz)\cap L^2 (\boz))$$
and there exists a
positive constant $C$ such that for all $f\in H_{\Phi,\,r}(\boz)\cap
L^2 (\boz)$,
$$\|f\|_{H_{\Phi,\,\cn_h}(\boz)}\le
C\|f\|_{H_{\Phi,\,r}(\boz)}.$$
\end{proposition}

To show Proposition \ref{p3.1}, we need the atomic decomposition
characterization of the Orlicz-Hardy space $H_{\bfai}(\rn)$
established by Viviani in \cite{vi87}. To state this, we begin with
the notions of $(\rz,\,q,\,s)$-atoms and the atomic Orlicz-Hardy
space $H^{\rz,\,q,\,s}(\rn)$.

\begin{definition}\rm\label{d3.1}
Let $\Phi$ be as in Definition \ref{d1.1} and $\rz$ as in \eqref{2.5},
$q\in(0,\fz]$ and $s\in\zz_+$. A function $a$ is called a {\it
$(\rz,\,q,\,s)$-atom} if
\begin{enumerate}
  \item[(i)] $\supp a\subset Q$, where $Q$ is a closed cube of $\rn;$
  
  \item[(ii)] $$\|a\|_{L^q (\rn)}\le|Q|^{1/q-1}[\rz(|Q|)]^{-1};$$
  
  \item[(iii)] for all $\bz:=(\bz_1,\,\bz_2,\,\cdots,\,\bz_n)\in\zz^n_+$
with $|\bz|\le s$,
$$\int_{\rn}a(x)x^{\bz}\,dx=0.$$
\end{enumerate}
\end{definition}

Obviously, when $\bfai(t):= t$ for all $t\in(0,\fz)$, the
$(\rz,\,q,\,s)$-atom is just the classical $(1,\,q,\,s)$-atom; see,
for example, \cite{s93}.

\begin{definition}\rm\label{d3.2}
Let $p_0$ be as in Definition \ref{d1.1} and $\Phi,\,q$ and $\rz$ as
in Definition \ref{d3.1}, and
$$s:=\lf n(1/p_0 -1)\rf.$$
The {\it atomic Orlicz-Hardy space $H^{\rz,\,q,\,s}(\rn)$} is defined to be
the space of all distributions $f\in\cs'(\rn)$ that can be written
as $f=\sum_j b_j$ in $\cs'(\rn)$, where $\{b_j\}_j$ is a sequence of
constant multiples of $(\rz,\,q,\,s)$-atoms, with the constant
depending on $j$,  such that for each $j$, $\supp b_j \subset Q_j$
and
$$\sum_j |Q_j|\Phi\lfz(\frac{\|b_j\|_{L^q (\rn)}}{|Q_j|^{1/q}}\r)<\fz.$$
Define
$$\blz_q (\{b_j\}_j):=\inf\lfz\{\lz\in (0,\fz):\ \sum_j |Q_j|\Phi
\lfz(\frac{\|b_j\|_{L^q (\rn)}}{\lz|Q_j|^{1/q}}\r)\le1\r\}$$
and
$$\|f\|_{H^{\rz,\,q,\,s}(\rn)}:=\inf\{\blz_q (\{b_j\}_j)\},$$
where the infimum is taken over all decompositions of $f$ as above.
\end{definition}

The $(\rz,\,q,\,s)$-atom and the atomic Orlicz-Hardy space
$H^{\rz,\,q,\,s}(\rn)$ were introduced by Viviani \cite{vi87}, in
which the following Lemma \ref{l3.1} was also obtained
(see \cite[Theorem 2.1]{vi87}).

\begin{lemma}\label{l3.1}
Let $p_0$ be as in Definition \ref{d1.1} and $\Phi,\,q$ and $\rz$ as
in Definition \ref{d3.1}, and
$$s:=\lf n(1/p_0 -1)\rf.$$
Then the spaces $H_{\Phi}(\rn)$ and $H^{\rz,\,q,\,s}(\rn)$ coincide with
equivalent norms.
\end{lemma}

Now we prove Proposition \ref{p3.1} by applying Lemma \ref{l3.1}.

\begin{proof}[\bf{Proof of Proposition \ref{p3.1}.}]
Let $f\in H_{\Phi,\,r}(\boz)\cap L^2 (\boz)$. By the definition of
$H_{\Phi,\,r}(\boz)$, we know that there exists $\wz{f}\in
H_{\Phi}(\rn)$ such that $\wz{f}\big|_{\boz}=f$ and
\begin{equation}\label{3.1}
\lfz\|\wz{f}\r\|_{H_{\Phi}(\rn)}\ls\|f\|_{H_{\Phi,\,r}(\boz)}.
\end{equation}

To show Proposition \ref{p3.1}, we only need prove that for any
constant multiple of a $(\rz,\,\fz,\,0)$-atom $b$ supported in the
closed cube $Q_0:= Q(x_0,r_0)$,
\begin{equation}\label{3.2}
\int_{\boz}\Phi(\cn_h
(b)(x))\,dx\ls|Q_0|\Phi\lfz(\|b\|_{L^{\fz}(\rn)}\r).
\end{equation}
Indeed, for $\wz{f}\in H_{\Phi}(\rn)$, by Lemma \ref{l3.1}, there
exists a sequence $\{b_i\}_i$ of constant multiples of
$(\rz,\,\fz,\,0)$-atoms, with the constant depending on $i$, such
that $\wz{f}=\sum_i b_i$ in $\cs'(\rn)$ and
$$\blz_\fz(\{b_i\}_i)\sim\|\wz f\|_{H_{\Phi}(\rn)}.$$
Moreover, by the proof of \cite[Theorem 2.1]{vi87} and (2.15) in \cite[Lemma
(2.9)]{ms79}, we know that the supports of $\{b_i\}_i$ are of {\it
finite intersection property}. By this, $f\in L^2 (\boz)$,
$\wz{f}=\sum_i b_i$ in $\cs'(\rn)$ and $\wz{f}\big|_{\boz}=f$, we
obtain that $f=\sum_i b_i$ almost everywhere on $\boz$, which
further implies that
\begin{eqnarray*}
\int_{\boz}K_{t^2} (x,y)f(y)\,dy=\sum_i
\int_{\boz}K_{t^2} (x,y)b_i (y)\,dy.
\end{eqnarray*}
From this, we deduce that for all $x\in\boz$,
$$\cn_h (f)(x)\le\sum_i\cn_h (b_i)(x).$$
By this and the fact that $\Phi$ is strictly
increasing, continuous and subadditive, if \eqref{3.2} holds,  we
then have
\begin{eqnarray*}
\int_{\boz}\Phi\lfz(\cn_h (f)(x)\r)\,dx\le\sum_i
\int_{\boz}\Phi\lfz(\cn_h (b_i)(x)\r)\,dx\ls\sum_i
|Q_i|\Phi\lfz(\|b_i\|_{L^{\fz}(\rn)}\r),
\end{eqnarray*}
where for each $i$, $\supp b_i\subset Q_i$. This, together with the
facts that for all $\lz\in(0,\fz)$,
$$\cn_h (f/\lz)=\cn_h (f)/\lz$$
and for each $i$,
$$\|b_i/\lz\|_{L^{\fz}(\rn)}=\|b_i\|_{L^{\fz}(\rn)}/\lz,$$
implies that for all $\lz\in(0,\fz)$,
\begin{eqnarray*}
\int_{\boz}\Phi\lfz(\frac{\cn_h (f)(x)}{\lz}\r)\,dx\ls\sum_i
|Q_i|\Phi\lfz(\frac{\|b_i\|_{L^{\fz}(\rn)}}{\lz}\r).
\end{eqnarray*}
By this and \eqref{3.1}, we obtain that
$$\|f\|_{H_{\bfai,\,\cn_h}(\boz)}\ls\blz_\fz(\{b_i\}_i) \sim\|\wz
f\|_{H_{\Phi}(\rn)}\ls\|f\|_{H_{\bfai,\,r}(\boz)},$$
which, together
with the arbitrariness of $f\in H_{\bfai,\,r}(\boz)\cap L^2 (\boz)$,
implies the conclusions of Proposition \ref{p3.1}.

It is easy to see that for all $x\in\boz$,
\begin{equation}\label{3.3}
e^{-t^2 L}(b)(x)=\int_{Q_0 \cap\boz}K_{t^2} (x,y)b(y)\,dy.
\end{equation}
Now we show \eqref{3.2} by considering the following three cases for
$Q_0$.

{\it Case 1)} $Q_0\cap\boz=\emptyset$. In this case, by \eqref{3.3},
we know that for all $x\in\boz$, $\cn_h (b)(x)=0$. From this, it
follows that \eqref{3.2} holds.

{\it Case 2)} $Q_0\subset\boz$. In this case, let $\wz{Q}_0:= 8Q_0$.
Then we have
\begin{equation}\label{3.4}
\int_{\boz}\Phi(\cn_h (b)(x))\,dx=\int_{\wz{Q}_0\cap\boz}\Phi(\cn_h
(b)(x))\,dx+
\int_{(\wz{Q}_0)^\complement\cap\boz}\cdots=:\mathrm{I_1}+\mathrm{I_2}.
\end{equation}

We first estimate $\mathrm{I_1}$. For any $x\in\wz{Q}_0$, by
\eqref{3.3} and \eqref{2.1}, we have
$$\cn_h (b)(x)\le\sup_{y\in\boz,\,t\in(0,d_{\boz}),\,|x-y|<t}
\int_{\boz}|K_{t^2} (y,z)||b(z)|\,dz\ls\|b\|_{L^{\fz}(\rn)},$$
which, together with the upper type 1 property of $\bfai$,
implies that
\begin{equation}\label{3.5}
\mathrm{I_1}\ls
\int_{\wz{Q}_0}\Phi(\|b\|_{L^{\fz}(\rn)})\,dx\ls|Q_0|
\Phi\lfz(\|b\|_{L^{\fz}(\rn)}\r).
\end{equation}

Now we estimate $\mathrm{I_2}$. Let
$x\in(\wz{Q}_0)^\complement\cap\boz$, $t\in(0,d_{\boz})$ and
$y\in\boz$ satisfy $|x-y|<t$. By the moment condition of $b$ and
\eqref{3.3}, we have
\begin{equation}\label{3.6}
e^{-t^2 L}(b)(y)=\int_{Q_0}[K_{t^2} (y,z)-K_{t^2} (y,x_0)]b(z)\,dz.
\end{equation}
Since $p_{\Phi}\in(\frac{n}{n+\mu},1]$, there exists
$\wz{\mu}\in(0,\mu)$ such that $p_{\Phi}>\frac{n}{n+\wz{\mu}}$. Now
we estimate $e^{-t^2 L}(b)(y)$ by considering the following two
cases for $t$.

(i) $t<\frac{1}{4} |x-x_0|$. In this case, let $z\in Q_0$. Then
$$|x-x_0|\le|x-y|+|y-x_0|<\frac{1}{4}|x-x_0|+|y-x_0|,$$
which deduces
that $|x-x_0|<\frac{4}{3}|y-x_0|$. Moreover,
$$|x-x_0|\ge 4r_0\ge 4|z-x_0|.$$
Thus, we have
$$|y-x_0|\ge\frac{3}{4}|x-x_0|
\ge3 |z-x_0|,$$
which implies that
\begin{equation}\label{3.7}
|y-z|\ge|y-x_0|-|z-x_0|\ge\frac{2}{3}|y-x_0|\ge\frac{1}{2} |x-x_0|.
\end{equation}
Thus, by \eqref{3.7}, \eqref{2.1} and \eqref{2.2}, we obtain
\begin{equation*}
|K_{t^2} (y,z)-K_{t^2} (y,x_0)|\ls
\frac{|z-x_0|^{\wz{\mu}}}{|x-x_0|^{n+\wz{\mu}}},
\end{equation*}
which, together with \eqref{3.6}, implies that
\begin{equation}\label{3.8}
\lfz|e^{-t^2 L}
(b)(y)\r|\ls\frac{r_0^{n+\wz{\mu}}}{|x-x_0|^{n+\wz{\mu}}}
\|b\|_{L^{\fz}(\rn)}.
\end{equation}

(ii) $t\ge \frac{1}{4}|x-x_0|$. In this case, by \eqref{2.2}, we
obtain
$$|K_{t^2} (y,z)-K_{t^2} (y,x_0)|\ls\frac{|z-x_0|^{\mu}}
{t^{n+\mu}}\ls\frac{|z-x_0|^{\mu}}{|x-x_0|^{n+\mu}}\ls
\frac{|z-x_0|^{\wz{\mu}}}{|x-x_0|^{n+\wz{\mu}}},$$
which, together with  \eqref{3.6}, implies that \eqref{3.8} also
holds in this case.

By the estimates obtained in (i) and (ii), and the arbitrariness of
$y\in\boz$ satisfying $|x-y|<t$, we obtain that
$$\cn_h (b)(x)\ls\frac{r_0^{n+\wz{\mu}}}{|x-x_0|^{n+\wz{\mu}}}
\|b\|_{L^{\fz}(\rn)},$$ which, together with the lower type
$p_{\Phi}$ property of $\Phi$ and $p_{\Phi}>\frac{n}{n+\wz{\mu}}$,
implies that
\begin{eqnarray}\label{3.9}
\hs\mathrm{I_2}&\ls&\int_{4
r_0}^{d_{\boz}}\Phi\lfz(\frac{r_0^{n+\wz{\mu}}}
{s^{n+\wz{\mu}}}\|b\|_{L^{\fz}(\rn)}\r)s^{n-1}\,ds\\ \nonumber
&\ls&\Phi(\|b\|_{L^{\fz}(\rn)})r_0^{(n+\wz{\mu})p_{\Phi}} \int_{4
r_0}^{\fz}s^{n-(n+\wz{\mu})p_{\Phi}-1}\,ds\sim
|Q_0|\Phi(\|b\|_{L^{\fz}(\rn)}).
\end{eqnarray}
Thus, by \eqref{3.4}, \eqref{3.5} and \eqref{3.9}, we know that
\eqref{3.2} holds in this case.

{\it Case 3)} $Q_0\cap\paz\boz\neq\emptyset$. In this case, recall
that for any $x\in\boz$, $t\in(0,\fz)$ and $y\in\paz\boz$, $K_t
(x,y)=0$ (see, for example, \cite[p.\,156]{ar03}). Take $y_0\in
Q_0\cap\paz\boz$. Then we have that for any $x\in\boz$ and
$t\in(0,d_{\boz})$, $K_{t^2} (x,y_0)=0$, which further implies that
for any $x\in\boz$,
$$e^{-t^2 L}(b)(x)=\int_{Q_0\cap\boz}[K_{t^2}
(x,y)-K_{t^2} (x,y_0)]b(y)\,dy.$$
The remaining estimates are similar
to those of Case 2). We omit the details, which completes the proof
of Proposition \ref{p3.1}.
\end{proof}

To show Theorem \ref{t1.1}, we need the following key proposition.

\begin{proposition}\label{p3.2}
Let $\Phi,\,\boz$ and $L$ be as in Proposition \ref{p3.1}. Then
under DBC, there exists a positive constant $C$ such that for all
$f\in H_{\Phi,\,\cn_h}(\boz)\cap L^2 (\boz),$
$$\lfz\|\wz{S}_h(f)\r\|_{L^{\Phi}(\boz)}\le C\|f\|_{H_{\Phi,\,\cn_h}(\boz)}.$$
\end{proposition}

To show Proposition \ref{p3.2}, we need the following Lemmas
\ref{l3.2} through \ref{l3.7}.

\begin{lemma}\label{l3.2}
Let $\boz$ be a strongly Lipschitz domain of $\rn$ and $L$ as in
\eqref{1.3}, and
$$I_r(x_0,t_0):=(Q(x_0,r)\cap\boz)\times[t_0-cr^2,t_0],$$
where
$(x_0,t_0)\in\boz\times(4cr^2,\fz)$, $r\in(0,\fz)$ and $c$ is a
positive constant. If
$$\paz_t u_t=-Lu_t$$
in $I_{2r} (x_0,t_0)$, then
there exists a positive constant $C$, depending only on $\boz$, $c$
and $\dz$ in \eqref{1.1}, such that
\begin{equation}\label{3.10}
\int_{I_r (x_0,t_0)}|\nabla u_t (x)|^2\,dx\,dt\le
\frac{C}{r^2}\int_{I_{2r} (x_0,t_0)}|u_t (x)|^2\,dx\,dt.
\end{equation}
\end{lemma}

Lemma \ref{l3.2} is usually called the {\it Caccioppoli inequality},
whose proof is similar to that of \cite[Lemma 3(a)]{ku00}. We omit
the details.

\begin{remark}\rm\label{r3.1}
Let $\boz,\,L,\,x_0,\,t_0,\,r$, $c$ and $u_t$ be as in Lemma
\ref{l3.1} but with $t_0^2\in(4cr^2,\fz)$. Then by making the change
of variables in \eqref{3.10}, we see that
\begin{eqnarray*}
&&\int_{\sqrt{t^2_0-cr^2}}^{t_0}\int_{Q(x_0,r)\cap\boz}t|\nabla
u_{t^2} (x)|^2\,dx\,dt\\
&&\hs\ls\frac{1}{r^2}\int_{\sqrt{t^2_0-4cr^2}}^{t_0}
\int_{Q(x_0,2r)\cap\boz}t| u_{t^2}(x)|^2\,dx\,dt.
\end{eqnarray*}
\end{remark}

In \cite[p.\,183]{ar03}, Auscher and Russ proved the following
geometric property of strongly Lipschitz domains, which plays an
important role in this paper.

\begin{lemma}\label{l3.3}
Let $\boz$ be a strongly Lipschitz domain of $\rn$. Then there
exists a constant $C\in(0, 1]$ such that for all cubes $Q$ centered
in $\boz$ with $l(Q)\in(0,\,\infty)\cap(0,\, d_{\boz}]$, $|Q\cap
\boz|\ge C|Q|$.
\end{lemma}

In what follows, we denote by $B((z,\tau),r)$ the {\it ball in
$\rn\times(0,\fz)$ with center $(z,\tau)$ and radius $r$}; namely,
$$B((z,\tau),r):=\{(x,t)\in\rn\times(0,\fz):\
\max(|x-z|,\,|t-\tau|)<r\}.$$

\begin{lemma}\label{l3.4}
Let $\boz$ be a strongly Lipschitz domain of $\rn$, $\az\in(0,1)$,
$\epz,\,R\in(0,d_{\boz})$ and $\epz< R$. Then there exists a
positive constant $C$, depending only on $\az$, $\boz$ and $n$, such
that for all $f\in L^2 (\boz)$ and $x\in\boz$,
\begin{equation}\label{3.11}
\wz{S}^{\epz,\,R,\,\az}_h (f)(x)\le C[1+\ln(R/\epz)]^{1/2}\cn_h
(f)(x).
\end{equation}
\end{lemma}

\begin{proof}[\bf{Proof.}]\rm
Fix $\az\in(0,1)$, $0<\epz<R<d_{\boz}$ and $x\in\boz$. Let $f\in L^2
(\boz)$ and for all $t\in(0,d_{\boz})$, $u_t := e^{-t^2 L}(f)$. For
all $(z,\tau)\in\bgz^{\epz,\,R}_{\az}(x)$, let
$$E_{(z,\,\tau)}:=
B((z,\tau),\gz\tau)\cap(\boz\times(0,d_{\boz})),$$
where $\gz\in(0,1)$ is a positive constant which is determined later. By
the Besicovitch covering lemma, there exists a subcollection
$\{E_{(z_j,\tau_j)}\}_j$ of
$\{E_{(z,\tau)}\}_{(z,\tau)\in\bgz^{\epz,\,R}_{\az}(x)}$ such that
\begin{equation}\label{3.12}
\bgz^{\epz,\,R}_{\az}(x)\subset\bigcup_j E_{(z_j,\tau_j)}\
\text{and}\ \sum_j \chi_{E_{(z_j,\tau_j)}}\le M,
\end{equation}
where $M$ is a positive integer depending only on $n$. For each $j$,
we denote $E_{(z_j,\tau_j)}$ simply by $E_j$. Then we have the
following two facts for $E_j$:

\smallskip

(i) \emph{For each $j$, if $(y,t)\in E_j$, then $t\sim\tau_j\sim d_j$,
where $d_j$ denotes the distance from $E_j$ to the bottom boundary
$\boz\times\{0\}$.}

\smallskip

Indeed, if $(y,t)\in E_j$, we then have
$$(1-\gz)\tau_j<t<(1+\gz)\tau_j,$$ 
which implies that $t\sim\tau_j$.
By $d_j=(1-\gz)\tau_j$, we obtain that
$$d_j<t<(1+\gz)\tau_j=\frac{1+\gz}{1-\gz}d_j.$$
Thus, $t\sim d_j$.

\smallskip

(ii) \emph{For each $j$, let
$$\wz{E}_j:=B((z_j,\tau_j),9\gz\tau_j)\cap(\boz\times(0,d_{\boz})).$$
If $\gz\in(0,\frac{1-\az}{18})$, then
$\wz{E}_j\subset\bgz^{\epz/2,\,2R}(x)$.}

\smallskip

Indeed, for all $(y,t)\in\wz{E}_j$, since
$(z_j,\tau_j)\in\bgz^{\epz,\,R}_{\az}(x)$, we have that
$$|y-z_j|<9\gz\tau_j$$
and $|x-z_j|<\az\tau_j$. From this, it follows that
\begin{equation}\label{3.13}
|x-y|<|x-z_j|+|z_j-y|<(9\gz+\az)\tau_j.
\end{equation}
Moreover, by $|t-\tau_j|<9\gz\tau_j$, we know that
$$(1-9\gz)\tau_j<t<(1+9\gz)\tau_j,$$
which implies that
$\tau_j<\frac{t}{1-9\gz}$ if $\gz\in(0,1/9)$. From this and
\eqref{3.13}, it follows that $|x-y|<\frac{9\gz+\az}{1-9\gz}t$.
Thus, to make that $\wz{E}_j\subset\bgz(x)$, it suffices to choose
$\gz\in(0,\frac{1-\az}{18}]$. Furthermore, by the facts that for any
$j$ and $(y,t)\in\wz{E}_j$,
$$(1-9\gz)\tau_j<t<(1+9\gz)\tau_j,$$
and $\epz<\tau_j<R$, to make that $t\in(\epz/2,2R)$, it suffices to take
$\gz\in(0,\frac{1}{18})$. Thus, if we choose
$\gz\in(0,\frac{1-\az}{18}]$, we then have that for each $j$,
$\wz{E}_j\subset\bgz^{\epz/2,\,2R}(x)$.

Now we show \eqref{3.11}. By the fact that $R\in(0,d_{\boz})$ and
Lemma \ref{l3.3}, we know that for all $t\in(\epz,R)$,
$$|Q(x,t)\cap\boz|\sim t^n.$$
From this, \eqref{3.12}, the above two
facts (i) and (ii), and Remark \ref{r3.1} (in which, if
$\tau_j\in(\epz,\frac{d_{\boz}}{1+\gz}]$, we choose
$t_0:=(1+\gz)\tau_j,\,r:=\gz\tau_j$ and $c:=\frac{4}{\gz}$, and if
$\tau_j\in(\frac{d_{\boz}}{1+\gz},d_{\boz})$, we then choose $t_0:=
d_{\boz},\,r:=\gz\tau_j$ and $c:=\frac{4}{\gz(1+\gz)^2}$, and in
both cases, we need choose
$\gz\in(0,\min\{\frac{2}{81},\,\frac{1-\az}{18}\}$), it follows that
\begin{eqnarray*}
\lfz[\wz{S}^{\epz,\,R,\,\az}_h
(f)(x)\r]^2&\sim&\int_{\bgz^{\epz,\,R}_{\az}(x)}|t\nabla u_t
(y)|^2\frac{dy\,dt}{t^{n+1}}\ls\sum_j \int_{E_j}|t\nabla u_t
(y)|^2\frac{dy\,dt}{t^{n+1}}\\
&\ls&\sum_j \int_{(1-\gz)\tau_j}^{\min\{(1+\gz)\tau_j,\,d_{\boz}\}}
\int_{Q(z_j,\gz\tau_j)\cap\boz}t|\nabla
u_t (y)|^2\frac{dy\,dt}{t^n}\\
&\ls&\sum_j d_j^{-n}\frac{1}{(\gz\tau_j)^2}\int_{\wz{E}_j}t|u_t
(y)|^2\,dy\,dt\\
&\ls&\lfz\{\sum_j d_j^{-n}(\gz\tau_j)^{-2}|\wz{E}_j|(1+9\gz)\tau_j\r\}
\lfz[\cn_h (f)(x)\r]^2\\
&\sim&\lfz\{\sum_j \int_{E_j}\frac{dy\,dt}{t^{n+1}}\r\}\lfz[\cn_h
(f)(x)\r]^2\\
&\ls&\int_{\bgz^{\epz/2,\,2R}(x)}\frac{dy\,dt}{t^{n+1}}\lfz[\cn_h
(f)(x)\r]^2\\
&\ls&\int_{\epz/2}^{2R}\lfz\{\int_{\rn} \chi_{Q(0,1)}\lfz(\frac{x-y}{
t}\r)\,dy\r\} t^{-(n+1)}\,dt\lfz[\cn_h
(f)(x)\r]^2\\
&\sim&[1+\ln(R/\epz)]\lfz[\cn_h (f)(x)\r]^2,
\end{eqnarray*}
which implies that
$$\wz{S}^{\epz,\,R,\,\az}_h
(f)(x)\ls[1+\ln(R/\epz)]^{1/2}\cn_h (f)(x).$$
Thus, \eqref{3.11}
holds, which completes the proof of Lemma \ref{l3.4}.
\end{proof}

\begin{lemma}\label{l3.5}
Let $\boz$ be a strongly Lipschitz domain of $\rn$ and $L$ as in
\eqref{1.3}, and $d_{\boz}:=2\diam(\boz)$. Then there exists a
positive constant $C$ such that for all $\gz\in(0,1]$,
$\lz\in(0,\fz)$, $\epz,\,R\in(0,d_{\boz})$ with $\epz<R$ and $f\in
H_{\Phi,\,\cn_h}(\boz)\cap L^2 (\boz)$,
\begin{eqnarray}\label{3.14}
&&\lfz|\lfz\{x\in\boz:\ \wz{S}^{\epz,\,R,\,1/20}_h (f)(x)>2\lz,\,\cn_h
(f)(x)\le\gz\lz\r\}\r|\\ \nonumber &&\hs\le C\gz^2\lfz|\lfz\{x\in\boz:\
\wz{S}^{\epz,\,R,\,1/2}_h (f)(x)>\lz\r\}\r|.
\end{eqnarray}
\end{lemma}

We point out that in the proof of Proposition \ref{p3.2}, Lemma
\ref{l3.5} plays a key role. The inequality \eqref{3.14} is usually
called  the ``\emph{good-$\lz$ inequality}" concerning the maximal function
$\cn_h (f)$ and the truncated area functions
$\wz{S}^{\epz,\,R,\,1/20}_h (f)$ and $\wz{S}^{\epz,\,R,\,1/2}_h
(f)$.

\begin{proof}[\bf{Proof of Lemma \ref{l3.5}.}]
To prove this lemma, we borrow some ideas from \cite{amr08} and
\cite{ar03}.

Fix $0<\epz<R<d_{\boz}$, $\gz\in(0,1]$ and
$\lz\in(0,\fz)$. Let $f\in H_{\Phi,\,\cn_h}(\boz)\cap L^2(\boz)$ and
$$O:=\lfz\{x\in\boz:\ \wz{S}^{\epz,\,R,\,1/2}_h (f)(x)>\lz\r\}.$$
It is easy to show that $O$ is an open subset of $\boz$.

Now we show \eqref{3.14} by considering the following two cases for
$O$.

{\it Case 1)} $O\neq\boz$. In this case, let
\begin{eqnarray}\label{3.15}
O=\bigcup_{k}(Q_k\cap\boz)
\end{eqnarray}
be the Whitney decomposition of $O$, where $\{Q_k\}_k$ are dyadic
cubes of $\rn$ with disjoint interiors and $(2Q_k)\cap\boz\subset
O\subset\boz$, but
$$((4Q_k)\cap\boz)\cap(\boz\setminus
O)\neq\emptyset.$$
To show \eqref{3.14}, by \eqref{3.15} and the
disjoint property of $\{Q_k\}_k$, it suffices to show that for all
$k$,
\begin{eqnarray}\label{3.16}
&&\lfz|\lfz\{x\in Q_k\cap\boz:\ \wz{S}^{\epz,\,R,\,1/20}_h
(f)(x)>2\lz,\, \cn_h (f)(x)\le\gz\lz\r\}\r|\\
\nonumber&&\hs\ls\gz^2
|Q_k\cap\boz|.
\end{eqnarray}
From now on, we fix $k$ and denote by $l_k$ the sidelength of $Q_k$.

If $x\in Q_k\cap\boz$, then
\begin{eqnarray}\label{3.17}
\wz{S}^{\max\{10l_k,\,\epz\},\,R,\,1/20}_h (f)(x)\le\lz.
\end{eqnarray}
Indeed, pick $x_k\in(4Q_k)\cap\boz$ with $x_k\not\in O$. For any
$(y,t)\in\boz\times(0,d_{\boz})$, if $|x-y|<\frac{t}{20}$ and
$t\ge\max\{10l_k,\,\epz\}$, then
$$|x_k-y|\le|x_k-x|+|x-y|<4l_k+\frac{t}{20}<\frac{t}{2},$$
which implies that
$$\bgz^{\max\{10l_k,\,\epz\},\,R}_{1/20}(x)\subset
\bgz^{\max\{10l_k,\,\epz\},\,R}_{1/2}(x_k).$$
By this, we obtain that
$$\wz{S}^{\max\{10l_k,\,\epz\},\,R,\,1/20}_h (f)(x)\le
\wz{S}^{\max\{10l_k,\,\epz\},\,R,\,1/2}_h (f)(x_k)\le\lz.$$
Thus, the claim \eqref{3.17} holds.

If $\epz\ge10l_k$, by \eqref{3.17}, we see that \eqref{3.16} holds.
If $\epz<10l_k$, to show \eqref{3.16}, by the fact that
$$\wz{S}^{\epz,\,R,\,1/20}_h (f)\le \wz{S}^{\epz,\,10l_k,\,1/20}_h
(f)+\wz{S}^{10l_k,\,R,\,1/20}_h (f)$$
and \eqref{3.17}, it remains to
show that
\begin{eqnarray}\label{3.18}
\lfz|\{x\in Q_k\cap F:\ g(x)>\lz\}\r|\ls\gz^2|Q_k\cap\boz|,
\end{eqnarray}
where $g:= \wz{S}^{\epz,\,10l_k,\,1/20}_h (f)$ and
$$F:=\{x\in\boz:\ \cn_h (f)(x)\le\gz\lz\}.$$
By Chebyshev's inequality, we see that
\eqref{3.18} is deduced from
\begin{eqnarray}\label{3.19}
\int_{Q_k\cap F}[g(x)]^2\,dx\ls(\gz\lz)^2|Q_k\cap\boz|.
\end{eqnarray}
Now we prove \eqref{3.19}. It is easy to see that $F$ is a closed
subset of $\boz$.

If $\epz\ge5l_k$, then by the definitions of $g$ and $F$ and Lemma
\ref{l3.4}, we have
\begin{eqnarray*}
\int_{Q_k\cap F}[g(x)]^2\,dx&&\ls\int_{Q_k\cap F} \lfz[\cn_h
(f)(x)\r]^2\,dx\\
&&\ls(\gz\lz)^2|Q_k\cap F|\ls(\gz\lz)^2|Q_k\cap \boz|,
\end{eqnarray*}
which shows \eqref{3.19} in this case.

Assume from now on that $\epz<5l_k$. Let
\begin{eqnarray}\label{3.20}
G:=\lfz\{(y,t)\in\boz\times(\epz,\min\{10l_k,\,d_{\boz}\}):\
\pz(y)<\frac{t}{20}\r\},
\end{eqnarray}
where
\begin{eqnarray}\label{3.21}
\pz(y):=\dist(y,Q_k\cap F).
\end{eqnarray}
By the geometric properties of $\boz$, we have
\begin{eqnarray}\label{3.22}
\int_{Q_k\cap F}[g(x)]^2\,dx\ls\int_{G}t|\nabla u_t (y)|^2\,dy\,dt.
\end{eqnarray}
Indeed, if $\boz$ is unbounded, by Lemma \ref{l3.3}, we know that
for all $x\in\boz$ and $t\in(0,\fz)$,
$$|Q(x,t)\cap\boz|\sim|Q(x,t)|.$$
Thus, in this case, we have
\begin{eqnarray*}
\int_{Q_k\cap F}[g(x)]^2\,dx&=&\int_{Q_k\cap
F}\lfz\{\int_{\bgz^{\epz,\,10l_k}_{1/20}(x)}\lfz|t\nabla u_t
(y)\r|^2\frac{dy\,dt}{t|Q(x,t)\cap\boz|}\r\}\,dx\\
&\ls&\int_{G}\lfz\{\int_{\boz}t^{1-n}\chi_{Q(0,1)}
\lfz(\frac{20(x-y)}{t}\r)\,dx\r\}\lfz|\nabla u_t (y)\r|^2\,dy\,dt\\
&\ls&\int_{G}t\lfz|\nabla u_t (y)\r|^2\,dy\,dt.
\end{eqnarray*}
That is, \eqref{3.22} holds in this case. If $\boz$ is bounded, we
first assume that $\diam(\boz)\le10l_k$. Then
\begin{eqnarray*}
\int_{Q_k\cap F}[g(x)]^2\,dx&=&\int_{Q_k\cap
F}\lfz\{\int_{\bgz^{\epz,\,\diam(\boz)}_{1/20}(x)}\lfz|t\nabla u_t
(y)\r|^2\frac{dy\,dt}{t|Q(x,t)\cap\boz|}\r\}\,dx\\
&&+\int_{Q_k\cap F}\lfz\{\int_{\bgz^{\diam(\boz),\,10l_k}_{1/20}}
\cdots\r\}\,dx\\
&\ls&\int_{G}\lfz\{\int_{\boz}t^{1-n}\chi_{Q(0,1)}
\lfz(\frac{20(x-y)}{t}\r)\,dx\r\}\lfz|\nabla u_t (y)\r|^2\,dy\,dt\\
&&+\int_{G}\lfz\{\frac{1}{|\boz|}\int_{\boz}t\chi_{Q(0,1)}
\lfz(\frac{20(x-y)}{t}\r)\,dx\r\}\lfz|\nabla u_t (y)\r|^2\,dy\,dt\\
&\ls&\int_{G}t\lfz|\nabla u_t (y)\r|^2\,dy\,dt,
\end{eqnarray*}
which is desired. If $\boz$ is bounded and $\diam(\boz)>10l_k$, then
$$g\le\wz{S}^{\epz,\,\diam(\boz),\,1/20}_h (f),$$
which, together with
an argument similar to the above, shows that \eqref{3.22} also holds
in this case. Thus, \eqref{3.22} is always true.

Let
$$E:=\lfz\{y\in\boz:\ \text{there exists}\
t\in(\epz,\min\{10l_k,\,d_{\boz}\})\ \text{such that}\
\pz(y)<\frac{t}{20}\r\}.$$ Then $E\subset 2Q_k\cap\boz$. Indeed, if
$y\in E$, then there exist $t\in(\epz,\min\{10l_k,\,d_{\boz}\})$
such that $(y,t)\in G$ and $x\in Q_k\cap F$ such that
$|x-y|<\frac{t}{20}$. By $t<10l_k$, we have
$|x-y|<\frac{10l_k}{20}=\frac{l_k}{2}$, which implies that
$E\subset2Q_k\cap\boz$.

Let
$$\wz{G}:=\lfz\{(y,t)\in\boz\times\lfz(\frac{\epz}{5},\min\{40l_k,\,
d_{\boz}\}\r):\ \pz(y)<t\r\}.$$
Then for all $(y,t)\in\wz{G}$,
\begin{equation}\label{3.23}
|u_t (y)|\le\gz\lz.
\end{equation}
Indeed, for any $(y,t)\in\wz{G}$, there exists $x\in Q_k\cap F$ such
that $|x-y|<t$ with $t\in(\frac{\epz}{5},\min\{40l_k,\,d_{\boz}\})$,
which implies that $(y,t)\in\bgz(x)$. Thus, by the definitions of
$F$ and $\cn_h (f)$, we have
$$|u_t (y)|\le \cn_h (f)(x)\le\gz\lz.$$

To finish the proof of Lemma \ref{l3.5}, we need the following
conclusion, which is just \cite[Lemma 3.5]{yys}.

\begin{lemma}\label{l3.6}
Let
$$D :=\lfz\{(y,t)\in\boz\times(\epz, 10l_k):\ \pz(y)<\frac{t}{20}\r\}$$
and
$$D_1 :=\lfz\{(y,t)\in\boz\times\lfz(\frac{\epz}{2}, 20l_k\r):\
\pz(y)<\frac{t}{10}\r\},$$
where $\pz$ is as in \eqref{3.21}. Then
there exists $\wz\zez\in C^{\fz}(D_1)\cap C(\ol{D}_1)$ satisfying
 that $0\le\wz\zez\le1$, $\wz\zez\equiv1$ on $D$,
$|\nabla\wz\zez(y,t)|\ls\frac{1}{t}$ for all $(y,t)\in D_1$, and
$$\supp\wz\zez\subset
D_1\cup\lfz\{\paz\boz\times\lfz(\frac{\epz}{2},20l_k\r)\r\},$$
where and in what
follows, $\ol{D}_1$ denotes the closure of $D_1$ in $\rr^{n+1}$.
\end{lemma}

Now we continue proving Lemma \ref{l3.5} by using Lemma \ref{l3.6}.
Let
$$G_1:=\lfz\{(y,t)\in\boz\times\lfz(\frac{\epz}{2},
\min\{20l_k,\,d_{\boz}\}\r):\ \pz(y)<\frac{t}{10}\r\}$$ and
$\wz\zez$ be as in Lemma \ref{l3.6}. Let
$\zez:=\wz\zez\chi_{\ol{\boz}\times(0,d_{\boz}]}$. Then $\zez\in
C^{\fz}(G_1)\cap C(\ol{G}_1)$, $0\le\zez\le1$, $\zez\equiv1$ on $G$,
 $|\nabla\zez(y,t)|\ls\frac{1}{t}$ for all $(y,t)\in G_1$, and
$$\supp\zez\subset
G_1\cup\lfz\{\paz\boz\times\lfz(\frac{\epz}{2},\min\{20l_k,d_{\boz}\}\r)\r\}.$$

Recall that $u_t := e^{-t^2 L}(f)$ for all $t\in(0,d_{\boz})$. By
$0\le\zez\le1$, $\zez\equiv1$ on $G$ and \eqref{1.1}, we have
\begin{eqnarray}\label{3.24}
\int_{G}t|\nabla u_t (y)|^2\,dy\,dt&\le&\int_{G_1} t|\nabla u_t
(y)|^2\zez(y,t)\,dy\,dt\\ \nonumber &\le& \dz^{-1}\Re\int_{G_1}t
A(y) \nabla u_t (y)\cdot\nabla\ol{u_t (y)}\zez(y,t)\,dy\,dt\\
\nonumber &=:& \dz^{-1}\Re\mathrm{I},
\end{eqnarray}
where $A(y)$ and $\dz$ are as in \eqref{1.1}. Let
$$\mathrm{J}:=\int_{G_1}tA(y)\nabla u_t (y)\cdot\nabla \zez(y,t)
\ol{u_t (y)}\,dy\,dt.$$ For all
$t\in(\epz/2,\min\{20l_k,\,d_{\boz}\})$ and all $y\in\boz$, let
$\zez_t (y):=\zez(y,t)$. Then $\zez_t \in C^{\fz}(\boz)$. By
\cite[p.\,23, (1.19)]{o04}, we know that for all $t\in(0,d_{\boz})$,
$$u_t\in D(L)\subset W^{1,\,2}_0 (\boz),$$
which, together with
$\zez_t\in C^{\fz}(\boz)$, implies that for all $t\in(0,d_{\boz})$,
$u_t \zez_t\in W^{1,\,2}_0 (\boz)$. From this, \eqref{1.2} and the
fact that
$$\paz_t u_t +2tLu_t=0$$
in $L^2 (\boz)$, it follows that
\begin{eqnarray}\label{3.25}
\hs\mathrm{I}&=&\int_{G_1}tA(y)\nabla u_t (y)\cdot\nabla \ol{u_t
(y)}\zez(y,t)\,dy\,dt\\ \nonumber &=&\int_{G_1}tA(y)\nabla u_t
(y)\cdot\nabla(\ol{u}_t\zez_t)(y)\,dy\,dt\\
\nonumber &&-\int_{G_1}tA(y)\nabla u_t (y)\cdot\nabla \zez_t (y)
\ol{u_t (y)}\,dy\,dt\\ \nonumber &=&\int_{G_1}tLu_t
(y)(\ol{u}_t\zez_t)(y)\,dy\,dt-\mathrm{J}\\ \nonumber
&=&-\frac{1}{2}\int_{G_1}\paz_t u_t
(y)(\ol{u}_t\zez_t)(y)\,dy\,dt-\mathrm{J}
=:-\frac{1}{2}\mathrm{I}_1-\mathrm{J}.
\end{eqnarray}

For $\mathrm{I}_1$, by the fact that $2\Re((\paz_t u_t)
\ol{u}_t)=\paz_t |u_t|^2$ and integral by parts, we have that
\begin{eqnarray*}
\Re\mathrm{I}_1 &=&\frac{1}{2}\int_{G_1}\paz_t |u_t (y)|^2\zez(y,t)\,dy\,dt\\
&=&\frac{1}{2}\lfz\{\int_{\paz G_1} |u_t (y)|^2
\zez(y,t)N(y,t)\cdot(0,\,0,\,\cdots,\,1)\,d\sz(y,t)\r.\\
&&-\lfz.\int_{G_1} |u_t (y)|^2\paz_t\zez(y,t)\,dy\,dt\r\},
\end{eqnarray*}
where $\paz G_1$ denotes the boundary of $G_1$, $N(y,t)$ the unit
normal vector outward $G_1$ and $d\sz$ the surface measure over
$\paz G_1$. This, combined with \eqref{3.25}, implies that
\begin{eqnarray}\label{3.26}
\Re\mathrm{I}&=&-\frac{1}{2}\Re\mathrm{I_1}-\Re\mathrm{J}\\
\nonumber &=&\frac{1}{4}\lfz\{\int_{G_1} |u_t
(y)|^2\paz_t\zez(y,t)\,dy\,dt\r.\\ \nonumber &&-\lfz.\int_{\paz G_1}
|u_t(y)|^2 \zez(y,t)N(y,t)\cdot(0,\,0,\,\cdots,\,1)\,d\sz(y,t)\r\}\\
\nonumber &&-\Re\lfz\{\int_{G_1}tA(y)\nabla u_t (y)\cdot\nabla \zez_t
(y) \ol{u_t (y)} \,dy\,dt\r\}.
\end{eqnarray}
By  $\supp\zez\subset
G_1\cup\{\paz\boz\times(\frac{\epz}{2},\min\{20l_k,\,d_{\boz}\})\}$
and  the fact that
$$N(y,t)\cdot(0,\cdots,0,1)=0$$
on $\paz\boz\times(\frac{\epz}{2},\min\{20l_k,\,d_{\boz}\})$, we
obtain
\begin{eqnarray}\label{3.27}
&&\int_{\paz G_1}|u_t
(y)|^2\zez(y,t)N(y,t)\cdot(0,\cdots,0,1)\,d\sz(y,t)=0.
\end{eqnarray}
From $\zez\equiv 1$ on $G$, we deduce that $\nabla\zez\equiv 0$ on $G$. Thus,
by this, \eqref{3.26} and \eqref{3.27}, we have
\begin{eqnarray}\label{3.28}
\hs\hs\Re\mathrm{I}&=&\frac{1}{4}\int_{G_1\setminus G}|u_t
(y)|^2\paz_t\zez(y,t)\,dy\,dt\\ \nonumber
&&-\Re\lfz\{\int_{G_1\setminus G}tA(y)\nabla u_t (y)\cdot\nabla \zez_t
(y) \ol{u_t (y)}\,dy\,dt\r\}\\ \nonumber
&=&:\mathrm{I_2}+\mathrm{I_3}.
\end{eqnarray}

First, we estimate $\mathrm{I_2}$. By $G_1\subset\wz{G}$ and
\eqref{3.23}, we obtain that for all $(y,t)\in G_1\setminus G$,
$|u_t (y)|\le\gz\lz$. Moreover,
\begin{eqnarray*}
&&G_1\setminus G\\
&&\hs=\lfz\{(y,t)\in\boz\times
\lfz(\frac{\epz}{2},\min\{20l_k,\,d_{\boz}\}\r):\ \frac{t}{20}
\le\pz(y)<\frac{t}{10}\r\}\\
&&\hs\hs\bigcup\lfz\{(y,t)\in\boz\times\lfz(\frac{\epz}{2},
\min\{20l_k,\,d_{\boz}\}\r):\
\pz(y)<\frac{t}{10},\,\frac{\epz}{2}\le
t<\epz\r\}\\
&&\hs\hs\bigcup\lfz\{(y,t)\in\boz\times\lfz(\frac{\epz}{2},
\min\{20l_k,\,d_{\boz}\}\r):\ \pz(y)<\frac{t}{10},\,10l_k\le
t<20l_k\r\}.
\end{eqnarray*}
From these observations and the fact that for all $(y,t)\in G_1$,
$$|\nabla\zez(y,t)|\ls\frac{1}{t},$$
we deduce that
\begin{eqnarray}\label{3.29}
&&\int_{G_1\setminus G}|u_t (y)|^2|\paz_t \zez(y,t)|\,dy\,dt\\
\nonumber &&\hs\ls(\gz\lz)^2\int_{G_1\setminus
G}\,\frac{dy\,dt}{t}\\
\nonumber&&\hs\ls(\gz\lz)^2\int_{H_1}\lfz\{\int_{\epz/2}^{\epz}
\frac{\,dt}{t}+\int_{10l_k}^{20l_k}
\frac{\,dt}{t}+\int_{10\pz(y)}^{20\pz(y)} \frac{\,dt}{t} \r\}\,dy\\
\nonumber
&&\hs\ls(\gz\lz)^2|H_1|,\nonumber
\end{eqnarray}
where
\begin{eqnarray*}
H_1:=&&\lfz\{y\in G_1:\ \text{there exists}\
t\in\lfz(\frac{\epz}{2},\min\{20l_k,\,d_{\boz}\}\r)\r.\\
&&\hspace{1cm}\ \text{such that}\
(y,t)\in G_1\Big\}.
\end{eqnarray*}
For all $y\in H_1$, we know that there exists
$t\in(\frac{\epz}{2},\min\{20l_k,\,d_{\boz}\})$ such that $(y,t)\in
G_1$. From this and the definition of $G_1$, it follows that there
exists $x\in Q_k\cap F$ such that $|x-y|<\frac{t}{10}$ with
$t\in(\frac{\epz}{2},\min\{20l_k,\,d_{\boz}\})$. Thus, $|x-y|<2l_k$,
which implies that $y\in(5Q_k)\cap\boz$. By this, we know that
$H_1\subset(5Q_k)\cap\boz$, which together with \eqref{3.19} implies
that
\begin{eqnarray}\label{3.30}
|\mathrm{I_2}|&&\ls\int_{G_1\setminus G}|u_t (y)|^2|\paz_t
\zez(y,t)|\,dy\,dt\ls(\gz\lz)^2|H_1|\\
\nonumber &&\ls(\gz\lz)^2|Q_k\cap\boz|.
\end{eqnarray}

To estimate $\mathrm{I_3}$, by the facts that
$|\nabla\zez(y,t)|\ls\frac{1}{t}$ for all $(y,t)\in G_1$ and that
$|u_t (y)|\le\gz\lz$ for all $(y,t)\in G_1$, we have
\begin{eqnarray}\label{3.31}
|\mathrm{I_3}|\ls\int_{G_1\setminus G}\lfz|\nabla u_t (y)\r||u_t
(y)|\,dy\,dt\ls\gz\lz\int_{G_1\setminus G}\lfz|\nabla u_t
(y)\r|\,dy\,dt.
\end{eqnarray}
Now, we need show
\begin{eqnarray}\label{3.32}
\int_{G_1\setminus G}\lfz|\nabla u_t
(y)\r|\,dy\,dt\ls\gz\lz|Q_k\cap\boz|.
\end{eqnarray}
For all $(y,t)\in (G_1\setminus G)$ and $\dz_1\in(0,1)$, let
$$E_{(y,t)}:= B((y,t),\dz_1 t)\cap(\boz\times(0,d_{\boz}))$$
and
$$\wz{E}_{(y,t)}:= B((y,t),9\dz_1 t)\cap(\boz\times(0,d_{\boz})).$$
Take $\dz_1$ small enough such that for all $(y,t)\in (G_1\setminus
G)$,
\begin{eqnarray*}
\wz{E}_{(y,t)}\subset&&
\lfz\{(y,t)\in\boz\times\lfz(\frac{\epz}{5},\min\{30l_k,\,d_{\boz}\}\r):\
\frac{t}{40}
<\pz(y)<\frac{t}{2}\r\}\\
&&\bigcup\lfz\{(y,t)\in\boz\times\lfz(\frac{\epz}{5},
\min\{30l_k,\,d_{\boz}\}\r):\ \pz(y)<\frac{t}{2},\,\frac{\epz}{5}<
t<2\epz\r\}\\
&&\bigcup\lfz\{(y,t)\in\boz\times\lfz(\frac{\epz}{5},
\min\{30l_k,\,d_{\boz}\}\r):\ \pz(y)<\frac{t}{2},\,5l_k<
t<30l_k\r\}\\
=&&: G_2.
\end{eqnarray*}
By the Besicovith covering lemma, there exists a sequence
$\{E_{(y_j,t_j)}\}_j$ of sets which are a bounded covering of
$G_1\setminus G$. Let $E_j:= E_{(y_j,t_j)}$ and $\wz{E}_j:=
\wz{E}_{(y_j,t_j)}$. Notice that for all $(y,t)\in E_j$, $t\sim
t_j\sim r(E_j)$, where $ r(E_j)$ denotes the radius of $E_j$. From
this, H\"older's inequality, Remark \ref{r3.1} (in which, if
$\tau_j\in(\epz,\frac{d_{\boz}}{1+\dz_1}]$, we choose
$t_0:=(1+\dz_1)\tau_j,\,r:=\dz_1\tau_j$ and $c:=\frac{4}{\dz_1}$,
and if $\tau_j\in(\frac{d_{\boz}}{1+\dz_1},d_{\boz})$, we then
choose $t_0:= d_{\boz},\,r:=\dz_1\tau_j$ and
$c:=\frac{4}{\dz_1(1+\dz_1)^2}$, and in both cases, we need
$\dz_1\in(0,2/81)$) and the fact that for all $j$ and $(y,t)\in
\wz{E}_j$, $|u_t (y)|\le\gz\lz$, it follows that
\begin{eqnarray}\label{3.33}
&&\int_{G_1\setminus G}\lfz|\nabla u_t (y)\r|\,dy\,dt\\ \nonumber
&&\hs\ls\sum_j \int_{E_j}\lfz|\nabla u_t
(y)\r|\,dy\,dt\\ \nonumber
&&\hs\ls\sum_j|E_j|^{1/2}\lfz\{\int_{E_j}\lfz|\nabla u_t
(y)\r|^2\,dy\,dt\r\}^{1/2}\\ \nonumber &&\hs\ls\sum_j|E_j|^{1/2}
[r(E_j)]^{-1}\lfz\{\int_{\wz{E}_j}|u_t (y)|^2\,dy\,dt\r\}^{1/2}\\
\nonumber &&\hs\ls\gz\lz\sum_j|E_j|[r(E_j)]^{-1}
\ls\gz\lz\int_{G_2}\frac{dy\,dt}{t}\\ \nonumber
&&\hs\ls\gz\lz\int_{H_2}\lfz\{\int_{\epz/5}^{\epz}
\frac{\,dt}{t}+\int_{5l_k}^{30l_k}
\frac{\,dt}{t}+\int_{2\pz(y)}^{40\pz(y)} \frac{\,dt}{t}
\r\}\,dy\\
&&\hs\ls\gz\lz|H_2|,\nonumber
\end{eqnarray}
where
$$H_2:=\lfz\{y\in\boz:\ \text{there exists}\
t\in\lfz(\frac{\epz}{5},30l_k\r)\  \text{such that}\ (y,t)\in G_2\r\}.$$
Similarly to the estimate of $H_1$, we also have
$|H_2|\ls|Q_k\cap\boz|$, which, together with \eqref{3.33}, implies
that \eqref{3.32} holds. Thus, by \eqref{3.31} and \eqref{3.32}, we
obtain that
$$|\mathrm{I_3}|\ls(\gz\lz)^2|Q_k\cap\boz|,$$
which, together with \eqref{3.22}, \eqref{3.28} and \eqref{3.30}, implies
that \eqref{3.19} holds. This finishes the proof of Lemma \ref{l3.5}
in Case 1).

{\it Case 2)} $O=\boz$. In this case, we claim that $\boz$ is
bounded. Otherwise, $|\boz|=\fz$. Indeed, if $\boz$ is unbounded,
then $\diam(\boz)=\fz$. By this and Lemma \ref{l3.3}, we know that
for any cube $Q$ with its center $x_Q \in\boz$,
$$|\boz|\ge|Q\cap\boz|\gs|Q|,$$
which, together with the arbitrariness
of $Q$, implies that $|\boz|=\fz$. Moreover, from $f\in
H_{\Phi,\,\cn_h}(\boz)$, we deduce that $\cn_h (f)\in
L^{\Phi}(\boz)$, which, together with Lemma \ref{l3.4}, implies that
$\wz{S}_h^{\epz,\,R,\,1/2}(f)\in L^{\Phi}(\boz)$. By this and the
definition of $O$, we have $|O|<\fz$, which conflicts with
$|O|=|\boz|=\fz$. Thus, the claim holds.

By Lemma \ref{l3.4}, we know that there exists a positive constant
$C_1$ such that for all $R\in(\diam(\boz),d_{\boz})$ and $x\in\boz$,
\begin{equation}\label{3.34}
\wz{S}^{\diam(\boz),\,R,\,1/20}_h (f)(x)\le C_1\cn_h (f)(x).
\end{equation}

Now we continue the proof of Lemma \ref{l3.5} by using \eqref{3.34}.
Without loss of generality, we may assume that $R\ge\diam(\boz)$.
Otherwise, we replace $R$ just by $\diam(\boz)$ in \eqref{3.14}. If
$\gz\ge\frac{1}{C_1}$, then
\begin{eqnarray*}
&&\lfz|\lfz\{x\in\boz:\ \wz{S}^{\epz,\,R,\,1/20}_h (f)(x)>2\lz,\,\cn_h
(f)(x)\le \gz\lz\r\}\r|\\
&&\hs\le|\boz|\le C_1^2\gz^2|O|\ls\gz^2|O|,
\end{eqnarray*}
which shows Lemma \ref{l3.5} in the case that $O=\boz$
and $\gz\ge\frac{1}{C_1}$.

If $\gz<\frac{1}{C_1}$, by the fact that $\cn_h (f)(x)\le\gz\lz$ for
all $x\in F$ and \eqref{3.34}, we have that for any $R\ge
\diam(\boz)$ and $x\in F$,
$$\wz{S}^{\diam(\boz),\,R,\,1/20}_h (f)(x)\le C_1 \cn_h
(f)(x)\le\frac{1}{\gz}\gz\lz=\lz,$$
which implies that
\begin{eqnarray*}
&&\lfz\{x\in\boz:\ \wz{S}^{\epz,\,R,\,1/20}_h (f)(x)>2\lz,\,\cn_h
(f)(x)\le \gz\lz\r\}\\
&&\hs\subset\lfz\{x\in\boz:\ \wz{S}^{\epz,\,\diam(\boz),\,1/20}_h
(f)(x)>\lz,\,\cn_h (f)(x)\le \gz\lz\r\}.
\end{eqnarray*}
Thus, to finish the proof of Lemma \ref{l3.5} in this case, it
suffices to show that
$$\lfz|\lfz\{x\in\boz:\ \wz{S}^{\epz,\,\diam(\boz),\,1/20}_h
(f)(x)>\lz,\,\cn_h (f)(x)\le \gz\lz\r\}\r|\ls\gz^2|O|,$$ whose proof
is similar to that of \eqref{3.18} with $10l_k$ and $Q_k\cap F$
respectively replaced by $\diam(\boz)$ and $\boz$. We omit the
details, which completes the proof of Lemma \ref{l3.5}.
\end{proof}

\begin{lemma}\label{l3.7}
Let $\Phi,\,\boz$ and $L$ be as in Proposition \ref{p3.1}. For all
$\az,\,\bz\in(0,\fz)$, $0\le\epz<R<d_{\boz}$ and all $f\in L^2
(\boz)$,
$$\int_{\boz}\Phi\lfz(\wz{S}^{\epz,\,R,\,\az}_h (f)(x)\r)\,dx\sim
\int_{\boz}\Phi\lfz(\wz{S}^{\epz,\,R,\,\bz}_h (f)(x)\r)\,dx,$$
where the implicit constants are independent of $\epz,\,R$ and $f$.
\end{lemma}

The proof of Lemma \ref{l3.7} is similar to that of
\cite[Proposition 4]{cms85}. We omit the details.

Now we show Proposition \ref{p3.2} by using Lemmas \ref{l3.5} and
\ref{l3.7}.

\begin{proof}[\bf{Proof of Proposition \ref{p3.2}.}]
Let $f\in H_{\Phi,\,\cn_h}(\boz)\cap L^2 (\boz)$. By the upper type
1 and the lower type $p_{\Phi}$ properties of $\Phi$, we know that
$$\Phi(t)\sim\int_0^t \frac{\Phi(s)}{s}\,ds$$
for all $t\in(0,\fz)$.
From this, Fubini's theorem and Lemma \ref{l3.5}, it follows that
for all $\epz,\,R\in(0,d_{\boz})$ with $\epz<R$ and $\gz\in(0,1]$,
\begin{eqnarray}\label{3.35}
&&\int_{\boz}\Phi\lfz(\wz{S}^{\epz,\,R,\,1/20}_h (f)(x)\r)\,dx\\
\nonumber &&\hs\sim\int_{\boz}\int_0^{\wz{S}^{\epz,\,R,\,1/20}_h
(f)(x)}\frac{\Phi(t)}{t}\,dt\,dx
\sim\int_0^{\fz}\frac{\Phi(t)}{t}\sz_{\wz{S}^{\epz,\,R,\,1/20}_h
(f)}(t)\,dt\\ \nonumber
&&\hs\ls\int_0^{\fz}\frac{\Phi(t)}{t}\sz_{\cn_h (f)}(\gz t)\,dt+
\gz^2\int_0^{\fz}\frac{\Phi(t)}{t}\sz_{\wz{S}^{\epz,\,R,\,1/2}_h
(f)}(t/2)\,dt\\ \nonumber
&&\hs\ls\frac{1}{\gz}\int_0^{\fz}\frac{\Phi(t)}{t}\sz_{\cn_h (f)}(
t)\,dt+\gz^2\int_0^{\fz}\frac{\Phi(t)}{t}\sz_{\wz{S}^{\epz,\,R,\,1/2}_h
(f)}(t)\,dt\\ \nonumber
&&\hs\sim\frac{1}{\gz}\int_{\boz}\Phi\lfz(\cn_h
(f)(x)\r)\,dx+\gz^2\int_{\boz}\Phi\lfz(\wz{S}^{\epz,\,R,\,1/2}_h
(f)(x)\r)\,dx,
\end{eqnarray}
where
$$\sz_{\wz{S}^{\epz,\,R,\,1/20}_h (f)}(t):=\lfz|\lfz\{x\in\boz:\
\wz{S}^{\epz,\,R,\,1/20}_h (f)(x)>t\r\}\r|.$$
Furthermore, by Lemma
\ref{l3.7}, \eqref{3.35} and $\wz{S}^{\epz,\,R,\,1/2}_h
(f)\le\wz{S}^{\epz,\,R}_h (f)$, we have that for all
$\epz,\,R\in(0,d_{\boz})$ with $\epz<R$ and $\gz\in(0,1]$,
\begin{eqnarray*}
\int_{\boz}\Phi\lfz(\wz{S}^{\epz,\,R}_h
(f)(x)\r)\,dx&\sim&\int_{\boz}\Phi\lfz(\wz{S}^{\epz,\,R,\,1/20}_h
(f)(x)\r)\,dx\\
&\ls&\frac{1}{\gz}\int_{\boz}\Phi\lfz(\cn_h
(f)(x)\r)\,dx+\gz^2\int_{\boz}\Phi\lfz(\wz{S}^{\epz,\,R}_h
(f)(x)\r)\,dx,
\end{eqnarray*}
which, together with the facts that for all $\lz\in(0,\fz)$,
$$\wz{S}^{\epz,\,R}_h (f/\lz)=\wz{S}^{\epz,\,R}_h (f)/\lz$$
and
$$\cn_h(f/\lz)=\cn_h (f)/\lz,$$
implies that there exists a positive
constant $C_2$ such that
\begin{eqnarray}\label{3.36}
&&\int_{\boz}\Phi\lfz(\frac{\wz{S}^{\epz,\,R}_h (f)(x)}{\lz}\r)\,dx\\
\nonumber &&\hs\le C_2\lfz\{\frac{1}{\gz}\int_{\boz}\Phi\lfz(\frac{\cn_h
(f)(x)}{\lz}\r)\,dx+\gz^2\int_{\boz}\Phi\lfz(\frac{\wz{S}^{\epz,\,R}_h
(f)(x)}{\lz}\r)\,dx\r\}.
\end{eqnarray}
Take $\gz\in(0,1]$ such that $C_2\gz^2=1/2$. Then by \eqref{3.36},
we obtain that for all $\lz\in(0,\fz)$,
\begin{equation*}
\int_{\boz}\Phi\lfz(\frac{\wz{S}^{\epz,\,R}_h
(f)(x)}{\lz}\r)\,dx\ls\int_{\boz}\Phi\lfz(\frac{\cn_h
(f)(x)}{\lz}\r)\,dx.
\end{equation*}
By the Fatou lemma and letting $\epz\to0$ and $R\to d_{\boz}$, we
obtain that for any $\lz\in(0,\fz)$,
$$\int_{\boz}\Phi\lfz(\frac{\wz{S}_h
(f)(x)}{\lz}\r)\,dx\ls\int_{\boz}\Phi\lfz(\frac{\cn_h
(f)(x)}{\lz}\r)\,dx,$$
which implies that
$$\|\wz{S}_h(f)\|_{L^{\Phi}(\boz)}\ls\|\cn_h (f)\|_{L^{\Phi}(\boz)}.$$
This finishes the proof of Proposition \ref{p3.2}.
\end{proof}

\begin{proposition}\label{p3.3}
Let $\Phi,\,\boz$ and $L$ be as in Proposition \ref{p3.1}. Then
under DBC, there exists a positive constant $C$ such that for all
$f\in L^2 (\boz)$,
$$\|S_h (f)\|_{L^{\Phi}(\boz)}\le C\lfz\|\wz{S}_h
(f)\r\|_{L^{\Phi}(\boz)}.$$
\end{proposition}

\begin{proof}[\bf{Proof.}]
To show this proposition, we borrow some ideas from \cite{hm09}. Fix
$\epz,\,R\in(0,d_{\boz})$ with $\epz<R$ and $x\in\boz$. Let $f\in
L^2 (\boz)$ and, for $\az\in(0,\fz)$,
$$\wz{\wz{\bgz}}^{\epz,\,R}_{\az}(x):=\{(y,t)\in\rn\times(\epz,R):\
|x-y|<\az t\}.$$
Take $\eta\in C^{\fz}_{c}(\rn\times(0,\fz))$ such
that $\eta\equiv1$ on $\wz{\wz{\bgz}}^{\epz,\,R}_1 (x)$, $0\le\eta\le1$,
$$\supp\eta\subset\wz{\wz{\bgz}}^{\epz/2,\,2R}_{3/2}(x)$$
and for all
$(y,t)\in\wz{\wz{\bgz}}^{\epz/2,\,2R}_{3/2} (x)$,
$|\nabla\eta(y,t)|\ls\frac{1}{t}$. By the choice of $\eta$, we have
that for all $t\in(\epz/2,2R)$, $\eta_t (\cdot):=\eta(\cdot,t)\in
C^{\fz}(\boz)$. In the rest part of this proof, we denote $e^{-t^2
L}(f)$ by $u_t$ for all $t\in(0,d_{\boz})$. Then by \cite[p.\,23,
(1.19)]{o04}, we know that for any given $t\in(0,d_{\boz})$, $u_t\in
D(L)\subset W^{1,\,2}_0 (\boz)$. Moreover, by the fact that for all
$t\in(0,d_{\boz})$,
$$Lu_t=e^{-\frac{t^2}{2}L}\lfz(Le^{-\frac{t^2}{2}L}(f)\r),$$
and \cite[p.\,23, (1.19)]{o04} again,
we have that $Lu_t\in D(L)\subset W^{1,\,2}_0 (\boz)$, which,
together with $\eta_t \in C^{\fz}(\boz)$, implies that for all
$t\in(0,d_{\boz})$, $(Lu_t)\eta_t \in W^{1,\,2}_0 (\boz)$. From
this, \eqref{1.2}, the facts that $0\le\eta\le1$ and $\eta\equiv1$ on
$\wz{\wz{\bgz}}^{\epz,\,R}_1 (x)$, and H\"older's inequality,  we
deduce that
\begin{eqnarray}\label{3.37}
\hs\hs S^{\epz,\,R}_h (f)(x)&=&\lfz\{\int_{\bgz^{\epz,\,R}_1 (x)}\lfz|t^2
Le^{-t^2 L}(f)(y)\r|^2\frac{dy\,dt}{t|Q(x,t)\cap\boz|}\r\}^{1/2}\\
\nonumber &\le&\Bigg\{\int_{\bgz^{\epz/2,\,2R}_{3/2} (x)} t^2
Le^{-t^2 L}(f)(y)\\ \nonumber &&\times\ol{t^2 Le^{-t^2
L}(f)(y)}\eta(y,t)\frac{dy\,dt}{t|Q(x,t)\cap\boz|}\Bigg\}^{1/2}\\
\nonumber &\le&\Bigg\{\int_{\bgz^{\epz/2,\,2R}_{3/2} (x)} t\lfz|
A(y)\nabla u_t (y)\cdot t\nabla \ol{(t^2 L
u_t)(y)}\r|\\ \nonumber
&&\times\eta(y,t)\frac{dy\,dt}{t|Q(x,t)\cap\boz|}\Bigg\}^{1/2}\\
\nonumber &&+\Bigg\{\int_{\bgz^{\epz/2,\,2R}_{3/2} (x)}
t\lfz|A(y)\nabla u_t (y)\cdot\nabla \eta(y,t) \ol{t^3 L u_t(y)}\r|\\
\nonumber &&\times\frac{dy\,dt}{t|Q(x,t)\cap\boz|}\Bigg\}^{1/2}\\
\nonumber &\ls&\lfz\{\int_{\bgz^{\epz/2,\,2R}_{3/2} (x)} |t\nabla u_t
(y)|^2\frac{dy\,dt}{t|Q(x,t)\cap\boz|}\r\}^{1/4}\\ \nonumber
&&\times\lfz\{\int_{\bgz^{\epz/2,\,2R}_{3/2} (x)} \lfz|t\nabla (t^2 L
u_t)(y)\r|^2\frac{dy\,dt}{t|Q(x,t)\cap\boz|}\r\}^{1/4}\\
\nonumber &&+\lfz\{\int_{\bgz^{\epz/2,\,2R}_{3/2} (x)} |t\nabla u_t
(y)|^2\frac{dy\,dt}{t|Q(x,t)\cap\boz|}\r\}^{1/4}\\
\nonumber &&\times\lfz\{\int_{\bgz^{\epz/2,\,2R}_{3/2} (x)} |t^2 L
u_t(y)|^2\frac{dy\,dt}{t|Q(x,t)\cap\boz|}\r\}^{1/4}.
\end{eqnarray}
For all $(z,\tau)\in\bgz^{\epz/2,\,2R}_{3/2}(x)$, let
$$E_{(z,\,\tau)}:= B((z,\tau),\gz\tau)\cap(\boz\times(0,d_{\boz})),$$
where $\gz$ is a positive constant which is determined later.
From the Besicovitch covering lemma, it follows that
there exists a subcollection
$\{E_{(z_j,\tau_j)}\}_j$ of
$\{E_{(z,\tau)}\}_{(z,\tau)\in\bgz^{\epz/2,\,2R}_{3/2}(x)}$ such
that \eqref{3.12} holds in this case. For each $j$, we denote
$E_{(z_j,\tau_j)}$ simply by $E_j$. Similarly to the facts (i) and
(ii) appearing in the proof of Lemma \ref{l3.4}, we have the
following two facts for $E_j$:

\begin{enumerate}
  \item[(i)] For each $j$, if $(y,t)\in E_j$, then $t\sim d_j\sim r(E_j)$,
where $d_j$ and $r(E_j)$ denote, respectively, the distance from
$E_j$ to the bottom boundary $\boz\times\{0\}$ and the radius of
$E_j$.
  \item[(ii)] For each $j$, let
$$\wz{E}_j:=B((z_j,\tau_j),9\gz\tau_j)\cap(\boz\times(0,d_{\boz})).$$
If $\gz\in(0,1/54)$, then $\wz{E}_j\subset\bgz^{\epz/4,\,4R}_2 (x)$.
\end{enumerate}

For all $t\in(0,d_{\boz})$, let $v_t:= Le^{-tL}(f)$. Then we have
that
$$\paz_t v_t +Lv_t=0.$$
Thus, from Remark \ref{r3.1} (in which,
if $\tau_j\in(\epz,\frac{d_{\boz}}{1+\gz}]$, we choose
$$t_0:=(1+\gz)\tau_j,$$
$r:=\gz\tau_j$ and $c:=\frac{4}{\gz}$, and if
$\tau_j\in(\frac{d_{\boz}}{1+\gz},d_{\boz})$, we then choose $t_0:=
d_{\boz},\,r:=\gz\tau_j$ and $c:=\frac{4}{\gz(1+\gz)^2}$, and in
both cases, we need choose $\gz\in(0,1/54)$), we deduce that for
each $j$,
$$\int_{E_j}t|\nabla(Lu_t)(y)|^2\,dy\,dt\ls\frac{1}{[r(E_j)]^2}
\int_{\wz{E}_j}t|Lu_t
(y)|^2\,dy\,dt.$$
By this, the above facts (i) and (ii), and \eqref{3.12},
we obtain that
\begin{eqnarray*}
&&\int_{\bgz^{\epz/2,\,2R}_{3/2}(x)}\lfz|t\nabla (t^2
Lu_t)(y)\r|^2\frac{dy\,dt}{t|Q(x,t)\cap\boz|}\\
&&\hs\le\sum_j \int_{E_j} \lfz|t\nabla (t^2
Lu_t)(y)\r|^2\frac{dy\,dt}{t|Q(x,t)\cap\boz|}\\
&&\hs\sim\sum_j \frac{[r(E_j)]^4}{|Q(x,r(E_j))\cap\boz|}\int_{E_j}t
\lfz|\nabla (Lu_t)(y)\r|^2\,dy\,dt\\
&&\hs\ls\sum_j
\frac{[r(E_j)]^2}{|Q(x,r(E_j))\cap\boz|}\int_{\wz{E}_j}t
\lfz|Lu_t(y)\r|^2\,dy\,dt\\
&&\hs\sim\sum_j \frac{1}{[r(E_j)]^2}\int_{\wz{E}_j} \lfz|t^3 Lu_t
(y)\r|^2\frac{dy\,dt}{t|Q(x,t)\cap\boz|}\\
&&\hs\ls\int_{\bgz^{\epz/4,\,4R}_{2}(x)}\lfz|t^2
Lu_t(y)\r|^2\frac{dy\,dt}{t|Q(x,t)\cap\boz|},
\end{eqnarray*}
which, together with \eqref{3.37}, implies that
$$S^{\epz,\,R}_h
(f)(x)\ls\lfz[\wz{S}^{\epz/2,\,2R,\,3/2}_h
(f)(x)\r]^{1/2}\lfz[S^{\epz/4,\,4R,\,2}_h (f)(x)\r]^{1/2}.$$
By the Fatou lemma, letting $\epz\to0$ and $R\to d_{\boz}$, we have
that
$$S_h (f)(x)\ls[\wz{S}^{3/2}_h (f)(x)]^{1/2}[S^{2}_h (f)(x)]^{1/2},$$
which, together with Cauchy's inequality, implies that there exists
a positive constant $C_2$ such that for all $\uc\in(0,1)$,
\begin{eqnarray}\label{3.38}
S_h (f)(x)\le\frac{C_2}{\uc}\wz{S}^{3/2}_h (f)(x)+\uc S^{2}_h
(f)(x).
\end{eqnarray}
Similarly to the proof of Lemma \ref{l3.7}, we have that there
exists a positive constant $C_3$ such that for all $g\in L^2
(\boz)$,
$$\int_{\boz}\Phi\lfz(S^{2}_h (g)(y)\r)\,dy\le C_3
\int_{\boz}\Phi\lfz(S_h (g)(y)\r)\,dy.$$
From this, \eqref{3.38}, the
strictly lower type $p_{\Phi}$ and the upper type 1 properties of
$\Phi$, it follows that there exists a positive constant $C$ such
that for all $x\in\boz$,
\begin{eqnarray}\label{3.39}
\hs\int_{\boz}\Phi\lfz(S_h
(f)(x)\r)\,dx&\le&\int_{\boz}\Phi\lfz(\frac{C_2}{\uc}\wz{S}^{3/2}_h
(f)(x)\r)\,dx\\ \nonumber
&&+\int_{\boz}\Phi\lfz(S^2_h (f)(x)\r)\,dx\\
\nonumber &\le&\frac{CC_2}{\uc}\int_{\boz}\Phi\lfz(\wz{S}^{3/2}_h
(f)(x)\r)\,dx\\ \nonumber &&+C_3 \uc^{p_{\Phi}}\int_{\boz}\Phi\lfz(S_h
(f)(x)\r)\,dx.
\end{eqnarray}
Take $\uc\in(0,1)$ small enough such that $C_3
\uc^{p_{\Phi}}\le\frac{1}{2}$. By this, \eqref{3.39} and Lemma
\ref{l3.7}, we obtain that
$$\int_{\boz}\Phi(S_h
(f)(x))\,dx\ls\int_{\boz}\Phi\lfz(\wz{S}_h (f)(x)\r)\,dx,$$
which, together
with the facts that for all $\lz\in(0,\fz)$,
$$S_h (f/\lz)=S_h(f)/\lz\ \ {\rm and}\ \ \wz{S}_h (f/\lz)=\wz{S}_h (f)/\lz,$$
implies that
$$\int_{\boz}\Phi\lfz(\frac{S_h
(f)(x)}{\lz}\r)\,dx\ls\int_{\boz}\Phi\lfz(\frac{\wz{S}_h
(f)(x)}{\lz}\r)\,dx.$$ From this, it follows that Proposition
\ref{p3.3} holds, which completes the proof of Proposition
\ref{p3.3}.
\end{proof}

To complete the proof of Theorem \ref{t1.1}, we need the following
key proposition.

\begin{proposition}\label{p3.4}
Let $\Phi,\,\boz$ and $L$ be as in Theorem \ref{t1.1}. Assume that
$L$ satisfies DBC and the semigroup generated by $L$ has the
Gaussian property $(G_{\diam(\boz)})$.
\begin{enumerate}
  \item[{\rm(i)}] If $\boz$ is unbounded, then
$$(H_{\Phi,\,S_h}(\boz)\cap
L^2(\boz))\subset (H_{\Phi,\,r}(\boz)\cap L^2(\boz))$$
and there is a positive constant $C$ such that for all $f\in
H_{\Phi,\,S_h}(\boz)\cap L^2 (\boz)$,
$$\|f\|_{H_{\Phi,\,r}(\boz)}\le
C\|f\|_{H_{\Phi,\,S_h}(\boz)}.$$
  \item[{\rm(ii)}] If $\boz$ is bounded, then
$$(H_{\Phi,\,S_h\,d_{\boz}}(\boz)\cap L^2(\boz))\subset
(H_{\Phi,\,r}(\boz)\cap L^2(\boz))$$
and there is a positive
constant $C$ such that for all $f\in
H_{\Phi,\,S_h,\,d_{\boz}}(\boz)\cap L^2 (\boz)$,
$$\|f\|_{H_{\Phi,\,r}(\boz)}\le
C\|f\|_{H_{\Phi,\,S_h,\,d_{\boz}}(\boz)}.$$
Moreover, if, in
addition, $n\ge3$ and $(G_{\fz})$ holds, then
\begin{eqnarray*}
(H_{\Phi,\,\wz{S}_h,\,d_{\boz}}(\boz)\cap L^2(\boz))
&&=(H_{\Phi,\,\wz{S}_h}(\boz)\cap L^2(\boz))\\
&&=(H_{\Phi,\,S_h,\,d_{\boz}}(\boz)\cap L^2(\boz))\\
&&=(H_{\Phi,\,S_h}(\boz)\cap L^2 (\boz))
\end{eqnarray*}
with equivalent norms.
\end{enumerate}
\end{proposition}

To show Proposition \ref{p3.4}, we need the atomic decomposition of
the tent space on $\boz$. Now we recall some definitions and notion
about the tent space, which was initially introduced by Coifman,
Meyer and Stein \cite{cms85} on $\rn$, and then generalized by Russ
\cite{ru07} to spaces of homogeneous type in the sense of Coifman
and Weiss \cite{cw71,cw77}. Recall that it is well known that the
strongly Lipschitz domain $\boz$ is a space of homogeneous type. For
all measurable functions $g$ on $\boz\times(0,\fz)$ and $x\in\boz$,
define
$$\ca(g)(x):=\lfz\{\int_{\wz{\bgz}(x)}|g(x,t)|^2
\frac{dy}{|Q(x,t)\cap\boz|}\frac{dt}{t}\r\}^{1/2},$$
where
$$\wz{\bgz}(x):=\{(y,t)\in\boz\times(0,\fz):\ |y-x|<t\}.$$
In what follows, we denote by $T_{\Phi}(\boz)$ the {\it space of all
measurable functions $g$ on $\boz\times(0,\fz)$ such that $\ca(g)\in
L^{\Phi}(\boz)$} and for any $g\in T_{\Phi}(\boz)$, define its {\it
quasi-norm} by
$$\|g\|_{T_{\Phi}(\boz)}:=\|\ca(g)\|_{L^{\Phi}(\boz)}:=\inf\lfz
\{\lz\in (0,\fz):\ \int_{\boz}\Phi\lfz(\frac{\ca(g)(x)}{\lz}\r)\,dx\le1\r\}.$$
When $\bfai(t):= t$ for all $t\in(0,\fz)$, we denote
$T_{\bfai}(\boz)$ simply by $T_1(\boz)$.

A function $a$ on $\boz\times(0,\fz)$ is called a {\it
$T_{\Phi}(\boz)$-atom} if
\begin{enumerate}
  \item[(i)] there exists a cube
$$Q:= Q(x_{Q},l(Q))\subset\rn$$
with $x_Q\in\boz$ and $l(Q)\in(0,\fz)\cap(0,d_{\boz}]$ such that $\supp
a\subset\wt{Q\cap\boz}$, where and in what follows,
$$\wt{Q\cap\boz}:=\lfz\{(y,t)\in\boz\times(0,\fz):\
|y-x_Q|<\frac{l(Q)}{2}-t\r\};$$
  \item[(ii)]
$$\int_{\wt{Q\cap\boz}}|a(y,t)|^2\frac{dy\,dt}{t}\le|Q\cap\boz|^{-1}
[\rz(|Q\cap\boz|)]^{-2}.$$
\end{enumerate}

Since $\Phi$ is of upper type 1, it is easy to see that there exists
a positive constant $C$ such that for all $T_{\Phi}(\boz)$-atoms
$a$, we have $\|a\|_{T_{\Phi}(\boz)}\le C$; see \cite{jy}. By a
slight modification on the proof of \cite[Theorem 3.1]{jy}, we have
the following atomic decomposition for functions in
$T_{\Phi}(\boz)$. We omit the details.
\begin{lemma}\label{l3.8}
Let $\boz$ be a strongly Lipschitz domain of $\rn$ and $\Phi$
satisfy Assumption (A). Then for any $f\in T_{\Phi}(\boz)$, there
exist a sequence $\{a_j\}_j$ of $T_{\Phi}(\boz)$-atoms  and a
sequence  $\{\lz_j\}_j $ of numbers such that for almost every
$(x,t)\in\boz\times(0,\fz)$,
\begin{equation}\label{3.40}
f(x,t)=\sum_j \lz_j a_j (x,t).
\end{equation}
Moreover, there exists a positive constant $C$ such that for all
$f\in T_{\Phi}(\boz)$,
\begin{eqnarray}\label{3.41}
\quad&&\blz(\{\lz_j a_j\}_j)\\ 
\nonumber &&\hs:=\!\inf\lfz\{\lz\in (0,\fz):\ \sum_j
|Q_j \cap\boz|\Phi\lfz(\frac{|\lz_j|\|a_j\|_{T^2_2
(\boz\times(0,\fz))}}{\lz|Q_j\cap\boz|^{1/2}}\r)\le1\r\}\\
\nonumber &&\hs\le C\|f\|_{T_{\Phi}(\boz)},
\end{eqnarray}
where $Q_j\cap\boz$ appears in the support of $a_j$ and
$$\|a_j\|_{T^2_2 (\boz\times(0,\fz))}:=\lfz\{\int_{\wt{Q_j\cap\boz}}
|a_j (y,t)|^2\frac{dy\,dt}{t}\r\}^{1/2}.$$
\end{lemma}

In \cite[p.\,183]{ar03}, Auscher and Russ showed the following
property of strongly Lipschitz domains, which plays an important
role in the proof of Proposition \ref{p3.4}.

\begin{lemma}\label{l3.9}
Let $\boz$ be a strongly  Lipschitz domain of $\rn$. Then there
exists $\rz(\boz)\in(0,\fz)$ such that for any cube $Q$ satisfying
$l(Q)<\rz(\boz)$ and $2Q\subset\boz$ but
$4Q\cap\paz\boz\neq\emptyset$, where $\paz\boz$ denotes the boundary
of $\boz$, there exists a cube $\wz{Q}\subset\boz^\complement$ such
that $l(\wz{Q})=l(Q)$ and the distance from $\wz{Q}$ to $Q$ is
comparable to $l(Q)$. Furthermore, $\rz(\boz)=\fz$ if
$\boz^\complement$ is unbounded.
\end{lemma}

Now we show Proposition \ref{p3.4} by applying Lemmas \ref{l3.1},
\ref{l3.8} and \ref{l3.9}.

\begin{proof}[\bf{Proof of Proposition \ref{p3.4}.}]
We first prove  Proposition \ref{p3.4}(i) by borrowing some ideas
from the proof of \cite[p.\,594,\,Theorem C]{cw77} (see also
\cite{hyz09} and \cite{jyz09}). Recall that in this case, since
$\boz$ is unbounded, we have $\diam(\boz)=\fz$. Let $f\in
H_{\Phi,\,S_h}(\boz)\cap L^2 (\boz)$. Then by the
$H^{\fz}$-functional calculus for $L$,  we know that
\begin{eqnarray}\label{3.42}
f=8\int_0^{\fz}(t^2 Le^{-t^2 L})(t^2 Le^{-t^2 L}) (f)\frac{dt}{t}
\end{eqnarray}
in $L^2 (\boz)$; see also \cite[(9)]{h}. Since $f\in
H_{\Phi,\,S_h}(\boz)$, we have that $S_h (f)\in L^{\Phi}(\boz)$,
which implies that $t^2 Le^{-t^2 L}(f)\in T_{\Phi}(\boz)$ and
$$\lfz\|f\r\|_{H_{\Phi,\,S_h}(\boz)}=\lfz\|t^2 Le^{-t^2 L}
(f)\r\|_{T_{\Phi}(\boz)}.$$
Then from Lemma \ref{l3.8}, we deduce
that there exist $\{\lz_j\}_j \subset\cc$ and a sequence $\{a_j\}_j$
of $T_{\Phi}(\boz)$-atoms such that for almost every
$(x,t)\in\boz\times(0,\fz)$,
\begin{eqnarray}\label{3.43}
t^2 Le^{-t^2 L}(f)(x)=\sum_j\lz_j a_j(x,t).
\end{eqnarray}
For each $j$, let
$$\az_j:=8\int_0^{\fz}t^2 Le^{-t^2 L}
(a_j)\frac{dt}{t}.$$
Then by \eqref{3.42} and \eqref{3.43}, similarly
to the proof of \cite[Proposition 4.2]{jy10}, we have that
\begin{equation}\label{3.44}
f=\sum_j\lz_j \az_j
\end{equation}
in $L^2 (\boz)$. For any $T_{\Phi}(\boz)$-atom $a$ supported in
$\wt{Q\cap\boz}$, let
\begin{eqnarray}\label{3.45}
\az:=8\int_0^{\fz}t^2 Le^{-t^2 L}(a)\frac{dt}{t}.
\end{eqnarray} To show Proposition
\ref{p3.4}, it suffices to show that there exist a function
$\wz{\az}$ on $\rn$ such that
\begin{equation}\label{3.46}
\wz{\az}|_{\boz}=\az
\end{equation}
and a sequence $\{b_i\}_i$ of harmless constant multiples of
$(\rz,\,2,\,0)$-atoms, with the constant depending on $i$,  such
that $\wz{\az}=\sum_i b_i$ in $L^2 (\rn)$ and
\begin{equation}\label{3.47}
\sum_i |Q_i|\Phi\lfz(\frac{\|b_i\|_{L^2
(\rn)}}{|Q_i|^{1/2}}\r)\ls|Q\cap\boz|\Phi\lfz(\frac{\|a\|_{T^2_2
(\boz\times(0,\fz))}}{|Q\cap\boz|^{1/2}}\r),
\end{equation}
where for each $i$, $\supp b_i\subset Q_i$ and  $Q\cap\boz$ appears
in the support of $a$. Indeed, if \eqref{3.46} and \eqref{3.47}
hold, then by \eqref{3.46}, we know that for each $j$, there exists
a function $\wz{\az}_j$ on $\rn$ such that
$\wz{\az}_j|_{\boz}=\az_j$. Let
$$\wz{f}:=\sum_j \lz_j \wz{\az}_j.$$
Then $\wz{f}|_{\boz}=f$. Furthermore, from \eqref{3.47}, we deduce
that there exists a sequence $\{b_{j,\,i}\}_{j,\,i}$ of harmless
constant multiples of $(\rz,\,2,\,0)$-atoms, with the constant
depending on $j$ and $i$, such that
$$\wz{f}=\sum_j \sum_{i} \lz_j
b_{j,\,i}$$
and
\begin{equation*}
\sum_{j,\,i} |Q_{j,\,i}|\Phi\lfz(\frac{|\lz_j|\|b_{j,\,i}\|_{L^2
(\rn)}}{|Q_{j,\,i}|^{1/2}}\r)\ls\sum_j |Q_j
\cap\boz|\Phi\lfz(\frac{|\lz_j|\|a_j\|_{T^2_2
(\boz\times(0,\fz))}}{|Q_j \cap\boz|^{1/2}}\r),
\end{equation*}
where for each $j$ and $i$, $\supp b_{j,\,i}\subset Q_{j,\,i}$ and
$Q_j \cap\boz$ appears in the support of $a_j$, which, together with
the facts that for all $\lz\in(0,\fz)$,
$$\|b_{i,\,j}/\lz\|_{L^2
(\rn)}=\|b_{i,\,j}\|_{L^2 (\rn)}/\lz$$
and
$$\|a_j /\lz\|_{T^2_2
(\boz\times(0,\fz))}=\|a_j\|_{T^2_2 (\boz\times(0,\fz))}/\lz,$$
implies that for all $\lz\in(0,\fz)$,
\begin{equation*}
\sum_{j,\,i} |Q_{j,\,i}|\Phi\lfz(\frac{|\lz_j|\|b_{j,\,i}\|_{L^2
(\rn)}}{\lz|Q_{j,\,i}|^{1/2}}\r)\ls\sum_j|Q_j
\cap\boz|\Phi\lfz(\frac{|\lz_j|\|a_j\|_{T^2_2
(\boz\times(0,\fz))}}{\lz|Q_j \cap\boz|}\r).
\end{equation*}
From this and Lemmas \ref{l3.1} and \ref{l3.8}, it follows that
$\wz{f}\in H_{\Phi}(\rn)$ and
\begin{equation}\label{3.48}
\lfz\|\wz{f}\r\|_{H_{\Phi}(\rn)}\sim\lfz\|\wz{f}\r\|_{H^{\rz,\,2,\,0}(\rn)}
\ls\lfz\|f\r\|_{H_{\Phi,\,S_h}(\boz)}.
\end{equation}
Thus, $f\in H_{\Phi,\,r}(\boz)$ and
$$\|f\|_{H_{\Phi,\,r}(\boz)}\ls\|f\|_{H_{\Phi,\,S_h}(\boz)},$$
which, together with the arbitrariness
of $f\in H_{\Phi,\,S_h}(\boz)\cap
L^2 (\boz)$, implies that the desired conclusion of Proposition
\ref{p3.4}(i).

Let $Q:= Q(x_0,r_0)$. Now we show \eqref{3.46} and \eqref{3.47} by
considering the following two cases for $Q$ which appears in the
support of $a$.

{\it Case 1)} $8Q\cap\boz^\complement\neq\emptyset$.
In this case, let
$$R_k (Q) :=(2^{k+1}Q\setminus2^{k}Q)\cap\boz$$
if $k\ge3$ and $R_0 (Q):=8Q\cap\boz$. Let
$$J_{\boz}:=\lfz\{k\in\nn:\ k\ge3,\,|R_k
(Q)|>0\r\}.$$
For $k\in J_{\boz}\cup\{0\}$, let $\chi_k:=\chi_{R_k
(Q)}$, $\wz{\chi}_k := |R_k (Q)|^{-1}\chi_k$ and
$$m_k :=\int_{R_k(Q)}\az(x)\,dx.$$
Then we have
\begin{equation}\label{3.49}
\az=\az\chi_0+\sum_{k\in J_{\boz}}\az\chi_k
\end{equation}
almost everywhere and also in $L^2 (\boz)$. Take the cube
$\wz{Q}\subset\rn$ such that the center $x_{\wz{Q}}$ of $\wz{Q}$
satisfying that $x_{\wz{Q}}\in\boz^\complement$, $l(\wz{Q})=l(Q)$ and
$\dist(Q,\wz{Q})\sim l(Q)$. Then there exists a cube $Q^{\ast}_0$
such that $(Q\cup\wz{Q})\subset Q^{\ast}_0$ and

\begin{equation}\label{3.50}
l(Q^{\ast}_0)\sim l(Q).
\end{equation}
Let
$$b_0:= \az\chi_0
-\frac{1}{|\wz{Q}\cap\boz^\complement|}\lfz\{\int_{R_0
(Q)}\az(x)\,dx\r\} \chi_{\wz{Q}\cap\boz^\complement}.$$
Then $\int_{\rn}b_0 (x)\,dx=0$ and $\supp b_0 \subset Q^{\ast}_0$.
Similarly to the proof of \cite[(3.36)]{yys}, we have
\begin{equation}\label{3.51}
\|\az\|_{L^2 (\boz)}\ls\|a\|_{T^2_2 (\boz\times(0,\fz))}.
\end{equation}
By the facts that $\boz^\complement$ is an unbounded strongly
Lipschitz domain and Lemma \ref{l3.3}, we know that $|\wz{Q}\cap
\boz^\complement|\sim|\wz{Q}|$. From this, H\"older's inequality,
\eqref{3.50} and \eqref{3.51}, we deduce that
\begin{eqnarray*}
\|b_0\|_{L^2 (\rn)}&\le&\|\az\|_{L^2
(\boz)}+\frac{1}{|\wz{Q}\cap\boz^\complement|^{1/2}}\lfz\{\int_{R_0
(\boz)}
|\az(x)|^2\,dx\r\}^{1/2}|Q\cap\boz|^{1/2}\\
&\ls&\|\az\|_{L^2 (\boz)}+\|\az\|_{L^2
(\boz)}\frac{|Q|^{1/2}}{|\wz{Q}|^{1/2}}\ls\|a\|_{T^2_2
(\boz\times(0,\fz))}\\
&\ls&\frac{1}{|Q\cap\boz|^{1/2}\rz(|Q\cap\boz|)}\sim
\frac{1}{|Q|^{1/2}\rz(|Q|)}\sim
\frac{1}{|Q^{\ast}_0|^{1/2}\rz(|Q^{\ast}_0|)}.
\end{eqnarray*}
Thus, we know that $b_0$ is a harmless constant multiple of a
$(\rz,\,2,\,0)$-atom and, by the upper 1 property of $\bfai$,
\begin{eqnarray}\label{3.52}
|Q^{\ast}_0|\Phi\lfz(\frac{\|b_0\|_{L^2
(\rn)}}{|Q^{\ast}_0|^{1/2}}\r)&&\ls|Q|\Phi\lfz(\frac{\|a\|_{T^2_2
(\boz\times(0,\fz))}}{|Q|^{1/2}}\r)\\
\nonumber &&\ls|Q\cap\boz|\Phi\lfz(\frac{\|a\|_{T^2_2
(\boz\times(0,\fz))}}{|Q\cap\boz|^{1/2}}\r).
\end{eqnarray}

To finish the proof in this case, we need the following Fact 1,
whose proof is similar to the usual Whitney decomposition of an open
set in $\rn$; see, for example, \cite{s93}. We omit the details.

\smallskip

{\bf Fact 1.} {\it For all $k\in J_{\boz}$, there exists the Whitney
decomposition $\{Q_{k,\,i}\}_i$ of $R_k (Q)$ about $\paz\boz$, where
$\{Q_{k,\,i}\}_i$ are dyadic cubes of $\rn$ with disjoint interiors,
and for each $i$, $2Q_{k,\,i}\subset\boz$ but
$4Q_{k,\,i}\cap\paz\boz\neq\emptyset$.}

\smallskip

Notice that Fact 1 was also used in \cite[pp.\,304-305]{cks93} and
\cite[p.\,167]{ar03}. Let $\{Q_{k,\,i}\}_{k\in J_{\boz},\,i}$ be as
in Fact 1. Then for each $k\in J_{\boz}$,
$$\az\chi_{R_k (Q)}=\sum_i\az\chi_{Q_{k,\,i}}$$
almost everywhere. In what follows, for all
$t\in(0,\fz)$, let
$$D_t:= s\paz_s K_s|_{s=t^2}.$$
Then for all $x\in
R_k (Q)$, by \eqref{3.45}, Lemmas \ref{l2.1} and \ref{l3.3}, and
H\"older's inequality, we have that
\begin{eqnarray}\label{3.53}
\lfz|\az(x)\r|&\ls&\int_0^{r_0}\int_{Q\cap\boz}|D_t (x,y)|
|a(y,t)|\,\frac{dy\,dt}{t}\\ \nonumber
&\ls&\int_0^{r_0}\int_{Q\cap\boz}\frac{e^{-\az\frac{|x-y|^2}{t^2}}}{t^{n}}
|a(y,t)|\frac{dy\,dt}{t}\\ \nonumber &\ls&\|a\|_{T^2_2
(\boz\times(0,\fz))} \lfz\{\int_0^{r_0}\int_{Q\cap\boz}\frac{t^2}{|x-y|
^{2(n+1)}}\frac{dy\,dt}{t}\r\}^{1/2}\\ \nonumber
&\ls&|x-x_0|^{-(n+1)}r_0 |Q\cap\boz|^{1/2} \|a\|_{T^2_2
(\boz\times(0,\fz))}\\ \nonumber &\ls&2^{-k(n+1)}|Q\cap\boz|^{-1/2}
\|a\|_{T^2_2 (\boz\times(0,\fz))}.
\end{eqnarray}
Moreover, by Lemma \ref{l3.9}, we know that for each $k$ and $i$,
there exists a cube $\wz{Q}_{k,\,i}\subset\boz^\complement$ such
that $l(\wz{Q}_{k,\,i})=l(Q_{k,\,i})$ and
$\dist(\wz{Q}_{k,\,i},Q_{k,\,i})\sim l(Q_{k,\,i})$. Then for each
$k$ and $i$, there exists a cube $Q^{\ast}_{k,\,i}$ such that
$(Q_{k,\,i}\cup\wz{Q}_{k,\,i})\subset Q^{\ast}_{k,\,i}$ and
$l(Q^{\ast}_{k,\,i})\sim l(Q_{k,\,i})$. For each $k$ and $i$, let
$$b_{k,\,i}:=\az\chi_{Q_{k,\,i}}-\frac{1}{|\wz{Q}_{k,\,i}|}\lfz\{
\int_{Q_{k,\,i}}\az(x)\,dx\r\}\chi_{\wz{Q}_{k,\,i}}.$$
Then
$$\int_{\rn}b_{k,\,i}(x)\,dx=0$$
and $\supp b_{k,\,i}\subset
Q^{\ast}_{k,\,i}$. Furthermore, by \eqref{3.53} and H\"older's
inequality, we have that
\begin{eqnarray}\label{3.54}
\|b_{k,\,i}\|_{L^2 (\rn)}&&\ls \|\az\|_{L^2 (Q_{k,\,i})}\\
\nonumber
&&\ls2^{-k(n+1)}|Q\cap\boz|^{-1/2}|Q^{\ast}_{k,\,i}|
^{1/2}\|a\|_{T^2_2 (\boz\times(0,\fz))}.
\end{eqnarray}
Thus, for each $k$ and $i$, $b_{k,\,i}$ is a constant multiple of
some $(\rz,\,2,\,0)$-atom with the constant depending on $k$ and
$i$. Let
$$\wz{\az}:= b_0 +\sum_{k\in J_{\boz}}\sum_{i}b_{k,\,i}.$$
Then by the constructions of $b_0$ and $\{b_{k,\,i}\}_{k\in
J_{\boz},\,i}$, we know that $\wz{\az}|_{\boz}=\az$. Moreover, we
claim that $\sum_{k\in J_{\boz}}\sum_{i}b_{k,\,i}$ converges in $L^2
(\rn)$. Indeed, let $M$ denote the usual Hardy-Littlewood maximal
operator. Then by \eqref{3.53}, the boundedness of the vector-valued
Hardy-Littlewood maximal operator established by Fefferman and Stein
in \cite[Theorem 1(1)]{fs71}, and the disjoint property of
$\{Q_{k,\,i}\}_i$, we have that for each $k\in J_{\boz}$,
\begin{eqnarray*}
&&\int_{\rn}\lfz[\sum_{i}b_{k,\,i}(x)\r]^2\,dx\\
&&\hs\le
\int_{\rn}\lfz[\sum_{i}2^{-k(n+1)}|Q\cap\boz|^{-1/2}\|a\|_{T^2_2
(\boz\times(0,\fz))}\chi_{Q^{\ast}_{k,\,i}}(x)\r]^2\,dx\\
&&\hs\ls2^{-2k(n+1)}|Q\cap\boz|^{-1}\|a\|^2_{T^2_2
(\boz\times(0,\fz))}\int_{\rn}\lfz\{\sum_i
\lfz[M(\chi_{Q_{k,\,i}})(x)\r]^2\r\}^2\,dx\\
&&\hs\ls2^{-2k(n+1)}|Q\cap\boz|^{-1}\|a\|^2_{T^2_2
(\boz\times(0,\fz))}\int_{\rn}\lfz\{\sum_i
\lfz[\chi_{Q_{k,\,i}}(x)\r]^2\r\}^2\,dx\\
&&\hs\sim2^{-2k(n+1)}|Q\cap\boz|^{-1} \|a\|^2_{T^2_2
(\boz\times(0,\fz))}|R_k (Q)|\\
&&\hs\ls2^{-k(n+2)}\|a\|^2_{T^2_2 (\boz\times(0,\fz))},
\end{eqnarray*}
which, together with Minkowski's inequality, implies that
\begin{eqnarray}\label{3.55}
\lfz\|\sum_{k\in J_{\boz}}\sum_{i}b_{k,\,i}\r\|_{L^2 (\rn)}&\le&
\sum_{k\in J_{\boz}}\lfz\|\sum_{i}b_{k,\,i}\r\|_{L^2 (\rn)}\\
\nonumber&\ls&\sum_{k\in J_{\boz}}2^{-k(n/2+1)}\|a\|_{T^2_2
(\boz\times(0,\fz))}\\
\nonumber&\ls&\|a\|_{T^2_2 (\boz\times(0,\fz))}.
\end{eqnarray}
Thus, the claim holds and hence
$$\wz{\az}= b_0 +\sum_{k\in
J_{\boz}}\sum_{i}b_{k,\,i}$$
in $L^2 (\rn)$. Furthermore, by
\eqref{3.52}, \eqref{3.54}, the lower type $p_{\Phi}$ property and
$p_{\Phi}\in(n/(n+1),1]$, we have that
\begin{eqnarray}\label{3.56}
\hs\hs&&|Q^{\ast}_0|\Phi\lfz(\frac{\|b_0\|_{L^2
(\rn)}}{|Q^{\ast}_0|^{1/2}}\r)+\sum_{k\in J_{\boz}}\sum_i
|Q^{\ast}_{k,\,i}|\Phi\lfz(\frac{\|b_{k,\,i}\|_{L^2
(\rn)}}{|Q^{\ast}_{k,\,i}|^{1/2}}\r)\\ \nonumber
&&\hs\ls|Q\cap\boz|\Phi\lfz(\frac{\|a\|_{T^2_2
(\boz\times(0,\fz))}}{|Q\cap\boz|^{1/2}}\r)\\
\nonumber &&\hs\hs+\sum_{k\in J_{\boz}}\sum_i
|Q^{\ast}_{k,\,i}|\Phi\lfz(\frac{2^{-k(n+1)}|Q\cap\boz|^{-1/2}
|Q^{\ast}_{k,\,i}|^{1/2}\|a\|_{T^2_2
(\boz\times(0,\fz))}}{|Q^{\ast}_{k,\,i}|^{1/2}}\r)\\
\nonumber &&\hs\ls|Q\cap\boz|\Phi\lfz(\frac{\|a\|_{T^2_2
(\boz\times(0,\fz))}}{|Q\cap\boz|^{1/2}}\r)\\
\nonumber &&\hs\hs+\sum^{\fz}_{k=3}|(2^{(k+1)n}Q)\cap\boz|
\Phi\lfz(\frac{2^{-k(n+1)}\|a\|_{T^2_2
(\boz\times(0,\fz))}}{|Q\cap\boz|^{1/2}}\r)\\
\nonumber &&\hs\ls|Q\cap\boz|\Phi\lfz(\frac{\|a\|_{T^2_2
(\boz\times(0,\fz))}}{|Q\cap\boz|^{1/2}}\r)\lfz\{1+\sum_{k=3}^{\fz}
2^{-[k(n+1)p_{\Phi}-kn]}\r\}\\
\nonumber &&\hs\ls|Q\cap\boz|\Phi\lfz(\frac{\|a\|_{T^2_2
(\boz\times(0,\fz))}}{|Q\cap\boz|^{1/2}}\r),
\end{eqnarray}
which implies that $\wz{\az}\in H_{\Phi}(\rn)$ and \eqref{3.47} in
Case 1).

{\it Case 2)} $8Q\subset\boz$. In this case, let $k_0\in\nn$ such
that $2^{k_0} Q\subset\boz$ but $(2^{k_0
+1}Q)\cap\paz\boz\neq\emptyset$. Then $k_0\ge3$. Let
$$R_k (Q):=(2^{k+1}Q\setminus2^{k}Q)\cap\boz$$
when $k\ge1$ and $R_0 (Q):=2Q$. Let
$$J_{\boz,\,k_0}:=\{k\in\nn:\  k\ge k_0 +1,\,|R_k (Q)|>0\}.$$
For $k\in\zz_+$, let $\chi_k:=\chi_{R_k (Q)}$, $\wz{\chi}_k := |R_k
(Q)|^{-1}\chi_k$,
$$m_k :=\int_{R_k (Q)}\az(x)\,dx,$$
$M_k :=\az\chi_k-m_k \wz{\chi}_k$ and $\wz{M}_k := \az\chi_k$. Then
\begin{equation*}
\az=\sum_{k=0}^{k_0}M_k +\sum_{k\in
J_{\boz,\,k_0}}\wz{M}_k+\sum_{k=0}^{k_0}m_k \wz{\chi}_k.
\end{equation*}
For $k\in\{0,\,\cdots,\,k_0\}$, by the definition of $M_k$, we
know that
$$\int_{\rn}M_k (x)\,dx=0$$
and $\supp M_k\subset2^{k+1} Q$. Moreover, if $k=0$, by H\"older's
inequality and \eqref{3.51}, we have
\begin{eqnarray}\label{3.57}
\|M_0\|_{L^2 (\rn)}&&\ls\|a\|_{T^2_2 (\boz\times(0,\fz))}\\
\nonumber &&
\ls|Q|^{-1/2}[\rz(Q)]^{-1}\ls|2Q|^{-1/2} [\rz(2Q)]^{-1};
\end{eqnarray}
and, if $k\in\{1,\,\cdots,\,k_0\}$, similarly to the proof of
\eqref{3.54}, we have
\begin{equation}\label{3.58}
\|M_k\|_{L^2 (\rn)}\ls \|\az\|_{L^2 (R_k
(Q))}\ls2^{-k(n/2+1)}\|a\|_{T^2_2 (\boz\times(0,\fz))}.
\end{equation}
Thus, for each $i\in\{0,\,\cdots,\,k_0\}$, $M_k$ is a constant
multiple of a $(\rz,\,2,\,0)$-atom with the constant depending on
$k$. Furthermore, from \eqref{3.57}, we deduce that
\begin{equation}\label{3.59}
|2Q|\Phi\lfz(\frac{\|M_0\|_{L^2
(\rn)}}{|2Q|^{1/2}}\r)\ls|Q|\Phi\lfz(\frac{\|a\|_ {T^2_2
(\boz\times(0,\fz))}}{|Q|^{1/2}}\r).
\end{equation}
By \eqref{3.58}, the lower type $p_{\Phi}$ property and
$p_{\Phi}\in(n/(n+1),1]$, we then obtain
\begin{eqnarray}\label{3.60}
\sum_{k=1}^{k_0}|2^{k+1} Q|\Phi\lfz(\frac{\|M_k\|_{L^2 (\rn)}}
{|2^{k+1} Q|^{1/2}}\r)&\ls&\sum_{k=1}^{k_0}|2^k
Q|\Phi\lfz(\frac{\|a\|_{T^2_2
(\boz\times(0,\fz))}} {2^{-k(n+1)}|Q|^{1/2}}\r)\\
\nonumber &\ls&|Q|\Phi\lfz(\frac{\|a\|_ {T^2_2
(\boz\times(0,\fz))}}{|Q|^{1/2}}\r).
\end{eqnarray}
For each $k\in J_{\boz,\,k_0}$, by Fact 1, there exists the Whitney
decomposition $\{Q_{k,\,i}\}_i$ of $R_k (Q)$ about $\paz\boz$ such
that $\cup_i Q_{k,\,i}=R_k (Q)$ and for each $i$, $Q_{k,\,i}$ satisfies
that $2Q_{k,\,i}\subset\boz$ and $4Q_{k,\,i}\cap\paz\boz\neq\emptyset$.
Then $\wz{M}_k =\sum_i \az\chi_{Q_{k,\,i}}$ almost everywhere.
Moreover, by Lemma \ref{l3.9}, for each $k$ and $i$, there exists a
cube $\wz{Q}_{k,\,i}\subset\boz^\complement$ such that
$l(\wz{Q}_{k,\,i})=l(Q_{k,\,i})$ and
$\dist(\wz{Q}_{k,\,i},Q_{k,\,i})\sim l(Q_{k,\,i})$. Then for each
$k$ and $i$, there exists a cube $Q^{\ast}_{k,\,i}$ such that
$(Q_{k,\,i}\cup\wz{Q}_{k,\,i})\subset Q^{\ast}_{k,\,i}$ and
$l(Q^{\ast}_{k,\,i})\sim l(Q_{k,\,i})$. For each $k$ and $i$, let
$$b_{k,\,i}:=\az\chi_{Q_{k,\,i}}-\frac{1}{|\wz{Q}_{k,\,i}|}\lfz\{
\int_{Q_{k,\,i}}\az(x)\,dx\r\}\chi_{\wz{Q}_{k,\,i}}.$$
Then
$$\int_{\rn}b_{k,\,i}(x)\,dx=0$$
and $\supp b_{k,\,i}\subset
Q^{\ast}_{k,\,i}$. Furthermore, similarly to the proof of
\eqref{3.56} and \eqref{3.55}, we obtain that for each $k\in
J_{\boz,\,k_0}$ and $i$, $b_{k,\,i}$ is a constant multiple of a
$(\rz,\,2,\,0)$-atom with the constant, depending on $k$ and $i$, and
\begin{eqnarray}\label{3.61}
\sum_{k\in J_{\boz,\,k_0}}\sum_i
|Q^{\ast}_{k,\,i}|\Phi\lfz(\frac{\|b_{k,\,i}\|_{L^2
(\rn)}}{|Q^{\ast}_{k,\,i}|^{1/2}}\r)\ls|Q|\Phi\lfz(\frac{\|a\|_{T^2_2
(\boz\times(0,\fz))}}{|Q|^{1/2}}\r).
\end{eqnarray}

For $j\in\{0,\,\cdots,\,k_0\}$, let $N_j:= \sum_{k=j}^{k_0}
m_k$. It is easy to see that
\begin{eqnarray}\label{3.62}
\sum_{k=0}^{k_0}m_k \wz{\chi}_k=\sum_{k=1}^{k_0}(\wz{\chi}_k
-\wz{\chi}_{k-1})N_k +N_0 \wz{\chi}_0.
\end{eqnarray}
For any $k\in\{1,\,\cdots,\,k_0\}$, by \eqref{3.53} and
$|\wz{\chi}_{k}-\wz{\chi}_{k-1}|\ls|2^k Q|^{-1}$, we have that
\begin{eqnarray}\label{3.63}
\hs\lfz\|(\wz{\chi}_{k}-\wz{\chi}_{k-1})N_{k}\r\|_{L^2 (\rn)}&\ls&
|2^k Q|^{-1/2}|N_k|\\ \nonumber &\ls&|2^k
Q|^{-1/2}\lfz(\sum_{j=k}^{\fz}2^{-j}\r)
|Q|^{1/2}\|a\|_{T^2_2 (\boz\times(0,\fz))}\\
\nonumber &\ls&2^{-k(n/2+1)} \|a\|_{T^2_2 (\boz\times(0,\fz))}.
\end{eqnarray}
This, together with
$$\int_{\rn}[\wz{\chi}_{k}(x)-\wz{\chi}_{k-1}(x)]\,dx=0$$
and
$$\supp(\wz{\chi}_{k}-\wz{\chi}_{k-1})\subset2^{k}Q,$$
yields that for
each $k\in\{1,\,\cdots,\,k_0\}$,
$(\wz{\chi}_{k}-\wz{\chi}_{k-1})N_k$ is a constant multiple of a
$(\rz,\,2,\,0)$-atom with the constant depending on $k$.
Furthermore, by \eqref{3.63}, the lower type $p_{\Phi}$ property of
$\Phi$ and $p_{\Phi}\in(n/(n+1),1]$, we have
\begin{eqnarray}\label{3.64}
&&\sum_{k=1}^{k_0}|2^k Q|\Phi\lfz(\frac{\|(\wz{\chi}_k
-\wz{\chi}_{k-1})N_k\|_{L^2 (\rn)}}{|2^k Q|^{1/2}}\r)\\
\nonumber&&\hs\ls\sum_{k=1}^{k_0}|2^k Q|\Phi\lfz(\frac{\|a\|_{T^2_2
(\boz\times(0,\fz))}}{2^{-k(n+1)}|Q|^{1/2}}\r)
\ls|Q|\Phi\lfz(\frac{\|a\|_{T^2_2
(\boz\times(0,\fz))}}{|Q|^{1/2}}\r).
\end{eqnarray}

Finally we deal with $N_0\wz{\chi}_0$. By
$$2^{k_0-1}r_0<\dist(x_0,\paz\boz)\le2^{k_0}r_0,$$ 
we know that there
exist a positive integer $M$ and a sequence
$\{Q_{0,\,i}\}_{i=1}^{M}$ of cubes such that
\begin{enumerate}
  \item[(i)] $M\sim2^{k_0};$
  \item[(ii)] for all $i\in\{1,\,\cdots,\,M\}$, $l(Q_{0,\,i})=2r_0$ and
$Q_{0,\,i}\subset\boz$;
  \item[(iii)] for all $i\in\{1,\,\cdots,\,M-1\}$, $Q_{0,\,i}\cap
Q_{0,\,i+1}\neq\emptyset$ and
$$\dist(Q_{0,\,i},\paz\boz)\ge \dist(Q_{0,\,i+1},\paz\boz);$$
  \item[(iv)] $2Q_{0,\,M}\cap\paz\boz\neq\emptyset$.
\end{enumerate}

\noindent Then by Lemma \ref{l3.9}, there exists a cube
$Q_{0,\,M+1}\subset\boz^\complement$ such that $l(Q_{0,\,M+1})=r_0$
and $\dist(Q_{0,\,M},Q_{0,\,M+1})\sim r_0$. Let
$$b_{0,\,1}:=N_0\wz{\chi}_0-\frac{N_0}{|Q_{0,\,1}|}
\chi_{Q_{0,\,1}}$$
and
$$b_{0,\,i}:= \frac{N_0}{|Q_{0,\,i-1}|}\chi_{Q_{0,\,i-1}}-
\frac{N_0}{|Q_{0,\,i}|}\chi_{Q_{0,\,i}}$$
with $i\in\{2,\,\cdots,\,M+1\}$. Obviously, for all
$i\in\{1,\,\cdots,\,M+1\}$, by the definition of $b_{0,\,i}$, we
have that $\int_{\rn}b_{0,\,i}(x)\,dx=0$ and there exists a cube
$Q^{\ast}_{0,\,i}\subset\rn$ such that $\supp b_{0,\,i}\subset
Q^{\ast}_{0,\,i}$ and
\begin{eqnarray}\label{3.65}
l(Q^{\ast}_{0,\,i})\sim l(Q).
\end{eqnarray}

To finish the proof of Proposition \ref{p3.4}(i) in this case, we
need another fact as follows.

\smallskip

{\bf Fact 2.} {\it Let $L$ be as in \eqref{1.3} and satisfy $DBC$.
Let $\boz$, $Q$ and $k_0$ be as the above. Assume that $(G_{\fz})$
holds. For all $x\in\boz$, let
$$\dz(x):=\dist(x,\paz\boz).$$
Then there exist positive constant $C$ and $\bz$, independent on $k_0$
and $Q$, such that for all $x\in Q$,
$$\lfz|\int_{2^{k_0}Q}\paz_t K_t (y,x)\,dy\r|\le\frac{C}{t}
e^{-\frac{\bz[\dz(x)]^2}{t}}.$$}

Now we continue the proof of Proposition \ref{p3.4}(i) by using Fact
2. By Fact 2, \eqref{3.45} and H\"older's inequality, we have that
\begin{eqnarray}\label{3.66}
|N_0|&=&\lfz|\int_{2^{k_0}Q}\az(x)\,dx\r|\\
\nonumber
&=&8\lfz|\int_{2^{k_0}Q} \lfz\{\int_{0}^{\fz}\int_{\boz}D_t
(x,y)a(y,t)\frac{dy\,dt}{t}\r\}dx\r|\\ \nonumber &\le& 8\int_0^{\fz}
\int_{\boz}\lfz|\int_{2^{k_0}Q}D_t
(x,y)\,dx\r||a(y,t)|\frac{dy\,dt}{t}\\ \nonumber &\ls&\|a\|_{T^2_2
(\boz\times(0,\fz))}\lfz\{\int_0^{r_0}\int_{Q}e^{-\frac{2\bz[\dz(y)]^2}{t^2}}
\frac{dy\,dt}{t}\r\}^{1/2}\\ \nonumber &\ls&\|a\|_{T^2_2
(\boz\times(0,\fz))}\lfz\{\int_0^{r_0}\int_{Q}
\lfz(\frac{t}{2^{k_0}r_0}\r)^{2(n+1)/n} \frac{dy\,dt}{t}\r\}^{1/2}\\
\nonumber &\ls&2^{-k_0 (n+1)/n}|Q|^{1/2}\|a\|_{T^2_2
(\boz\times(0,\fz))}.
\end{eqnarray}
For each $i\in\{1,\,\cdots,\,M+1\}$, from the definition of
$b_{0,\,i}$, \eqref{3.65} and \eqref{3.66}, it follows that
\begin{eqnarray}\label{3.67}
\|b_{0,\,i}\|_{L^2 (\rn)}&\ls&|N_0||Q|^{-1/2}\ls2^{-k_0
(n+1)/n}\|a\|_{T^2_2 (\boz\times(0,\fz))}\\ \nonumber &\ls&2^{-k_0
(n+1)/n}|Q|^{-1/2}[\rz(|Q|)]^{-1}\\ \nonumber
&\sim&2^{-k_0
(n+1)/n}|Q_{0,\,i}|^{-1/2}[\rz(|Q_{0,\,i}|)]^{-1},
\end{eqnarray}
which, together with the facts that $\int_{\rn}b_{0,\,i}(x)\,dx=0$
and $\supp b_{0,\,i}\subset Q^{\ast}_{0,\,i}$, implies that
$b_{0,\,i}$ is a constant multiple of a $(\rz,\,2,\,0)$-atom with
the constant depending on $i$. Furthermore, by \eqref{3.67}, the
fact that $M\sim2^{k_0}$ and \eqref{3.65}, we have
\begin{eqnarray}\label{3.68}
&&\sum_{i=1}^{M+1}|Q^{\ast}_{0,\,i}|\Phi\lfz(\frac{\|b_{0,\,i}\|_{L^2
(\rn)}}{|Q_{0,\,i}|^{1/2}}\r)\\
 \nonumber
&&\hs\ls\sum_{i=1}^{M+1}| Q|\Phi\lfz(\frac{\|a\|_{T^2_2
(\boz\times(0,\fz))}}{2^{k_0 (n+1)/n}|Q|^{1/2}}\r)\\
\nonumber&&\hs\ls M2^{-\frac{k_0
(n+1)p_{\Phi}}{n}}|Q|\Phi\lfz(\frac{\|a\|_{T^2_2
(\boz\times(0,\fz))}}{|Q|^{1/2}}\r)\\ \nonumber&&\hs\ls 2^{k_0
[1-(n+1)p_{\Phi}/n]}|Q|\Phi\lfz(\frac{\|a\|_{T^2_2
(\boz\times(0,\fz))}}{|Q|^{1/2}}\r)\\ \nonumber &&\hs\ls
|Q|\Phi\lfz(\frac{\|a\|_{T^2_2 (\boz\times(0,\fz))}}{|Q|^{1/2}}\r).
\end{eqnarray}
Let
$$\wz{\az}:= \sum_{i=1}^{k_0}M_k +\sum_{k\in
J_{\boz,\,k_0}}\sum_i b_{k,\,i}+\sum_{k=1}^{k_0}(\wz{\chi}_k
-\wz{\chi}_{k-1})N_k+\sum_{i=1}^{M+1}b_{0,\,i}.$$
Similarly to the proof of \eqref{3.55}, we know that the
above equality holds in $L^2
(\rn)$. It is easy to see that $\wz{\az}|_{\boz}=\az$.  Furthermore,
from \eqref{3.59}, \eqref{3.60}, \eqref{3.61}, \eqref{3.63} and
\eqref{3.70}, it follows that $\wz{\az}\in H_{\Phi}(\rn)$ and
\eqref{3.47} holds.

To finish the proof of Proposition \ref{p3.4}(i), we need show Fact 2.

Fix $x\in Q$. Choose $\pz_1\in C^{\fz}_{c}(\boz)$ such that
$0\le\pz_1\le1$, $\pz_1\equiv1$ on $Q(x,\frac{\dz(x)}{8})$,
$\supp\pz_1\subset Q(x,\frac{\dz(x)}{4})$, and $|\nabla \pz_1
(z)|\ls\frac{1}{\dz(x)}$ for all $z\in\boz$. Then we have that
\begin{eqnarray}\label{3.69}
\quad&&\lfz|\int_{2^{k_0}Q}\paz_t K_t (y,x)\,dy\r|\\
\nonumber &&\hs\le\lfz|\int_{2^{k_0}Q}\paz_t K_t (y,x)\pz_1
(y)\,dy\r|+\lfz|\int_{2^{k_0}Q}\paz_t K_t (y,x)\lfz[1-\pz_1 (y)\r]\,dy\r|\\
\nonumber &&\hs\le\lfz|\int_{2^{k_0}Q}\paz_t K_t (y,x)\pz_1
(y)\,dy\r|+\int_{2^{k_0}Q\setminus Q(x,\dz (x)/8)}|\paz_t
K(y,x)|\,dy\\
\nonumber &&\hs=:\mathrm{I_1}+\mathrm{I_2}.
\end{eqnarray}
We first estimate $\mathrm{I_1}$. It was proved by Auscher and Russ
in \cite[Proposition A.4]{ar03} that for all $y\in\boz$,
$t\in(0,\fz)$ and all $r\in(0,\fz)$,
\begin{eqnarray}\label{3.70}
&&\lfz\{\int_{\{z\in\boz:\ r\le|y-z|\le2r\}}|\nabla_z K_t
(z,y)|^2\,dz\r\}^{1/2}\\ \nonumber &&\hs\ls
t^{-\frac{1}{2}-\frac{n}{4}}\lfz(\frac{r}{\sqrt{t}}\r)
^{\frac{n-2}{2}}e^{-\gz\frac{r^2}{t}},
\end{eqnarray}
where $\gz$ is a positive constant independent of $y,\,t$ and $r$.
Notice that
$$\paz_t K_t (\cdot,x)+LK_t (\cdot,x)=0$$
and
$Q(x,\frac{\dz(x)}{4})\subset 2^{k_0}Q$. From this, the facts that
$\pz_1\equiv1$ on $Q(x,\frac{\dz(x)}{8})$, $\supp\pz_1\subset
Q(x,\frac{\dz(x)}{4})$, $|\nabla\pz_1 (y)|\ls\frac{1}{\dz(x)}$ for
all $y\in\boz$, H\"older's inequality and \eqref{3.70}, it follows
that
\begin{eqnarray}\label{3.71}
\mathrm{I_1}&=&\lfz|\int_{\boz}L_y K_t (y,x)\pz_1 (y)\,dy\r|\\
\nonumber
&=&\lfz|\int_{\boz}A(y)\nabla_y K_t (y,x)
\cdot\nabla_y\pz_1 (y)\,dy\r|\\
\nonumber &\ls&\int_{\{y\in\boz:\
\frac{\dz(x)}{8}\le|x-y|\le\frac{\dz(x)}{4}\}}|\nabla_y K_t
(y,x)||\nabla_y\pz_1 (y)|\,dy\\ \nonumber &\ls&
\lfz\{\int_{\{y\in\boz:\
\frac{\dz(x)}{8}\le|x-y|\le\frac{\dz(x)}{4}\}}|\nabla_y K_t
(y,x)|^2\,dy\r\}^{1/2}\\ \nonumber &&\times\lfz\{\int_{\{y\in\boz:\
\frac{\dz(x)}{8}\le|x-y|\le\frac{\dz(x)}{4}\}}|\nabla_y\pz_1
(y)|^2\,dy\r\}^{1/2}\\ \nonumber &\ls&
t^{-\frac{1}{2}-\frac{n}{4}}\lfz[\frac{\dz(x)}{\sqrt{t}}\r]
^{\frac{n-2}{2}}e^{-\frac{\gz[\dz(x)]^2}{16t}}[\dz(x)]^{\frac{n-2}{2}}\\
\nonumber &\sim& \frac{1}{t}\lfz[\frac{\dz(x)}{\sqrt{t}}\r]
^{n-2}e^{-\frac{\gz[\dz(x)]^2}{16t}}\ls\frac{1}{t}
e^{-\frac{\gz[\dz(x)]^2}{32t}}.
\end{eqnarray}
For $\mathrm{I_2}$, by Lemma \ref{l2.1}, we have that
\begin{eqnarray*}
\mathrm{I_2}&\ls&\int_{\boz\setminus Q(x,\dz
(x)/8)}\frac{1}{t^{n/2+1}}e^{-\frac{\az|x-y|^2}{t}}\,dy\\&\ls&
\frac{1}{t}
e^{-\frac{\az[\dz(x)]^2}{2^9 t}}\int_{\rn}\frac{1}{t^{n/2}}
e^{-\frac{\az|y|^2}{t}}\,dy\ls\frac{1}{t}
e^{-\frac{\az[\dz(x)]^2}{2^9 t}}.
\end{eqnarray*}
From this, \eqref{3.69} and \eqref{3.71}, it follows that Fact 2
holds, which completes the proof of Proposition \ref{p3.4}(i).

Now we prove (ii) of Proposition \ref{p3.4}. To this end, let $f\in
H_{\Phi,\,S_h,\,d_{\boz}}(\boz)\cap L^2 (\boz)$. Recall that
$d_{\boz}:=2\diam(\boz)$  and we write $(d_{\boz})^2$ simply by
$d^2_{\boz}$. It is easy to see that for all $z\in\cc$ satisfying
$z\neq0$ and $|\arg z|\in(0,\pi/2)$,
$$8\int_0^{d_{\boz}} (t^2ze^{-t^2 z})(t^2ze^{-t^2 z})\frac{dt}{t}
+(2d_{\boz}^2 z+1)e^{-2d_{\boz}^2 z}=1,$$
which, together with the
$H^{\fz}$-functional calculus for $L$,  implies that for all $f\in
L^2 (\boz)$,
\begin{eqnarray}\label{3.72}
f&&=8\int_0^{d_{\boz}} (t^2 Le^{-t^2 L})(t^2 Le^{-t^2
L})(f)\,\frac{dt}{t}\\ \nonumber
&&\hs+\lfz[2d_{\boz}^2Le^{-2d_{\boz}^2 L}(f) +e^{-2d_{\boz}^2
L}(f)\r]=: f_1+f_2.
\end{eqnarray}
We first deal with $f_1$. By the fact that $f\in
H_{\Phi,\,\cn_h}(\boz)\cap L^2 (\boz)$, Propositions \ref{p3.2} and
\ref{p3.3}, we know that $S_h (f)\in L^{\bfai}(\boz)$. From this and
the definition of the space $T_{\bfai}(\boz)$, it follows that
$$t^2
Le^{-t^2 L}(f)\chi_{\boz\times(0,{d_{\boz}})}\in T_{\bfai}(\boz)$$
and
$$\lfz\|f\r\|_{H_{\bfai,\,S_h}(\boz)}=\lfz\|t^2
Le^{-t^2
L}(f)\chi_{\boz\times(0,{d_{\boz}})}\r\|_{T_{\bfai}(\boz)}.$$ Then
by Lemma \ref{l3.8}, there exist $\{\lz_j\}_j \subset\cc$ and a
sequence $\{a_j\}_j$ of $T_{\Phi}(\boz)$-atoms such that for almost
every $(x,t)\in\boz\times(0,\fz)$,
\begin{eqnarray}\label{3.73}
t^2 Le^{-t^2
L}(f)(x)\chi_{\boz\times(0,{d_{\boz}})}(x,t)=\sum_j\lz_j a_j(x,t).
\end{eqnarray}
For each $j$, let
$$\az_j:=8\int_0^{\fz}t^2 Le^{-t^2 L}
(a_j)\frac{dt}{t}.$$
Then by the fact that
$$f_1=8\int_0^{\fz} \lfz(t^2
Le^{-t^2 L}\r)\lfz(t^2 Le^{-t^2
L}(f)\chi_{\boz\times(0,{d_{\boz}})}\r)\,\frac{dt}{t}$$ and
\eqref{3.73}, similarly to the proof of \cite[Proposition
4.2]{jy10}, we have that $f_1=\sum_j\lz_j \az_j$ in $L^2 (\boz)$.
Also, similarly to the proof of \eqref{3.48}, there exists
$\wz{f}_1\in H_{\bfai}(\rn)$ such that $\wz{f}_1|_{\boz}=f_1$ and
$$\|\wz{f}_1\|_{H_{\bfai}(\rn)}\ls\|f\|_{H_{\bfai,\,S_h}(\boz)},$$
which implies that $f_1\in H_{\bfai,\,r}(\boz)$ and
\begin{eqnarray}\label{3.74}
\|f_1\|_{H_{\bfai,\,r}(\boz)}\ls\|f\|_{H_{\bfai,\,S_h}(\boz)}.
\end{eqnarray}

Now we deal with $f_2$. Since $\boz$ is bounded, there exists a
closed cube $\wz{Q}_0\subset\rn$ such that $x_{\wz{Q}_0}\in\boz$,
$l(\wz{Q}_0)\sim {d_{\boz}}$ and $\boz\subset\wz{Q}_0$. Take cubes
$\wz{Q}_1,\,\wz{Q}_2$ such that $\wz{Q}_1\subset\boz^\complement$,
$l(\wz{Q}_1)\sim {d_{\boz}}$,
$(\wz{Q}_0\cup\wz{Q}_1)\subset\wz{Q}_2$ and $l(\wz{Q}_2)\sim
{d_{\boz}}$. Let
$$\wz{f}_2:=f_2\chi_{\wz{Q}_0}-\frac{1}{|\wz{Q}_1|}\lfz[\int_{\boz}f_2
(y)\,dy\r]\chi_{\wz{Q}_1}.$$
Then $\wz{f}_2|_{\boz}=f_2$. It is easy to
see that $\supp \wz{f}_2\subset\wz{Q}_2$, $\int_{\rn}\wz{f}_2
(y)\,dy=0$ and
$$\|\wz{f}_2\|_{L^2 (\rn)}\ls\|f_2\|_{L^2 (\boz)}.$$
Thus, we have that $\wz{f}_2$ is a harmless constant multiple of
some $(\rz,\,2,\,0)$-atom. Denote by $\wz{K}$ the kernel of
$2d_{\boz}^2 Le^{-d_{\boz}^2 L}+e^{-d_{\boz}^2 L}$. Then by Lemma
\ref{l2.1}, we know that for all $x,\,y\in\boz$,
$$|\wz{K}(x,y)|\ls\frac{1}{d_{\boz}^n}
e^{-\frac{\az|x-y|^2}{d_{\boz}^2}},$$
where $\az$ is as in \eqref{2.1}, which, together with
the fact that $\boz$ is
bounded, implies that
\begin{eqnarray*}
\sup_{z\in\boz}|f_2 (z)|&=&\sup_{z\in\boz}\lfz|\int_{\boz}\wz{K}(z,y)
e^{-d_{\boz}^2 L}(f)(y)\,dy\r|\ls\|e^{-d_{\boz}^2 L}(f)\|_{L^1
(\boz)}.
\end{eqnarray*}
From this, the upper type $1$ property of $\bfai$, and the facts
that $\boz\subset\wz{Q}_2$ and $l(\wz{Q}_2)\sim d_{\boz}$, we deduce
that for all $\lz\in(0,\fz)$,
\begin{eqnarray*}
|\wz{Q}_2|\bfai\lfz(\frac{\|\wz{f}_2\|_{L^2 (\rn)}}
{\lz|\wz{Q}_2|^{1/2}}\r)&\ls&|\wz{Q}_2|\bfai\lfz(\frac{\|f_2\|_{L^2
(\boz)}} {\lz|\wz{Q}_2|^{1/2}}\r)\ls
|\wz{Q}_2|\bfai\lfz(\frac{\sup_{z\in\boz}|f_2 (z)|} {\lz}\r)\\
&\ls& |\wz{Q}_2|\bfai\lfz(\frac{\|e^{-d_{\boz}^2 L}(f)\|_{L^1
(\boz)}}{\lz}\r)\\
&\sim&\bfai\lfz(\frac{\|e^{-d_{\boz}^2 L}(f)\|_{L^1 (\boz)}}{\lz}\r).
\end{eqnarray*}
By this, Lemma \ref{l3.1} and the definition of
$H_{\bfai,\,r}(\boz)$, we know that $f_2\in H_{\bfai,\,r}(\boz)$ and
\begin{eqnarray}\label{3.75}
\|f_2\|_{H_{\bfai,\,r}(\boz)}&&\le\|\wz{f}_2\|_{H_{\bfai}(\rn)}\\
\nonumber &&\ls \inf\lfz\{\lz\in(0,\fz):\
\bfai\lfz(\frac{\|e^{-d_{\boz}^2 L}(f)\|_{L^1
(\boz)}}{\lz}\r)\le1\r\}.
\end{eqnarray}
Thus, from \eqref{3.72}, \eqref{3.74} and \eqref{3.75}, it follows
that $f\in H_{\bfai,\,r}(\boz)$ and
$$\|f\|_{H_{\bfai,\,r}(\boz)}\ls\|f\|_{H_{\bfai,\,S_h,\,d_{\boz}}(\boz)},$$
which, together with the arbitrariness of $f\in
H_{\bfai,\,S_h,\,d_{\boz}}(\boz)\cap L^2 (\boz)$, implies that the
first part of Proposition \ref{p3.4}(ii) holds.

We now show the second part of Proposition \ref{p3.3}(ii). We first
prove that
\begin{eqnarray}\label{3.76}
(H_{\bfai,\,S_h,\,d_{\boz}}(\boz)\cap L^2 (\boz))=
(H_{\bfai,\,S_h}(\boz)\cap L^2 (\boz))
\end{eqnarray}
with equivalent norms. Obviously, we have
$$(H_{\bfai,\,S_h,\,d_{\boz}}(\boz)\cap L^2 (\boz))\subset
(H_{\bfai,\,S_h}(\boz)\cap L^2 (\boz))$$
by their definitions. To
show the converse, let $f\in H_{\bfai,\,S_h}(\boz) \cap L^2 (\boz)$.
By Lemma \ref{l3.2}, the contraction property of
$\{e^{-tL}\}_{t\ge0}$ on $L^2 (\boz)$ and H\"older's inequality, we
have that for all $x\in\boz$,
\begin{eqnarray*}
\lfz[S_h (f)(x)\r]^2&\gs&\int_{d_{\boz}/2}^{d_{\boz}}\int_{\boz}\lfz|t^2
Le^{-t^2 L}(f)(y)\r|^2\,\frac{dy\,dt}{t^{n+1}}\gs
d_{\boz}^{-n}\|Le^{-d_{\boz}^2 L}(f)\|_{L^2 (\boz)}^2,
\end{eqnarray*}
which implies that
\begin{equation}\label{3.77}
\inf_{x\in\boz}S_h (f)(x)\gs d_{\boz}^{-n/2}\|Le^{-d_{\boz}^2
L}(f)\|_{L^2 (\boz)}.
\end{equation}

 To continue the proof,
we need the following fact, whose proof is similar to the proof of
\cite[p.\,42, Proposition 5.3]{a07}. We omit the details.

\smallskip

{\bf Fact 3.} {\it Let $1<p<q<\fz$ and
$\az:=\frac{1}{2}(\frac{n}{p}-\frac{n}{q})$. Assume that $(G_{\fz})$
holds. Then $L^{-\az}$ is bounded from $L^p (\boz)$ to $L^q
(\boz)$.}

\smallskip

By $n\ge3$, we know that there exists $p_0\in(1,2]$ and
$q_0\in(1,\fz)$ such that $\frac{1}{p_0}=\frac{2}{n}+\frac{1}{q_0}$.
Then $\frac{1}{2}(\frac{n}{p_0}-\frac{n}{q_0})=1$. By this, Fact 3,
\eqref{3.77} and H\"older's inequality, we obtain that
\begin{eqnarray*}
\lfz\|e^{-d_{\boz}^2 L}(f)\r\|_{L^1 (\boz)}&\ls&\lfz\|e^{-d_{\boz}^2
L}(f)\r\|_{L^{q_0}
(\boz)}\sim\lfz\|L^{-1}Le^{-d_{\boz}^2 L}(f)\r\|_{L^{q_0} (\boz)}\\
\nonumber &\ls&\lfz\|Le^{-d_{\boz}^2 L}(f)\r\|_{L^{p_0} (\boz)}\ls
d_{\boz}^{n/p_0}\inf_{x\in\boz}S_h (f)(x),
\end{eqnarray*}
which, together with the upper type 1 property of $\bfai$,
\eqref{3.77} and \eqref{3.65}, implies that for all $\lz\in(0,\fz)$,
$$\bfai\lfz(\frac{\|e^{-d_{\boz}^2 L}(f)\|_{L^1
(\boz)}}{\lz}\r)\ls\bfai\lfz(\frac{\inf_{x\in\boz}S_h
(f)(x)}{\lz}\r)\ls\int_{\boz}\bfai\lfz(\frac{S_h
(f)(x)}{\lz}\r)\,dx.$$
From this, it follows that $f\in
H_{\bfai,\,S_h,\,d_{\boz}}(\boz)$ and
$$\|f\|_{H_{\bfai,\,S_h,\,d_{\boz}}(\boz)}\ls\|f\|_{H_{\bfai,\,S_h}(\boz)},$$
which, together with the arbitrariness of $f\in
H_{\bfai,\,S_h}(\boz)\cap L^2 (\boz)$, implies that
$$(H_{\bfai,\,S_h}(\boz)\cap L^2 (\boz))\subset
(H_{\bfai,\,S_h,\,d_{\boz}}(\boz)\cap L^2 (\boz)).$$
Thus, \eqref{3.76} holds.

By Propositions \ref{p3.1}, \ref{p3.2}, \ref{p3.3} and \ref{p3.4},
we have
\begin{equation}\label{3.78}
(H_{\bfai,\,\cn_h}(\boz)\cap L^2
(\boz))=(H_{\bfai,\,S_h,\,d_{\boz}}(\boz)\cap L^2 (\boz))
\end{equation}
with equivalent norms. To finish the proof of the second part of
Proposition \ref{p3.4}(ii), let $f\in H_{\bfai,\,\cn_h}(\boz)\cap
L^2 (\boz)$. By
$$e^{-d_{\boz}^2L}(f)=e^{-\frac{d_{\boz}^2}{2}L}
\lfz(e^{-\frac{d_{\boz}^2}{2}L}(f)\r),$$
\eqref{2.1} and the fact that $|\boz|<\fz$,  we know that for all
$x\in\boz$,
\begin{eqnarray*}
\lfz\|e^{-d_{\boz}^2 L}(f)\r\|_{L^1
(\boz)}&\ls&\int_{\boz}\sup_{y\in\boz}\lfz|e^{-\frac{d_{\boz}^2}{2}L}
(f)(y)\r|\,dz\\
&\ls& d_{\boz}^n
\sup_{y\in\boz,\,t\in(0,d_{\boz}),\,|x-y|<t}\lfz|e^{-t^2
L}(f)(y)\r|\sim \cn_h (f)(x).
\end{eqnarray*}
From this, it follows that
$$\|e^{-d_{\boz}^2 L}(f)\|_{L^1
(\boz)}\ls\inf_{x\in\boz}\cn_h (f)(x),$$
which implies that for all $\lz\in(0,\fz)$,
\begin{eqnarray}\label{3.79}
\bfai\lfz(\frac{\|e^{-d_{\boz}^2 L}(f)\|_{L^1
(\boz)}}{\lz}\r)&\ls&\bfai\lfz(\frac{\inf_{x\in\boz}\cn_h
(f)(x)}{\lz}\r)\\ \nonumber &\ls&\frac{1}{|\boz|}\int_{\boz}
\bfai\lfz(\frac{\cn_h (f)(x)}{\lz}\r)\,dx\\ \nonumber
&\sim&\int_{\boz} \bfai\lfz(\frac{\cn_h (f)(x)}{\lz}\r)\,dx.
\end{eqnarray}
By Proposition \ref{3.1} and \eqref{3.79}, we obtain that
$$(H_{\bfai,\,\cn_h}(\boz)\cap L^2
(\boz))=(H_{\bfai,\,\wz{S}_h,\,d_{\boz}}(\boz)\cap L^2 (\boz)),$$
which, together with Proposition \ref{p3.3}, \eqref{3.76},
\eqref{3.78} and the obvious facts that
$$(H_{\bfai,\,\wz{S}_h,\,d_{\boz}}(\boz)\cap L^2 (\boz))\subset
(H_{\bfai,\,\wz{S}_h}(\boz)\cap L^2 (\boz)),$$
implies that
\begin{eqnarray*}
(H_{\bfai,\,\wz{S}_h}(\boz)\cap L^2 (\boz))&&\subset
(H_{\bfai,\,S_h}(\boz)\cap L^2
(\boz))=(H_{\bfai,\,S_h,\,d_{\boz}}(\boz)\cap L^2
(\boz))\\
&&=(H_{\bfai,\,\cn_h}(\boz)\cap L^2 (\boz))
\subset(H_{\bfai,\,\wz{S}_h,\,d_{\boz}}(\boz)\cap L^2
(\boz))\\
&&\subset(H_{\bfai,\,\wz{S}_h}(\boz)\cap L^2 (\boz)).
\end{eqnarray*}
From this, we deduce that
\begin{eqnarray*}
(H_{\bfai,\,\wz{S}_h}(\boz)\cap L^2
(\boz))&&=(H_{\bfai,\,\wz{S}_h,\,d_{\boz}}(\boz)\cap L^2
(\boz))=(H_{\bfai,\,S_h}(\boz)\cap L^2
(\boz))\\
&&=(H_{\bfai,\,\wz{S}_h,\,d_{\boz}}(\boz)\cap L^2 (\boz))
\end{eqnarray*}
with equivalent norms. This finishes the proof of Proposition
\ref{p3.4}(ii) and hence Proposition \ref{p3.4}.
\end{proof}

Combining Propositions \ref{p3.1}, \ref{p3.2}, \ref{p3.3} with
\ref{p3.4}, we then obtain Theorem \ref{t1.1}.

\begin{proof}[\bf{Proof of Theorem \ref{t1.1}.}] We first prove
Theorem \ref{t1.1}(i). By Propositions \ref{p3.1}, \ref{p3.2},
\ref{p3.3} and \ref{p3.4}(i), we know that
\begin{eqnarray*}
(H_{\bfai,\,r}(\boz)\cap L^2 (\boz))&&=(H_{\bfai,\,\cn_h}(\boz)\cap
L^2 (\boz))=(H_{\bfai,\,\wz{S}_h}(\boz)\cap L^2
(\boz))\\
&&=(H_{\bfai,\,S_h}(\boz)\cap L^2 (\boz))
\end{eqnarray*} with
equivalent norms, which, together with the fact that
$H_{\bfai,\,r}(\boz)\cap L^2 (\boz)$, $H_{\bfai,\,\cn_h}(\boz)\cap
L^2 (\boz)$, $H_{\bfai,\,\wz{S}_h}(\boz)\cap L^2 (\boz)$ and
$H_{\bfai,\,S_h}(\boz)\cap L^2 (\boz)$ are, respectively, dense in
$H_{\bfai,\,r}(\boz)$, $H_{\bfai,\,\cn_h}(\boz)$,
$H_{\bfai,\,\wz{S}_h}(\boz)$ and $H_{\bfai,\,S_h}(\boz)$, and a
density argument, implies that the spaces $H_{\Phi,\,r}(\boz)$,
$H_{\Phi,\,\cn_h}(\boz)$, $H_{\bfai,\,\wz{S}_h}(\boz)$ and
$H_{\Phi,\,S_h}(\boz)$ coincide with equivalent norms, which
completes the proof of Theorem \ref{t1.1}(i).

Now we prove  Theorem \ref{t1.1}(ii). By Proposition
\ref{p3.2}, we know that for all 
$f\in H_{\bfai,\,\cn_h}(\boz)\cap L^2 (\boz)$,
$$\|f\|_{H_{\bfai,\,\wz{S}_h}(\boz)}\ls\|f\|_{H_{\bfai,\,\cn_h}(\boz)},$$
which, together with \eqref{3.79}, implies that
$$\|f\|_{H_{\bfai,\,\wz{S}_h,\,d_{\boz}}(\boz)}\ls
\|f\|_{H_{\bfai,\,\cn_h}(\boz)}$$ for all $f\in
H_{\bfai,\,\cn_h}(\boz)\cap L^2 (\boz)$. By the arbitrariness of
$f\in H_{\bfai,\,\cn_h}(\boz)\cap L^2 (\boz)$, we know that
$$(H_{\bfai,\,\cn_h}(\boz)\cap L^2 (\boz))\subset
(H_{\bfai,\,\wz{S}_h,\,d_{\boz}}(\boz)\cap L^2 (\boz)).$$
From this, Propositions \ref{p3.1} and \ref{p3.3} and \ref{p3.4}(ii), it
follows that \begin{eqnarray*} (H_{\bfai,\,r}(\boz)\cap L^2
(\boz))&&=(H_{\bfai,\,\cn_h}(\boz)\cap L^2
(\boz))=(H_{\bfai,\,\wz{S}_h,\,d_{\boz}}(\boz)\cap L^2
(\boz))\\
&&=(H_{\bfai,\,S_h,\,d_{\boz}}(\boz)\cap L^2 (\boz))
\end{eqnarray*}
with equivalent norms, which, together with the fact
that $H_{\bfai,\,r}(\boz)\cap L^2 (\boz)$,
$H_{\bfai,\,\cn_h}(\boz)\cap L^2 (\boz)$,
$(H_{\bfai,\,\wz{S}_h,\,d_{\boz}}(\boz)\cap L^2 (\boz))$ and
$H_{\bfai,\,S_h,\,d_{\boz}}(\boz)\cap L^2 (\boz)$ are, respectively,
dense in $H_{\bfai,\,r}(\boz)$, $H_{\bfai,\,\cn_h}(\boz)$,
$H_{\bfai,\,\wz{S}_h,\,d_{\boz}}(\boz)$ and
$H_{\bfai,\,S_h,\,d_{\boz}}(\boz)$, and a density argument, then
implies that the spaces $H_{\Phi,\,r}(\boz)$,
$H_{\Phi,\,\cn_h}(\boz)$, $H_{\bfai,\,\wz{S}_h,\,d_{\boz}}(\boz)$
and $H_{\Phi,\,S_h,\,d_{\boz}}(\boz)$ coincide with equivalent
norms, which is desired.

Moreover, if $n\ge3$ and $(G_{\fz})$ holds, by the second part of
Proposition \ref{p3.4}(ii) and the fact that
$H_{\bfai,\,\wz{S}_h,\,d_{\boz}}(\boz)\cap L^2 (\boz)$,
$H_{\bfai,\,\wz{S}_h}(\boz)\cap L^2 (\boz)$,
$$H_{\bfai,\,S_h,\,d_{\boz}}(\boz)\cap L^2 (\boz)$$
and $H_{\bfai,\,S_h}(\boz)\cap L^2 (\boz)$ are, respectively, dense in
$H_{\bfai,\,\wz{S}_h,\,d_{\boz}}(\boz)$,
$H_{\bfai,\,\wz{S}_h}(\boz)$, $H_{\bfai,\,S_h,\,d_{\boz}}(\boz)$ and
$H_{\bfai,\,S_h}(\boz)$, together with a density argument,
we obtain that the spaces
$H_{\bfai,\,\wz{S}_h,\,d_{\boz}}(\boz)$,
$H_{\bfai,\,\wz{S}_h}(\boz)$, $H_{\Phi,\,S_h,\,d_{\boz}}(\boz)$ and
$H_{\Phi,\,S_h}(\boz)$ coincide with equivalent norms, which
completes the proof of Theorem \ref{t1.1}(ii) and hence Theorem
\ref{t1.1}.
\end{proof}

Assume that $(G_{\fz})$ holds. For all $t\in(0,\fz)$, let $P_t:=
e^{-t\sqrt L}$. For all $f\in L^2 (\boz)$ and $x\in\boz$, let
$$S_P (f)(x):=\lfz\{\int_{\wz{\bgz}(x)}|t\paz_t P_t
(f)(y)|^2\frac{dy\,dt}{t|Q(x,t)\cap\boz|}\r\}^{1/2}
$$
and
$$\wz{H}^1_{S_P}(\boz):=\lfz\{f\in L^2 (\boz):\
\|f\|_{H^1_{S_P}(\boz)}<\fz\r\},$$
where
$$\|f\|_{H^1_{S_P}(\boz)}:=\|S_P (f)\|_{L^1 (\boz)}.$$
The {\it Hardy space} $H^1_{S_P}(\boz)$ is then defined to be the completion of
$\wz{H}^1_{S_P}(\boz)$ in the norm $\|\cdot\|_{H^1_{S_P}(\boz)}$.

\begin{proposition}\label{p3.5}
Let $\boz$ and $L$ be as in Theorem \ref{t1.1}. Assume that
$(G_{\fz})$ holds. Then $H^1_r (\boz)=H^1_{S_P}(\boz)$ with
equivalent norms.
\end{proposition}

\begin{proof}[\bf{Proof.}]
Similarly to the proof of Proposition \ref{p3.1}, we have that
\begin{equation}\label{3.80}
\lfz(H^1_r (\boz)\cap L^2 (\boz)\r)\subset\lfz(H^1_{S_P}(\boz)\cap
L^2(\boz)\r).
\end{equation}
To finish the proof of Proposition \ref{p3.5}, it suffices to show
\begin{equation}\label{3.81}
\lfz(H^1_{S_P}(\boz)\cap L^2 (\boz)\r)\subset\lfz(H^1_r (\boz)\cap L^2 (\boz)\r).
\end{equation}
Indeed, if \eqref{3.81} holds, by \eqref{3.80}, we have
$$(H^1_r(\boz)\cap L^2 (\boz))
=(H^1_{S_P}(\boz)\cap L^2 (\boz))$$
with equivalent norms, which,
together with the fact that $H^1_r (\boz)\cap L^2 (\boz)$ and
$H^1_{S_P}(\boz)\cap L^2 (\boz)$ are respectively dense in $H^1_r
(\boz)$ and $H^1_{S_P}(\boz)$, and a density argument, implies that
the spaces $H^1_r (\boz)$ and $H^1_{S_P}(\boz)$ coincide with
equivalent norms, which completes the proof of Proposition
\ref{p3.5}.

To show \eqref{3.81}, let $f\in H^1_{S_P}(\boz)\cap L^2 (\boz)$.
Then by the $H^{\fz}$-functional calculus for $L$ and \cite[p.\,164,
(13)]{ar03}, we have that
\begin{equation}\label{3.82}
f=4\int_0^{\fz}t\sqrt{L}P_t (t\sqrt{L}P_t)(f)\frac{dt}{t}
\end{equation}
in $L^2 (\boz)$. By $f\in H^1_{S_P}(\boz)$, we see that $S_P (f)\in
L^1 (\boz)$, which, together with $t\sqrt{L}P_t (f)=t\paz_t P_t
(f)$, implies that $t\sqrt{L}P_t (f)\in T_1 (\boz)$ and
$$\|f\|_{H^1_{S_P}(\boz)}=\|t\sqrt{L}P_t (f)\|_{T_1 (\boz)}.$$
Then from Lemma \ref{l3.8}, it follows that there exist
$\{\lz_j\}_j\subset \cc$ and a sequence $\{a_j\}_j$ of $T_1
(\boz)$-atoms such that
\begin{equation}\label{3.83}
t\sqrt{L} P_t (f)=\sum_j \lz_j a_j.
\end{equation}
For each $j$, let
$$\az_j:=4\int_0^{\fz}t\sqrt{L}P_t(a_j)\frac{dt}{t}.$$
Then by \eqref{3.82} and \eqref{3.83}, similarly
to the proof of \cite[Proposition 4.2]{jy10}, we have
$f=\sum_j\lz_j\az_j$ in $L^2 (\boz)$. For any $T_1 (\boz)$-atom $a$
supported in $\wt{Q\cap\boz}$, let
$$\az:=4\int_0^{\fz}t\sqrt{L}P_t(a)\frac{dt}{t}.$$
To show \eqref{3.81}, similarly to the proof of
Proposition \ref{p3.4}(i), it suffices to show that for $\az$ as the
above, there exist a function $\wz{\az}$ on $\rn$ such that
\begin{equation}\label{3.84}
\wz{\az}|_{\boz}=\az
\end{equation}
and a sequence $\{b_i\}_i$ of harmless constant multiples of
$(1,\,2,\,0)$-atoms, with the constant depending on $i$,  such that
$\wz{\az}=\sum_i b_i$ in $L^2 (\rn)$ and
\begin{equation}\label{3.85}
\sum_i |Q_i|^{1/2}\|b_i\|_{L^2
(\rn)}\ls|Q\cap\boz|^{1/2}\|a\|_{T^2_2 (\boz\times(0,\fz))},
\end{equation}
where for each $i$, $\supp b_i\subset Q_i$.

Let $Q:= Q(x_0,r_0)$. Now we show \eqref{3.84} and \eqref{3.85} by
considering the following two cases for $Q$.

{\it Case 1)} $8Q\cap\boz^\complement\neq\emptyset$. In this case,
the proofs of \eqref{3.84} and \eqref{3.85} are similar to Case 1) of
the proof of Proposition \ref{p3.4}. We omit the details.

{\it Case 2)} $8Q\subset\boz$. In this case, let
$k_0$, $J_{\boz,\,k_0}$, $R_0 (Q)$, $R_k(Q)$,
$\chi_k$, $\wz{\chi}_k$, $m_k$, $M_k$
and $\wz{M}_k$ be as in
Case 2) of the proof of Proposition \ref{p3.4}. Then
$$\az=\sum_{k=0}^{k_0}M_k+\sum_{k\in
J_{\boz,\,k_0}}\wz{M}_k+\sum_{k=0}^{k_0}m_k \wz{\chi}_k.
$$
Similarly to the proof of \eqref{3.59} and \eqref{3.60}, we obtain
that
\begin{equation}\label{3.86}
|2Q|^{1/2}\|M_0\|_{L^2 (\rn)}+\sum_{k=0}^{k_0}|2^k
Q|^{1/2}\|M_k\|_{L^2 (\rn)}\ls|Q|^{1/2}\|a\|_{T^2_2
(\boz\times(0,\fz))}.
\end{equation}
For each $k\in J_{\boz,\,k_0}$, by Fact 1, there exists the Whitney
decomposition $\{Q_{k,\,i}\}_i$ of $R_k (Q)$ about $\paz\boz$ such
that $\cup_i Q_{k,\,i}=R_k (Q)$ and for each $i$, $Q_{k,\,i}$ satisfies
that $2Q_{k,\,i}\subset\boz$ and $4Q_{k,\,i}\cap\paz\boz\neq\emptyset$.
Then $\wz{M}_k =\sum_i \az\chi_{Q_{k,\,i}}$ almost everywhere.
Moreover, by Lemma \ref{l3.9}, for each $k$ and $i$, there exists a
cube $\wz{Q}_{k,\,i}\subset\boz^\complement$ such that
$l(\wz{Q}_{k,\,i})=l(Q_{k,\,i})$ and
$\dist(\wz{Q}_{k,\,i},Q_{k,\,i})\sim l(Q_{k,\,i})$. Then for each
$k$ and $i$, there exists a cube $Q^{\ast}_{k,\,i}$ such that
$(Q_{k,\,i}\cup\wz{Q}_{k,\,i})\subset Q^{\ast}_{k,\,i}$ and
$l(Q^{\ast}_{k,\,i})\sim l(Q_{k,\,i})$. For any $k$ and $i$, let
$$
b_{k,\,i}:=\az\chi_{Q_{k,\,i}}-\frac{1}{|\wz{Q}_{k,\,i}|}\lfz\{
\int_{Q_{k,\,i}}\az(x)\,dx\r\}\chi_{\wz{Q}_{k,\,i}}.$$ 
Then
$$\int_{\rn}b_{k,\,i}(x)\,dx=0$$ 
and $\supp b_{k,\,i}\subset
Q^{\ast}_{k,\,i}$. Furthermore, similarly to \eqref{3.59} and
\eqref{3.56}, we know that for each $k\in J_{\boz,\,k_0}$ and $i$,
$b_{k,\,i}$ is a constant multiple of some $(1,\,2,\,0)$-atom, with
the constant depending on $k$ and $i$, and
\begin{eqnarray}\label{3.87}
\sum_{k\in J_{\boz,\,k_0}}\sum_i
|Q^{\ast}_{k,\,i}|^{1/2}\|b_{k,\,i}\|_{L^2
(\rn)}\ls|Q|^{1/2}\|a\|_{T^2_2 (\boz\times(0,\fz))}.
\end{eqnarray}

For $j\in\{0,\,\cdots,\,k_0\}$, let $N_j:= \sum_{k=j}^{k_0}
m_k$. It is easy to see that
\begin{eqnarray}\label{3.88}
\sum_{k=0}^{k_0}m_k
\wz{\chi}_k=\sum_{k=1}^{k_0}(\wz{\chi}_k-\wz{\chi}_{k-1})N_k +N_0
\wz{\chi}_0.
\end{eqnarray}
Similarly to the proofs of \eqref{3.63} and \eqref{3.66}, we see
that for each $k\in\{1,\,\cdots,\,k_0\}$,
$(\wz{\chi}_k-\wz{\chi}_{k-1})N_k$ is a constant multiple of a
$(1,\,2,\,0)$-atom, with the constant depending on $k$, and
\begin{eqnarray}\label{3.89}
\sum_{k=1}^{k_0}|2^k Q|^{1/2}\|(\wz{\chi}_k
-\wz{\chi}_{k-1})N_k\|_{L^2 (\rn)}\ls|Q|^{1/2}\|a\|_{T^2_2
(\boz\times(0,\fz))}.
\end{eqnarray}

Finally we deal with $N_0\wz{\chi}_0$. Let
$M,\,\{Q_{0,\,i}\}_{i=0}^{M+1}$ and $\{b_{0,\,i}\}_{i=0}^{M+1}$ be
as in Case 2) of the proof of Proposition \ref{p3.4}(i). For all
$t\in(0,\fz)$, we denote the kernel of $P_t$ by $p_t$. Then by Fact
2 and the  subordination formula associated with $L$ that
\begin{equation*}
e^{-t\sqrt{L}}=\frac{1}{\sqrt{\pi}}\int_0^{\fz}
e^{-\frac{t^2}{4u}L}e^{-u}u^{-1/2}\,du
\end{equation*}
(see \cite[p.\,180, (A.1)]{ar03}), we have that for all $x\in Q$,
\begin{equation*}
\lfz|\int_{2^{k_0}Q}\paz_t p_t (y,x)\,dy\r|\ls\frac{1}{t}
\lfz\{1+\frac{\dz(x)}{t}\r\}^{-1},
\end{equation*}
where $\dz(x)$ for $x\in Q$ is as in Fact 2 of
the proof of Proposition \ref{p3.4}.
From this and H\"older's inequality, it follows that
\begin{eqnarray}\label{3.90}
|N_0|&=&\lfz|\int_{2^{k_0}Q}\az(x)\,dx\r|\\
\nonumber &=&4\lfz|\int_{2^{k_0}Q}
\lfz\{\int_{0}^{\fz}\int_{\boz}t\paz_t p_t (x,y)a(y,t)\frac{dy\,dt}{t}\r\}dx\r|\\
\nonumber &\le& 4\int_0^{\fz} \int_{\boz}\lfz|\int_{2^{k_0}Q}t\paz_t
p_t (x,y)\,dx\r||a(y,t)|\frac{dy\,dt}{t}\\ \nonumber
&\ls&\|a\|_{T^2_2
(\boz\times(0,\fz))}\lfz\{\int_0^{r_0}\int_{Q}\lfz[1+\frac{\dz(y)}{t}\r]^{-2}
\frac{dy\,dt}{t}\r\}^{1/2}\\ \nonumber &\ls&\|a\|_{T^2_2
(\boz\times(0,\fz))}\lfz\{\int_0^{r_0}\int_{Q}
\lfz(\frac{t}{2^{k_0}r_0}\r)^2 \frac{dy\,dt}{t}\r\}^{1/2}\\
\nonumber &\ls&2^{-k_0 }|Q|^{1/2}\|a\|_{T^2_2 (\boz\times(0,\fz))}.
\end{eqnarray}
For each $i\in\{1,\,\cdots,\,M+1\}$, by the definition of
$b_{0,\,i}$, \eqref{3.90} and the fact that $l(Q_{0,\,i})\sim l(Q)$,
we have
\begin{eqnarray}\label{3.91}
\|b_{0,\,i}\|_{L^2 (\rn)}&\ls&|N_0||Q|^{-1/2}\ls2^{-k_0
}\|a\|_{T^2_2 (\boz\times(0,\fz))}\\ \nonumber
&\ls&2^{-k_0}|Q|^{-1/2}[\rz(|Q|)]^{-1}\\
\nonumber&\sim&2^{-k_0}|Q_{0,\,i}|^{-1/2} [\rz(|Q_{0,\,i}|)]^{-1},
\end{eqnarray}
which, together with the facts that
$$\int_{\rn}b_{0,\,i}(x)\,dx=0$$
and $\supp b_{0,\,i}\subset Q^{\ast}_{0,\,i}$, implies that
$b_{0,\,i}$ is a constant multiple of some $(1,\,2,\,0)$-atom with
the constant depending on $i$. Furthermore, by \eqref{3.91} and the
fact that $M\sim2^{k_0}$, we obtain
\begin{eqnarray}\label{3.92}
\hs\hs\sum_{i=1}^{M+1}|Q^{\ast}_{0,\,i}|^{1/2}\|b_{0,\,i}\|_{L^2
(\rn)}&\sim&\sum_{i=1}^{M+1}2^{-k_0}|Q|^{1/2}\|a\|_{T^2_2
(\boz\times(0,\fz))}\\ \nonumber
&\sim&|Q|^{1/2}\|a\|_{T^2_2
(\boz\times(0,\fz))}.
\end{eqnarray}
Let
$$\wz{\az}:= \sum_{i=1}^{k_0}M_k +\sum_{k\in
J_{\boz,\,k_0}}\sum_i b_{k,\,i}+\sum_{k=1}^{k_0}(\wz{\chi}_k
-\wz{\chi}_{k-1})N_k+\sum_{i=1}^{M+1}b_{0,\,i}.$$
By an argument
similar to that used in the estimate \eqref{3.55}, we know that the
series in the definition of $\wz{\az}$ converges in $L^2 (\rn)$. It
is easy to see that $\wz{\az}|_{\boz}=\az$. Furthermore, by
\eqref{3.86}, \eqref{3.87}, \eqref{3.89} and \eqref{3.92}, we have
that $\wz{\az}\in H_{\Phi}(\rn)$ and \eqref{3.85} holds, which
completes the proof of Proposition \ref{p3.5}.
\end{proof}

From Proposition \ref{p3.5}, we deduce that for any given $f\in
H^1_{S_P}(\boz)$, there exists an atomic decomposition, which gives
a positive answer to the question asked by Duong and Yan
\cite[p.\,485, Remarks (iii)]{dy04} in the case that $p=1$. However,
it is still unknown whether this method also works for $p<1$ but
near to $1$, which seems to need nicer estimate than \eqref{3.90}.

\medskip

{\bf Acknowledgements.} Both authors would like to thank the referee
for her/his careful reading and several valuable remarks which
improve the presentation of this article.

\vspace{1cm}

\hfill Dachun  Yang

\hfill School of Mathematical Sciences

\hfill Beijing Normal University

\hfill Laboratory of Mathematics and Complex Systems

\hfill Ministry of Education,

\hfill Beijing 100875

\hfill People's Republic of China

\hfill \texttt{dcyang@bnu.edu.cn}

\bigskip

\hfill Sibei  Yang

\hfill School of Mathematical Sciences

\hfill Beijing Normal University

\hfill Laboratory of Mathematics and Complex Systems

\hfill Ministry of Education,

\hfill Beijing 100875

\hfill People's Republic of China

\hfill \texttt{yangsibei@mail.bnu.edu.cn}
\footnotetext{\hspace{-0.89cm}
The first author is supported by the
National Natural Science Foundation (Grant No. 10871025) of China
and Program for Changjiang Scholars and Innovative Research Team in
University of China.}

\end{document}